\documentclass{siamart}

\usepackage{lipsum}
\usepackage{amsfonts,amssymb}
\usepackage{graphicx}
\usepackage{algorithmic}
\usepackage{bm,bbm}
\usepackage{pstricks}
\usepackage{tikz}
\usepackage{pgfplots}
\pgfplotsset{width=10cm,compat=newest}

\usetikzlibrary{arrows}

\usetikzlibrary{decorations.markings}
\tikzset{->-/.style={decoration={
  markings,
  mark=at position #1 with {\arrow{>}}},postaction={decorate}}}

\tikzset{middlearrow/.style={
        decoration={markings,
            mark= at position 0.5 with {\arrow{#1}} ,
        },
        postaction={decorate}
    }
}

\DeclareGraphicsExtensions{.pdf,.png,.jpg}
\newcommand{\TheTitle}{A primer on noise-induced transitions in applied dynamical systems}
\newcommand{\TheAuthors}{E.\ Forgoston and R.O.\ Moore}

\newcommand{\cL}{{\mathcal L}}

\newcommand{\bE}{{\mathbb E}}
\newcommand{\ti}{t_{\mathrm i}}
\newcommand{\tf}{t_{\mathrm f}}
\newcommand{\xinit}{x_{\mathrm i}}
\newcommand{\xf}{x_{\mathrm f}}

\def\csch{\mathop{\rm csch}\nolimits}

\newcommand{\partialderiv}[3][]{\frac{\partial^{#1}#2}{\partial {#3}^{#1}}}

\def\_#1{{\underline{#1}}}
\def\Re{\mathop{\rm Re}}

\newcommand{\bp}{\begin{pmatrix}}
\newcommand{\ep}{\end{pmatrix}}

\newcommand{\be}{\begin{equation}}
\newcommand{\ee}{\end{equation}}

\newcommand{\sech} { {\rm sech} \hskip 0.01in}

\definecolor{myblue}{rgb}{.8, .8, 1}

\newrgbcolor{myeqn}{1.0 0.0 0.0}
\newrgbcolor{orange}{1.00 0.65 0.00}
\newrgbcolor{mygreen}{0.1 0.5 0.0}
\newrgbcolor{firebrick}{.7 .13 .13}
\newrgbcolor{mycyan}{0 .75 .75}
\newrgbcolor{myeqn}{1.0 0.0 0.0}
\newrgbcolor{orange}{1.00 0.65 0.00}
\newrgbcolor{mygreen}{0.1 0.5 0.0}
\newrgbcolor{firebrick}{.7 .13 .13}
\newrgbcolor{mycyan}{0 .75 .75}

\def\Expect#1{\mathbb{E}\hspace{-0.03in}\left[ #1 \right]}

\def\Var#1{\mathbb{V}\hspace{-0.03in}\left[ #1 \right]}
\def\CV#1{\mathbb{C}_{\mathrm{var}}\hspace{-0.03in}\left[ #1 \right]}

\newcommand\Des[2]{#1_{\mathrm{#2}}}

\headers{\TheTitle}{\TheAuthors}

\title{{\TheTitle}}

\author{Eric Forgoston\thanks{Department of Mathematical Sciences, Montclair
    State University, Montclair, NJ 07043,
    USA (\email{eric.forgoston@montclair.edu}). }
    \and Richard O.\ Moore\thanks{Department of Mathematical
    Sciences, New Jersey Institute of Technology, Newark, NJ 07102, USA (\email{rmoore@njit.edu}).} }

\usepackage{amsopn}

\hypersetup{
  pdftitle={\TheTitle},
  pdfauthor={\TheAuthors}
}

\begin{document}
\maketitle

\begin{abstract}
 Noise plays a fundamental role in a wide variety of physical and biological
dynamical systems. It can arise from an external forcing or due to random dynamics internal to the system.
It is well established that even weak noise can result in large behavioral changes such as
transitions between or escapes from quasi-stable states.  These transitions can correspond to critical events such as failures or extinctions that make them essential phenomena to understand and quantify, despite the fact that their occurrence is rare.
This
article will provide an overview of the theory underlying the dynamics of rare
events for stochastic models along with some example applications.
\end{abstract}



\section{Introduction}
\label{s:intro}

Recent years have seen a dramatic rise in the use of stochastic systems to
model a wide variety of important physical and biological phenomena.
Studies of sub-cellular processes and tissue dynamics~\cite{Tsimring_2014}, large-scale population
dynamics~\cite{ovamee2010}, genetic switching~\cite{Assaf2011},
magnetic devices~\cite{kohn2005}, optical devices\cite{jones_carrier-envelope_2000,shtaif_m._carrier_2011,paschotta_optical_2005,desurvire_erbium-doped_2002}, Josephson
junctions~\cite{Fulton1974}, fluid mechanics\cite{forsch09,fbys11,heckman2014}, weather and climate\cite{vanden-eijnden_data_2012}, and geosciences\cite{smith_uncertainty_2013} have included numerous investigations into how noise
affects physical and biological phenomena at a wide variety of scales.
One often sees rare transition events in these systems that are induced by noise which
may be internal or  external to the system. These noise-induced rare events
may be associated with a desirable outcome, such as the extinction of an
infectious disease outbreak~\cite{Assaf2010,bilfor17} or eradication of a pest~\cite{Nieddu2014}, or an undesirable outcome, such as the sudden
clustering of cancerous cells~\cite{Khain2014}, the outbreak of an infectious
disease~\cite{nieddu2017}, or a bit error in a digital communication
system~\cite{falkovich_statistics_2001}. Another important class of systems
are those that exhibit rare events due to intermittently activated finite-time
instabilities~\cite{mohamad2016a,farazmand2017}. These types of rare events include freak water waves in
the ocean~\cite{cousins2016} and ship rolling and capsizing~\cite{mohamad2016}.
The need to understand these phenomena has fueled a demand for methods that
quantify the impact of noise on complex systems.

In
these stochastic systems, noise can affect the system in a variety of
ways. Assessing the full impact of noise is rarely possible due to the ``curse
of dimensionality''.  Analysis and computations must often concentrate on the most
important noise-induced events, which include spontaneous switching between
coexisting stable states, critical failure of a system, and nucleation of
coherent structures.
One important feature of interest when studying
noise-induced transitions is the optimal transition pathway of escape from a
metastable state either to another metastable state or to a stable (absorbing) state.
We define the optimal path as the path that is most
likely to occur among all possible paths, recognizing that this path may not be unique or even exist (i.e., belong to the set of admissible paths).
For systems out of equilibrium such that detailed balance is not satisfied, knowledge of the optimal path is required for the computation of
the mean switching time between states or the mean time to exit.

This article provides a tutorial overview of noise-induced rare events. By providing theoretical and numerical background along with several
example applications of recently developed methods to compute the most likely
noise pathways to critical events, we hope to enable the reader to quickly apply
the methods to original applications of interest.

\subsection{Noise sources}
\label{s:noise}
Stochasticity manifests itself as either external or internal noise. External
noise comes from a source outside the system being considered (e.g. population
growth under the influence of climatic effects, or a random signal fed into a
transmission line), and often is modeled by replacing an
external parameter with a random process.
Internal noise is inherent to the system itself and is caused by the random
interactions of discrete particles
(e.g. individuals in a population, or molecules undergoing a chemical
reaction)~\cite{hanggi1990reaction,gar03,vanKampen_book,touchette2009}. Both types of noise can lead to a
rare switching event between metastable states or a rare escape event from a
metastable state.

\begin{figure}[t!]
\begin{minipage}{0.49\linewidth}
\begin{center}
\includegraphics[scale=0.41]{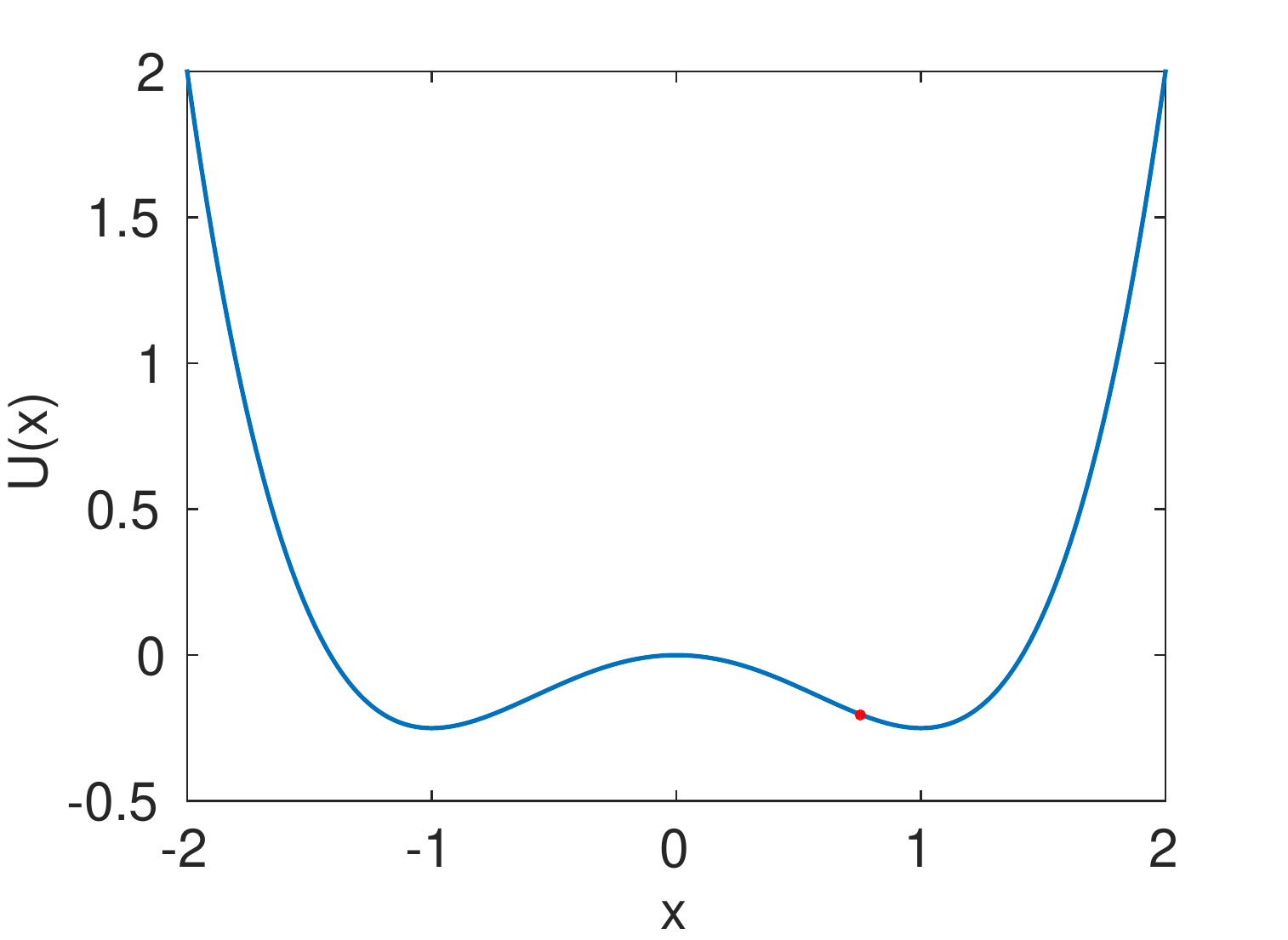}
\end{center}
\end{minipage}
\begin{minipage}{0.49\linewidth}
\begin{center}
\includegraphics[scale=0.41]{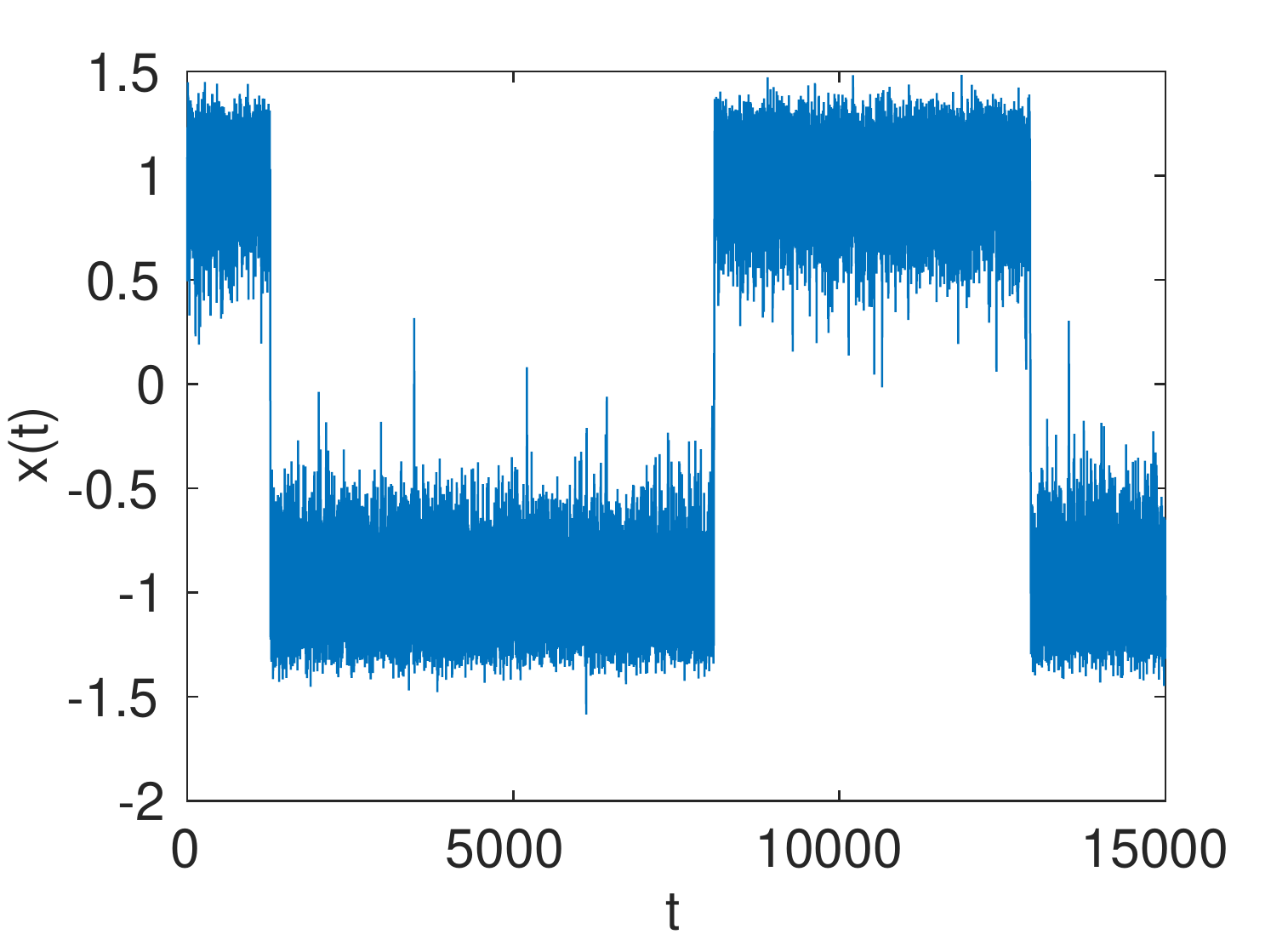}
\end{center}
\end{minipage}
\caption{\label{fig:ext-switch}(a)Particle moving in a potential
well under the influence of external noise. (b)Stochastic realization of the particle position in time.}
\end{figure}

\begin{figure}[t!]
\begin{center}
\includegraphics[scale=0.5]{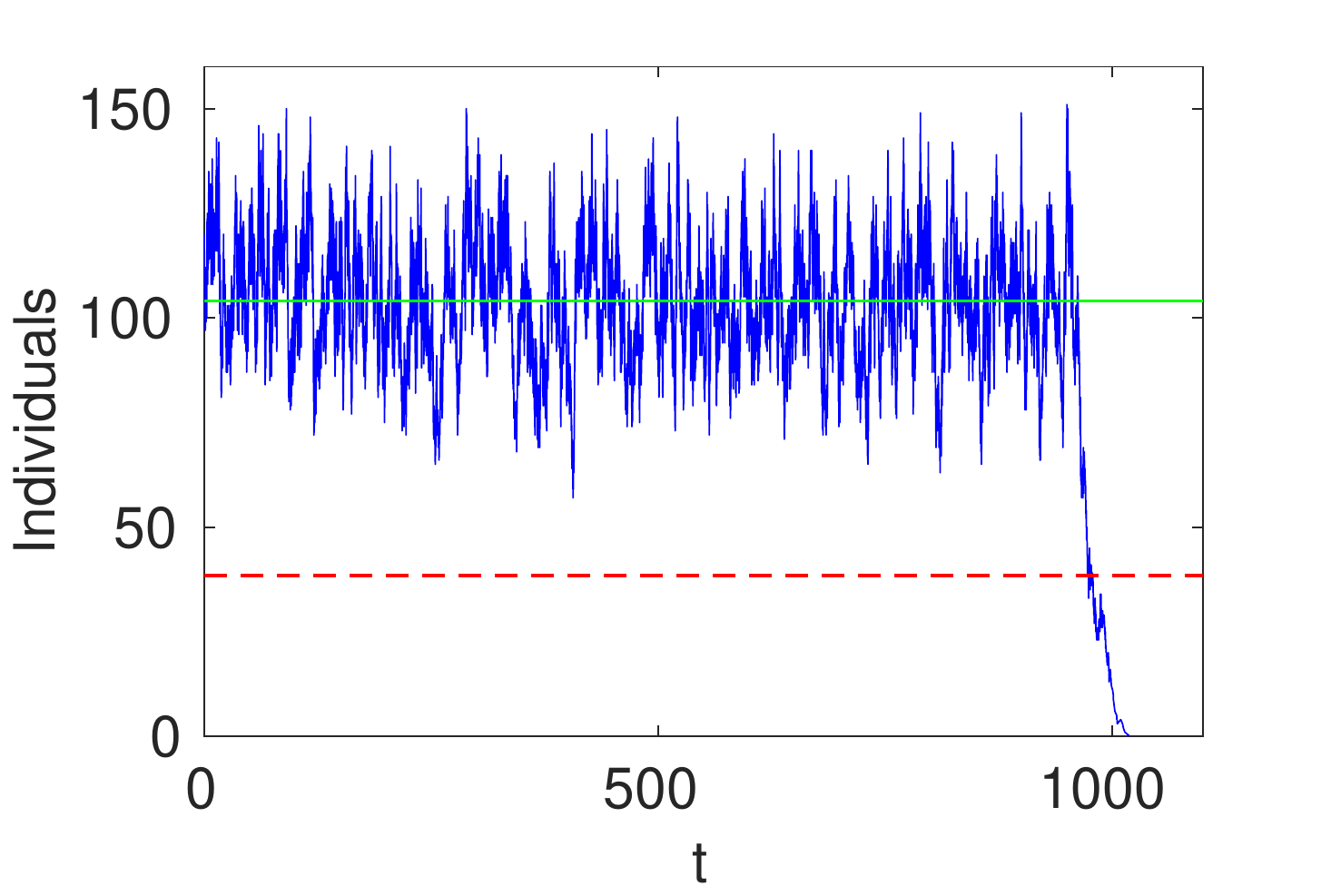}
\end{center}
\caption{\label{fig:int-switch}Stochastic realization of the population size
  in time where the local dynamics of the population exhibit the Allee
  effect.  The green, solid line denotes the carrying capacity, while the red, dashed line denotes the Allee threshold.}
\end{figure}

Figure~\ref{fig:ext-switch}(a) shows a snapshot in time of
a particle moving in a quartic potential
well under the influence of external noise.  Figure~\ref{fig:ext-switch}(b)
shows a single realization of the particle position in time. One sees multiple rare
stochastic switches from one well to another.
Figure~\ref{fig:int-switch} shows a single realization of the
population size in time for an Allee effect~\cite{Allee1931} problem in population biology. In
this example, the noise is internal and arises from the
random interactions of individuals in a population. One sees the population fluctuating about the carrying
capacity for a long period of time until the rare extinction event occurs.

There are many possible escape/switching paths, but there is
a path along which switching or escape is most likely to occur. We call this
most likely path of escape or switching the optimal path. It is of great
importance in a variety of applied problems to determine this optimal path
since knowledge of the path then enables the determination of the mean time to
escape from a metastable state or to switch from one metastable state to
another metastable state.

Mathematically, the effect of external noise is often described using a Langevin
equation or the associated Fokker-Planck equation (though the dynamics of
external noise may sometimes be described by a master equation~\cite{Roberts2015}). Feynman famously pointed out that each
noise realization corresponds to a particular trajectory of the system, and
therefore the probability density of realizations of trajectories is
determined by the probability density of noise realizations~\cite{feyhib65}. This idea can be
used to formulate a variational problem to find the optimal path that ultimately reduces to considering trajectories of an auxiliary Hamiltonian dynamical system.
One can solve for the Hamiltonian dynamics, either analytically or numerically, for
the most probable (i.e., optimal) path of escape or switching~\cite{FW84,Dykman1994d}.

The effect of internal
noise due to the random interactions of individuals within the system is
described mathematically using a master equation.
The master equation is a
large, or even infinite, set of differential equations that cannot in general be solved analytically. However, the optimal fluctuation path is found
to solve an eikonal equation obtained via Wentzel-Kramers-Brillouin-Jeffreys (WKBJ) analysis of the master equation in the limit of weak internal noise~\cite{Kubo1973,gang87,Dykman1994d,elgkam04,Kessler07,Assaf2010,fbss11,sfbs11,Nieddu2014,bilfor17,nieddu2017}.
As in the case of external noise, the optimal path is found by
transforming the original stochastic problem into a new deterministic system
described by a Hamiltonian $\mathcal{H}(x,\lambda)$.  The dimension of the Hamiltonian is
twice the dimension of the original system due to
the conjugate momentum variables $\lambda$.
An optimal path starting at a metastable state is identified with $\mathcal{H}(x,\lambda)=0$; therefore, the method amounts to finding a zero-energy trajectory
of an effective mechanical system, and at least one of the solutions to
the zero-energy Hamiltonian is the optimal path. There
may be other escape, switching, or extinction paths, but the optimal path is
the path that maximizes the probability of escape, switching, or extinction~\cite{bauver16}.

\subsection{Infinite-dimensional equations and large deviation theory}
\label{s:infiniteDim}

Even though it is often straightforward to extend the framework discussed in this
article to infinite-dimensional systems with states defined on Hilbert spaces~\cite{da_prato_stochastic_2008}, we generally restrict our discussion to finite-dimensional examples.  The exception to this is the formal introduction of stochastic PDEs used to model mode-locked lasers and communication systems introduced in Sections~\ref{s:phaseSlips} and~\ref{s:IS}.  The simulations performed in these sections, particularly those using importance sampling applied to~(\ref{e:filter}), effectively limited the bandwidth of noise that was otherwise taken to be delta-correlated in both time and space.  To the authors' knowledge no proof of well-posedness exists for this equation, and indeed when $c_2=0$ in~(\ref{e:filter}), there is evidence to suggest that the stochastic PDE is ill-posed~\cite{cargill_incorporating_2016,cargill_noise_2011}.

Finally, we note that large deviation results for stochastic PDEs are becoming more common in the literature, particularly in the context of parabolic equations with damping at high wavenumbers~\cite{hairer_large_2015,hairer_central_2017}  Large deviation principles for the nonlinear Schr\"odinger equation with additive and multiplicative noise have been shown when the noise is colored in space~\cite{de_bouard_stochastic_2003,gautier_large_2005,gautier_uniform_2005}.

\section{Transitions, exits, and extinctions}
\label{s:exitProblems}

\subsection{Mean first exit time due to external noise}
\label{s:background}

The effect of noise on a physical system can be described in a variety of
ways, all of which aim to capture how a single deterministic solution has been
transformed by noise into a distribution of solutions and what the
characteristics are of that distribution.  This can take the form of a
boundary value problem with random data or random coefficients, as one might
find in the forward (scattering) problem associated with imaging of the earth
or of living tissue, or it can take the form of an initial value problem with
random data, as one might find in numerical weather prediction.  We will focus
here on initial value problems with time-dependent random forcing,
specifically on stochastic dynamical systems described by a well-defined state
evolving in time, where the evolution is probabilistic due to the presence of
noise.  This section provides some definitions to set the context of what
follows. Although we will consider a one-dimensional (1D) framework, one can
naturally extend the results to the multi-dimensional case.

We begin by considering the stochastic differential equation (SDE) described by
\be
\dot{x} = f(x) + \eta(t),
\label{e:SDE}
\ee
where $\eta(t)$ expresses an additive randomly distributed noise term influencing the state $x(t)$.
Whether resulting from truly random (e.g., quantum) phenomena or simply an assortment of physical phenomena not included in the model for reasons of scale or uncertainty, $\eta(t)$ can often be equated with the formal time derivative of a Brownian motion, resulting in paths $x(t)$ that are continuous and Markovian (i.e., ``memoryless'').  Since these paths, referred to as diffusion processes, are generally not differentiable, it is helpful to consider the statistics of finite increments $\Delta\eta = \eta(t+\Delta t)-\eta(t)$.  These increments are drawn from a Gaussian distribution with mean and variance given by
\be
\bE[\Delta\eta] = 0\quad\mbox{and}\quad
\bE[(\Delta \eta)^2] = D\Delta t,
\ee
where $\bE[\cdot]$ denotes the expectation or expected value of a
quantity with respect to its distribution.
A natural way to think of~(\ref{e:SDE}) is therefore as the limit of the Euler-Maruyama method~\cite{oksendal_stochastic_2013}, which is simply the expression of the comments above in the form of a numerical integration method to obtain $\Delta x = x(t+\Delta t)-x(t)$:
\be
\Delta x = f(x(t))\Delta t + \Delta\eta.
\ee
Then~(\ref{e:SDE}) is understood as the limiting process as $\Delta t\rightarrow 0$, such that
\be
\bE[\eta(t)] = 0\quad\mbox{and}\quad \bE[\eta(t)\eta(t')] = D\delta(t-t'),
\label{e:noiseProps}
\ee
where $D$ is the noise intensity.  Details of the Euler-Maruyama method can be found in Appendix~\ref{s:appA}.

Some of the SDEs considered in the following sections have state-dependent functions multiplying the random term, i.e.,
\be
\dot{x} = f(x) + \sigma(x)\eta(t).
\label{e:SDECov}
\ee
This expression is subject to interpretation depending on how one determines the limit described above.
Two standard interpretations of this stochastic integral predominate, including the prepoint or It\^o interpretation,
\be
\int_t^{t+\Delta t}\sigma(x(s))\eta(s)\,ds \approx \sigma(x(t))\Delta\eta,
\ee
and the midpoint or Stratonovich interpretation,
\be
\int_t^{t+\Delta t}\sigma(x(s))\eta(s)\,ds \approx \frac12\left[\sigma(x(t))+\sigma(x(t+\Delta t))\right]\Delta\eta.
\ee
The examples to follow will use only the It\^o interpretation for simplicity, although for many physical models the Stratonovich interpretation is more appropriate~\cite{vanKampen_book}.

Reflecting the fact that~(\ref{e:SDE}) is best understood as a distribution of
paths determined by possible realizations of the random driving term
$\eta(t)$, its solution can be alternatively expressed as a transition density
$p(x,t|y,s)$, such that the probability of finding $x(t)$ in a measurable set $\Omega$ conditioned on $x(s) = y$ with $s\leq t$ is given by
\be
P(x(t)\in \Omega|x(s)=y) = \int_{\Omega}p(x,t|y,s)\,dx.
\ee
The time evolution of this transition density satisfies the forward Fokker-Planck equation,
\be
\partialderiv{p}{t} = -\partialderiv{}{x}(fp) + \frac12D\partialderiv[2]{p}{x}
\label{e:fFPE}
\ee
with initial condition $p(x,s|y,s) = \delta(x-y)$.
Conditional expectations can then be computed against this density using
\be
\overline{g}(y,t)=\bE[g(x(t))|y,s] = \int g(x)p(x,t|y,s)\,dx,
\ee
or the quantity $\overline{g}(y,t)$ can be evolved directly using the backward Fokker-Planck equation,
\be
\partialderiv{\overline{g}}{t} = f(y)\partialderiv{\overline{g}}{y} + \frac12D\partialderiv[2]{\overline{g}}{y}
\label{e:bFP}
\ee
with initial condition $\overline{g}(y,s) = g(y)$.

An important application of the Fokker-Planck equation is in computing exit times of state $x$ from a domain $\Omega$.  If we define $g$ to be the indicator function associated with set $\Omega$, $g(x) = {\mathbbm 1}_{x\in \Omega}$, then $\overline{g}(y,t)$ is the probability of being in $\Omega$ at time $t$ having started at $x(0)=y$ (we set $s=0$ without loss of generality for time-homogenous processes).  If we focus only on the first exit from $\Omega$, we note that this probability satisfies~(\ref{e:bFP}) with absorbing boundary condition
$\overline{g}\vert_{\partial \Omega} = 0$~\cite{gar03,pavliotis_stochastic_2011}.  Now we note that over a time interval $(t,t+\Delta t)$, the probability of exit from a trajectory starting at $y$ is given by
\be
\overline{g}(y,t) - \overline{g}(y,t+\Delta t) \approx -\partialderiv{\overline{g}}{t}\Delta t.
\ee
In other words, $-\partial\overline{g}/\partial t$ is the exit time density.  Under the assumption that all trajectories exit with probability one, the density can be used to compute the mean first exit time $\tau(y)$ for any $y\in \Omega$:
\be
\tau(y) = \int_0^\infty t(-\partialderiv{\overline{g}}{t})\,dt = \int_0^\infty \overline{g}(y,t)\,dt,
\ee
where we have integrated by parts and assumed that $t\overline{g}\rightarrow 0$ as $t\rightarrow\infty$.  Integrating~(\ref{e:bFP}) thus yields a boundary value problem referred to as Dynkin's Equation for the mean first exit time~\cite{gar03},
\be
f(y)\partialderiv{\tau}{y} + \frac12D\partialderiv[2]{\tau}{y} = -1,\qquad
\tau\left\vert\right ._{\delta \Omega} = 0.
\label{e:Dynkin}
\ee

\subsection{Large deviation theory, action and optimal paths}
\label{s:LDT}

 Section~\ref{s:background} poses two broad alternatives for computing
 expectations of functions of a diffusion process.  The first  is to
 solve the Fokker-Planck equation~(\ref{e:fFPE}) using, for example, eigenfunction expansions
 (with an inner product
 defined with respect to a weighting that renders the operator self-adjoint).  This
 approach is viable in one dimension as represented in~(\ref{e:fFPE}), but quickly becomes numerically intractable if the state is high- or infinite-dimensional.  This is partly due to the fact that, in the limit of small noise where
 $D\rightarrow 0$, the differential operator in~(\ref{e:fFPE}) is
 singular. This feature implies the existence of boundary layers, rendering the computation numerically stiff.  Analytically, the existence of a small parameter suggests the use of asymptotic methods such as singular perturbation or WKBJ theory to approximate the expectation~\cite{bobrovsky_singular_1982,matkowsky_exit_1977}.

The second broad approach is to simulate the random walk in~(\ref{e:SDE})
multiple times and to compute the expectations empirically using, for example,
a Monte Carlo method.  The difficulty in using this approach for rare events in high-dimensional systems stems primarily from the fact that the desired
expectation is dominated by a small set of events of interest. This problem is
exacerbated in the small-noise limit due to an exponential drop-off in
likelihood as the random walker moves away from this set.  In this
scenario of high dimension and small noise, it is thus critically important to
bias the Monte Carlo simulations such that a substantial number of simulated
trajectories coincide with those that sample from the set of interest with
high likelihood.
This search for a good bias function is made all the more
important by the inherent inefficiency of the Monte Carlo method, which
produces estimators with relative error that scales as $1/\sqrt{N}$ where $N$
is the sample size of simulations.  The importance of this issue has led to the development of many techniques appropriate for rare-event sampling in a variety of physical contexts~\cite{hartmann_efficient_2012,vanden-eijnden_rare_2012,zhang_applications_2014,adams_barrier_2010}.

Regardless of whether one uses asymptotic methods to approximate an
expectation or uses importance sampling to concentrate Monte Carlo
simulations on states of interest, both approaches lead to the theory of large
deviations described in detail in the text by Freidlin and Wentzell~\cite{FW84}.  We will now
present a  description of this theory, focusing on providing a practical
introduction to exit problems on finite and infinite time
horizons.

We return to~(\ref{e:SDE}),
\begin{equation}
\dot{x}=f(x) + \eta (t),\label{e:SDE2}
\end{equation}
where the noise $\eta (t)$ is as described in~(\ref{e:noiseProps}).
Each noise realization $\eta (t)$ produces an associated trajectory $x(t)$. Therefore, the probability of
generating a given path $x(t)$ is proportional to that of generating the corresponding noise realization $\eta(t)$~\cite{feyhib65}.
The formal probability density for this Gaussian noise process satisfies
\begin{equation}
P[\eta (t)]\propto\exp{\left (-\frac{1}{2D}\int_{t_i}^{t_f}\eta ^2\,dt\right )}.\label{eq:prob}
\end{equation}
Substitution of~(\ref{e:SDE2}) into~(\ref{eq:prob}) gives~\cite{feyhib65,dykkri79}
\begin{equation}
P[x(t)]\propto\exp{\left ( -\frac{1}{2D}\int_{t_i}^{t_f} (\dot{x}-f(x))^2 \, dt\right )}.\label{eq:prob2}
\end{equation}

Since we are concerned with maximizing this probability density over a set of eligible paths satisfying, for example, boundary conditions $x(t_i)=x_i$ and $x(t_f)=x_f$, we express~(\ref{eq:prob2}) using the conventions of variational calculus, where
\begin{equation}
P[x(t)]\propto\exp{\left ( -\frac{1}{2D}\int [\dot{x}-f(x)]^2 \, dt\right )}\propto\exp{\left ( -\mathcal{S}[x(t)]/D\right )},
\end{equation}
with action ${\mathcal S}$ and Langrangian ${\mathcal L}$~\cite{frewen70,dykkri79,FW84,lumcdy98} defined by
\begin{equation}
\mathcal{S}[x(t)]= \int\limits_{t_i}^{t_f}\mathcal{L}(\dot{x},x;t)\,
  dt, \qquad \mathcal{L}(\dot{x},x;t)=\frac{1}{2}[\dot{x}-f(x)]^2.
  \label{e:action}
\end{equation}

Next, we calculate the probability of a transition
between two states. In the limit as $D\rightarrow 0$, this probability will be
increasingly dominated by the path that minimizes ${\mathcal S}[x(t)]$ so that
$P[x(t)]$ is a maximum.
This optimal fluctuational path $x_{\rm opt}$
is found as a solution to the
variational problem  $\delta\mathcal{S}[x]=0$, where
\begin{flalign}
\delta\mathcal{S}[x]&=\mathcal{S}[x+\delta x]-\mathcal{S}[x]=\int\limits_{t_i}^{t_f}\mathcal{L}(\dot{x}+\delta\dot{x},x+\delta x;t)\, dt
- \int\limits_{t_i}^{t_f}\mathcal{L}(\dot{x},x;t)\, dt\nonumber\\
&\nonumber \\
&=\int\limits_{t_i}^{t_f}\left [\frac{\partial \mathcal{L}}{\partial
    x}-\frac{d}{dt}\left(\frac{\partial\mathcal{L}}{\partial\dot{x}}\right)\right ]\delta x\, dt =0.
\end{flalign}
Since the variation $\delta x$ is arbitrary one is left with the
Euler-Lagrange equation
\begin{equation}
\frac{\partial \mathcal{L}}{\partial
    x}-\frac{d}{dt}(\frac{\partial\mathcal{L}}{\partial\dot{x}})=
    f(x)f'(x)-\ddot{x} = 0,\label{e:EL}
\end{equation}
which is then solved for the optimal path $x_{\rm opt}$. It is worth noting
that the optimal path has been found by
transforming the original stochastic problem into a new deterministic system
described by the Euler-Lagrange equation~(\ref{e:EL}).  The dimensions of the Euler-Lagrange equation are
twice the dimensions of the original system.

Following the standard progression of classical mechanics, it is sometimes of
use to move from the Lagrangian ${\mathcal L}$ expressed in terms of the state
$x$ and associated velocity $\dot{x}$ to the Hamiltonian ${\mathcal H}$
formulated as a function of $x$ and the conjugate momentum $\lambda$ defined as
\begin{equation}
\lambda=\frac{\partial\mathcal{L}}{\partial \dot{x}} = \dot{x}-f(x).
\end{equation}
Note that we are denoting the conjugate momentum variable as $\lambda$ rather
than the usual $p$ to avoid conflict with our use of $p$ for probability density and to emphasize the connection between conjugate momentum and the costate (i.e., Lagrange multiplier) from the optimal control formulation below.
One way to move from the Lagrangian ${\mathcal L}$ to the Hamiltonian ${\mathcal H}$ is through the Legendre transformation
\begin{equation}
\mathcal{L}=\lambda\dot{x}-\mathcal{H}(x,\lambda),
\end{equation}
where $\dot{x}$ is found as a function of $x$ and $\lambda$ from $\lambda=\frac{\partial\mathcal{L}}{\partial \dot{x}}$
using the inverse function theorem. Starting with~(\ref{e:SDE}), the
Legendre transformation yields
\begin{equation}
\dot{x} = f(x) + \lambda,\quad\mbox{and}\quad{\mathcal H} = \frac{\lambda^2}2
+ \lambda f(x),
\label{e:Lagrange}
\end{equation}
illustrating the connection between the classical momentum $\lambda$ and the optimal fluctuation $\eta$.

One can then write the action $\mathcal{S}$ as
\begin{equation}
\mathcal{S}[x(t)]= \int\limits_{t_i}^{t_f}\mathcal{L}(\dot{x},x;t)\,  dt =
\int\limits_{t_i}^{t_f}\left ( \lambda\dot{x}-\mathcal{H}(x,\lambda)\right )\,  dt.
\end{equation}
The optimal path is found as a solution to the variational problem
  $\delta\mathcal{S}[x,\lambda]=0$, where variations in both $x$ and $\lambda$ must be
  considered. This leads to
\begin{flalign}
&\delta\mathcal{S}[x,\lambda]=\mathcal{S}[x+\delta x,\lambda+\delta \lambda]-\mathcal{S}[x,\lambda]\nonumber\\
&\nonumber\\
&=\int\limits_{t_i}^{t_f}\left [(\lambda+\delta \lambda)(\dot{x}+\delta
  \dot{x})-\mathcal{H}(x+\delta x,\lambda +\delta \lambda) \right ]\, dt
- \int\limits_{t_i}^{t_f}\left [\lambda\dot{x}-\mathcal{H}(x,\lambda)\right ]\, dt\nonumber\\
&\nonumber\\
&=\int\limits_{t_i}^{t_f} \left [\left (\dot{x}-\frac{\partial\mathcal{H}}{\partial \lambda}\right )\delta \lambda - \left (\dot{\lambda}+\frac{\partial\mathcal{H}}{\partial
    x}\right )\delta x\right ]\, dt =0.
\end{flalign}
As with the Lagrangian formulation, both variations are arbitrary. Therefore
one is left with Hamilton's equations
\begin{equation}
\dot{x}=\frac{\partial\mathcal{H}}{\partial \lambda},\qquad\mbox{and} \qquad \dot{\lambda}=-\frac{\partial\mathcal{H}}{\partial x}.
\label{e:Hamilt}
\end{equation}

To avoid potential inverse function theorem issues, one could find the
Hamiltonian from the outset. Thus,
\begin{equation}
P[\eta (t)]\propto\exp\left (-\frac{1}{2D}\int_{t_i}^{t_f}\eta ^2\,dt\right
    )=\exp(-\mathcal{R}_{\eta}/D)\quad\mbox{and}\quad \mathcal{R}_{\eta}[\eta
  (t)]=\frac{1}{2}\int_{t_i}^{t_f}\eta ^2 \, dt,
\end{equation}
where we regard $R_\eta$ as an objective function.
To find the optimal trajectory, one must minimize the functional
\begin{equation}
\mathcal{R}[x,\eta,\lambda]=\mathcal{R}_{\eta}[\eta
  (t)]+\int_{t_i}^{t_f}\lambda [\dot{x}-f(x)-\eta (t)] \, dt,
  \label{e:OptCont}
\end{equation}
where $\lambda$ is the Lagrange multiplier enforcing the dynamic constraint~(\ref{e:SDE}). Therefore
\begin{equation}
P[x(t)]\propto\exp(-R/D),
\end{equation}
where the value $R$ is obtained from
\begin{equation}
R= {\rm{min}} \mathcal{R}[x,\eta,\lambda].
\end{equation}
We extremize $\mathcal{R}[x,\eta,\lambda]$
with respect to $\eta$, $\lambda$, and $x$ as follows:
\begin{equation}
\frac{\delta\mathcal{R}}{\delta\eta}=\mathcal{R}[x,\eta +\delta\eta
,\lambda]-\mathcal{R}[x,\eta,\lambda]=\int_{t_i}^{t_f}(\eta -\lambda)\delta \eta \, dt =0,
\end{equation}
so that
\begin{equation}
\lambda=\eta;
\end{equation}

\begin{equation}
\frac{\delta\mathcal{R}}{\delta\lambda}=\mathcal{R}[x,\eta,\lambda +\delta
\lambda]-\mathcal{R}[x,\eta,\lambda]=\int_{t_i}^{t_f}(\dot{x}-f(x)-\eta
(t))\delta \lambda \, dt =0,
\end{equation}
so that
\begin{equation}
\dot{x}=f(x)+\eta (t);\qquad\mbox{and}
\end{equation}

\begin{equation}
\frac{\delta\mathcal{R}}{\delta x}=\mathcal{R}[x +\delta
x,\eta,\lambda]-\mathcal{R}[x,\eta,\lambda]=\int_{t_i}^{t_f}(-\dot{\lambda}-\lambda\frac{\partial
f}{\partial x})\delta x \, dt =0,
\end{equation}
so that
\begin{equation}
\dot{\lambda}=-\lambda\frac{\partial
f}{\partial x}.
\end{equation}
Taken together we have
\begin{equation}
\dot{x}=f(x)+\lambda, \quad
\dot{\lambda}=-\lambda\frac{\partial f}{\partial x}, \quad\mbox{and}\quad
\mathcal{H}(x,\lambda)=\frac{\lambda ^2}{2}+\lambda f(x).
\end{equation}
By comparison with the Lagrangian formulation, it is clear that the Lagrange
multiplier plays the role of the conjugate momentum.

As mentioned in Sec.~\ref{s:noise}, in both formulations the optimal path is found by
transforming the original stochastic problem into a new deterministic system
described by a Hamiltonian with $\mathcal{H}(x,\lambda)=0$ when $x_i$ is a metastable state and $t_i=-\infty$.  The dimensions of the Hamiltonian are
twice the dimensions of the original system due to
the conjugate momenta $\lambda$.
The method amounts to finding a zero-energy trajectory
of an effective mechanical system, and at least one of the solutions to
the zero-energy Hamiltonian is in the set of optimal paths.

\subsection{Escape in finite and infinite time}
\label{s:escapeTime}

In Sec.~\ref{s:LDT} we formulated a criterion satisfied by optimal (i.e., most likely) diffusion paths produced by~(\ref{e:SDE}).  These paths must also satisfy time constraints and boundary conditions reflecting, for example, a transition between quasi-stable states either in finite time or over an arbitrarily long time.
The minimal action and associated paths strongly depend on the time horizon relevant to the computations.
Computing mean first exit times from $\Omega$ generated by~(\ref{e:SDECov})
requires minimizing the following action over all exit times $\tf$:
\be
{\mathcal S}_{\tf} =  \frac12\int_{\ti}^{\tf} (\dot{x}-f(x))^\dag a^{-1}(x)(\dot{x}-f(x))\,dt,
\label{e:Sinf}
\ee
where we generalize the presentation in Sec.~\ref{s:background} to higher dimension so that $x\in {\mathbb R}^n$ and $a(x) = \sigma(x)\sigma^\dag(x)$.
Performing this minimization over paths from points $\xinit$ to $\xf$ yields the function
\be
Q(\xinit,\xf) = \inf_{\tf} {\mathcal S}_{\tf}[x(t)\vert x(\ti)=\xinit, x(\tf)=\xf],
\label{e:quasipot}
\ee
known as the quasi-potential for its similarity in form to a potential defined
between any two points.  Indeed, in the case of overdamped Langevin dynamics~\cite{kra40,Risken96,gar03,vanKampen_book} when~(\ref{e:SDECov}) is a noise-driven
gradient flow so that
\be
f(x) = -\left(\partialderiv{U}x\right)^\dag,
\ee
the quasi-potential is intimately related to the potential $U$, and under additional assumptions $Q(\xinit,\xf) = 2(U(\xf)-U(\xinit))$.
Regardless of whether $f$ is a gradient, the quasi-potential determines the scaling law of exit rates $\tau(x):=\inf\{t:X_t\notin \Omega|X_i=x\}$ as $D\rightarrow 0$, with
\[
\bE[\tau(x)] \sim \exp\left(\inf_{y\notin \Omega}Q(x,y)/D\right).
\]

SDEs with sufficient complexity require a numerical approach to compute the quasi-potential~\cite{bryson_applied_2001}, starting either from the optimal control formulation in~(\ref{e:OptCont}) or from the Hamiltonian framework presented in~(\ref{e:Hamilt}).  To include the matrix $a(x)$ in~(\ref{e:Sinf}), assumed positive-definite to reflect a non-degenerate noise process in (\ref{e:SDECov}), the Hamiltonian is changed slightly to
\be
{\mathcal H}(x,\lambda) = \frac12\lambda^\dag a\lambda + \lambda^\dag f.
\ee
Finite-time paths can then be solved using a shooting method to solve the two-point boundary value problem expressed by
\begin{align}
& \dot{x} = f(x) + a\lambda\quad\mbox{and}\label{e:shootForward}\\
& \dot{\lambda} = -\left(\partialderiv{f}{x}\right)^\dag \lambda - \frac12\left[\partialderiv{}{x}(\lambda^\dag a\lambda)\right]^\dag,\label{e:shootBackward}
\end{align}
where we have used the convention of representing vectors as columns and derivatives with respect to those vectors as rows.  Given a typical exit problem with fixed initial condition $\xinit$ and exit criterion forcing $x(\tf)$ to lie on the boundary $\partial\Omega$ of set $\Omega$, for instance, a straightforward shooting method would start with an initial guess for $\lambda$ (e.g., $\lambda(t)\equiv 0$), integrate~(\ref{e:shootForward}), adjust $\lambda(\tf)$ and integrate~(\ref{e:shootBackward}) backwards, then iterate this process until the terminal condition is satisfied~\cite{chernykh2001,elgkam04}.

Different alternatives exist for finding optimal paths in the case of gradient
systems, including the nudged elastic band method~\cite{jonsson_nudged} and
the string method~\cite{e_simplified_2007}.  A powerful method for more
general dynamical systems is the minimum action method (MAM), which can be
adapted via a rescaling in time to accommodate infinite-time optimal paths,
e.g., where either $t_i=-\infty$ in~(\ref{e:Sinf}) or~(\ref{e:quasipot}) or
the minimization over $t_f$ in those expressions yields $t_f=\infty$.  It sets
up the minimization as the gradient descent of the action ${\mathcal S}_{\tf}$.  This can
be accomplished in a number of ways, including by differentiating out the dependence of Hamilton's equations on the quasi-momentum $\lambda$ and advancing a PDE whose steady state is the optimal path sought.  The adaptation of this method to infinite-time
problems, referred to as the geometric minimum action method
(GMAM)~\cite{heymann_geometric_2008}, reparameterizes the optimal path to
depend on arclength rather than time, removing the obvious difficulty of
enforcing boundary conditions at infinity.  The iterative action minimization
method (IAMM) is a similar method that employs a direct, fully explicit
iterative scheme~\cite{LindSch13}. Basic details of both the GMAM and the IAMM can be found
in Appendix~\ref{s:appB}.

\subsubsection{Escape from potential well}
\label{s:potentialWell}

As an example of an infinite-time calculation, consider a particle in a potential well under the influence of external
Gaussian noise. The noise causes the particle to escape from
the well in which it resides, an event that grows increasingly rare as the strength of the noise decreases. To determine the mean time to escape (MTE) in this limit, we
consider again the Langevin equation
\begin{equation}
\dot{x}=-\frac{\partial U}{\partial x} + \eta (t),\label{eq:Langevin2}
\end{equation}
where $x$ is the particle's position, $U$ determines the potential well, and
$\eta$ is additive Gaussian noise as described by~(\ref{e:noiseProps}).
The associated Fokker-Planck equation is
\begin{equation}
\frac{\partial p(x,t)}{\partial t}= \frac{\partial}{\partial x}\left [
  \frac{\partial U(x)}{\partial x}p(x,t) \right ] +\frac12 D\frac{\partial ^2
  p(x,t)}{\partial x^2},\label{eq:FP}
\end{equation}
where $p$ is the probability density. The first term on the right-hand side is
the drift or transport term, and the second term on the right-hand side is the
diffusion or fluctuation term. It is possible to rewrite~(\ref{eq:FP}) as
\begin{equation}
\frac{\partial p}{\partial t}=-\frac{\partial J}{\partial
  x}, \qquad  J=-\frac{\partial U}{\partial x}p-\frac12 D\frac{\partial p}{\partial x}.\label{eq:PC}
\end{equation}
Therefore, $\frac{\partial p}{\partial t}+\frac{\partial J}{\partial x}=0$ is a
  continuity equation for the probability density $p$, and $J$ is interpreted as a probability current. For a stationary process, $J=\,$constant.

\begin{figure}[t!]
\begin{center}
\includegraphics[scale=0.5]{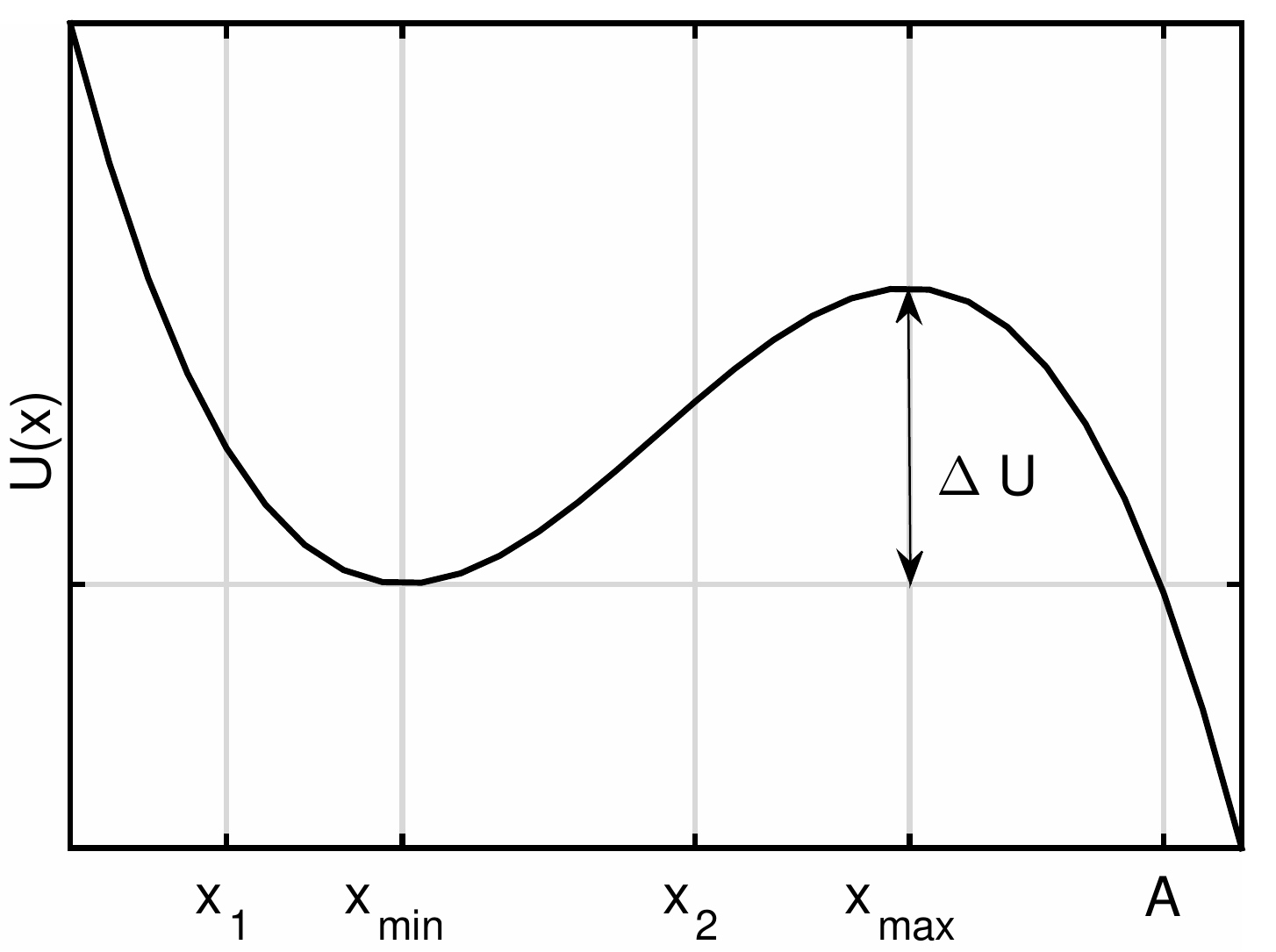}
\end{center}
\caption{\label{fig:potential}Schematic of a potential well $U(x)$. The local minimum
of the well is located at $x=x_{\rm min}$, while the top of the barrier is
located at $x=x_{\rm max}$. }
\end{figure}

Figure~\ref{fig:potential} shows a schematic of a potential well $U(x)$. We
assume the well height is much larger than the noise intensity so that $\Delta
U / D \gg 1$. We will determine the escape rate/time for particles sitting in
a deep well near $x=x_{\rm min}$.

The probability current given by~(\ref{eq:PC}) can be rewritten as
\begin{equation}
J=-\frac12D\exp{(-2U(x)/D)}\frac{\partial}{\partial x}\left [ \exp{(2U(x)/D)}p(x,t)
\right ].
\end{equation}
Integration from $x_{\rm min}$ to $A$ gives, under the assumption of quasi-stationarity (i.e., $J$ is constant),
\begin{equation}
J=\frac{D\exp{(2U(x_{\rm min})/D)}p(x_{\rm min},t)}{2\int\limits_{x_{\rm min}}^A\exp{(2U(x)/D)}\,dx},
\end{equation}
where $x=A$ is the location of an arbitrary point far to the right of the
location of the top
of the well barrier (see Fig.~\ref{fig:potential}).
For a high barrier the quasi-stationary distribution of the Fokker-Planck equation satisfies
\begin{equation}
0\approx \frac{\partial}{\partial x}\left [
  \frac{\partial U(x)}{\partial x}p(x,t) \right ] + \frac12 D\frac{\partial ^2 p(x,t)}{\partial x^2},
\end{equation}
where the time dependence of $p(x,t)$ arises from a slow leakage of probability from the well.
The quasi-stationary distribution function near $x_{\rm min}$ is then
\begin{equation}
p(x,t)=p(x_{\rm min},t)\exp{(-2[U(x)-U(x_{\rm min})]/D)}.
\end{equation}
The probability $P(t)$ to find the particle in the well is
\begin{equation}
P=\int\limits_{x_1}^{x_2} p(x,t) \, dx = p(x_{\rm min},t)\exp{(2U(x_{\rm
    min})/D)}\int\limits_{x_1}^{x_2} \exp{(-2U(x)/D)}\, dx,
\end{equation}
where $x=x_1$ and $x=x_2$ are the locations of arbitrary points in the
potential well found respectively to the left and the right of $x=x_{\rm min}$
(see Fig.~\ref{fig:potential}).
Since probability $P$ times the escape rate $r$  is the probability current
$J$, the inverse escape rate (escape time) is
\begin{equation}
\tau =
\frac{1}{r}=\frac{P}{J}=\frac{2}{D}\int\limits_{x_1}^{x_2}\exp{(-2U(x)/D)}\,dx\int\limits_{x_{\rm min}}^{A}\exp{(2U(x)/D)}\,dx.
\end{equation}

The first integral is largest near $x_{\rm min}$. We expand
\begin{equation}
U(x)\approx U(x_{\rm min})+\frac{1}{2}U^{\prime\prime}(x_{\rm min})(x-x_{\rm
  min})^2.
\end{equation}
The second integral is largest near $x_{\rm max}$. Again, we expand
\begin{equation}
U(x)\approx U(x_{\rm max})-\frac{1}{2}|U^{\prime\prime}(x_{\rm max})|(x-x_{\rm
    max})^2.
\end{equation}
By extending the integral boundaries to $\pm \infty$ in both
  directions, one obtains the expression for the MTE as
\begin{equation}
\tau =\frac{2\pi}{\sqrt{U^{\prime\prime}(x_{\rm min})|U^{\prime\prime}(x_{\rm
      max})|}}\exp{(2[U(x_{\rm max})-U(x_{\rm min})]/D)} =K\exp{(2\Delta U/D)}.\label{eq:MTE}
\end{equation}
It is important to note in~(\ref{eq:MTE}) that the escape time increases
exponentially with increasing barrier height and decreasing noise intensity, and the prefactor $K$ depends on the curvatures at
$x_{\rm min}$ and $x_{\rm max}$.

\begin{figure}[t!]
\begin{minipage}{0.49\linewidth}
\begin{center}
\includegraphics[scale=0.42]{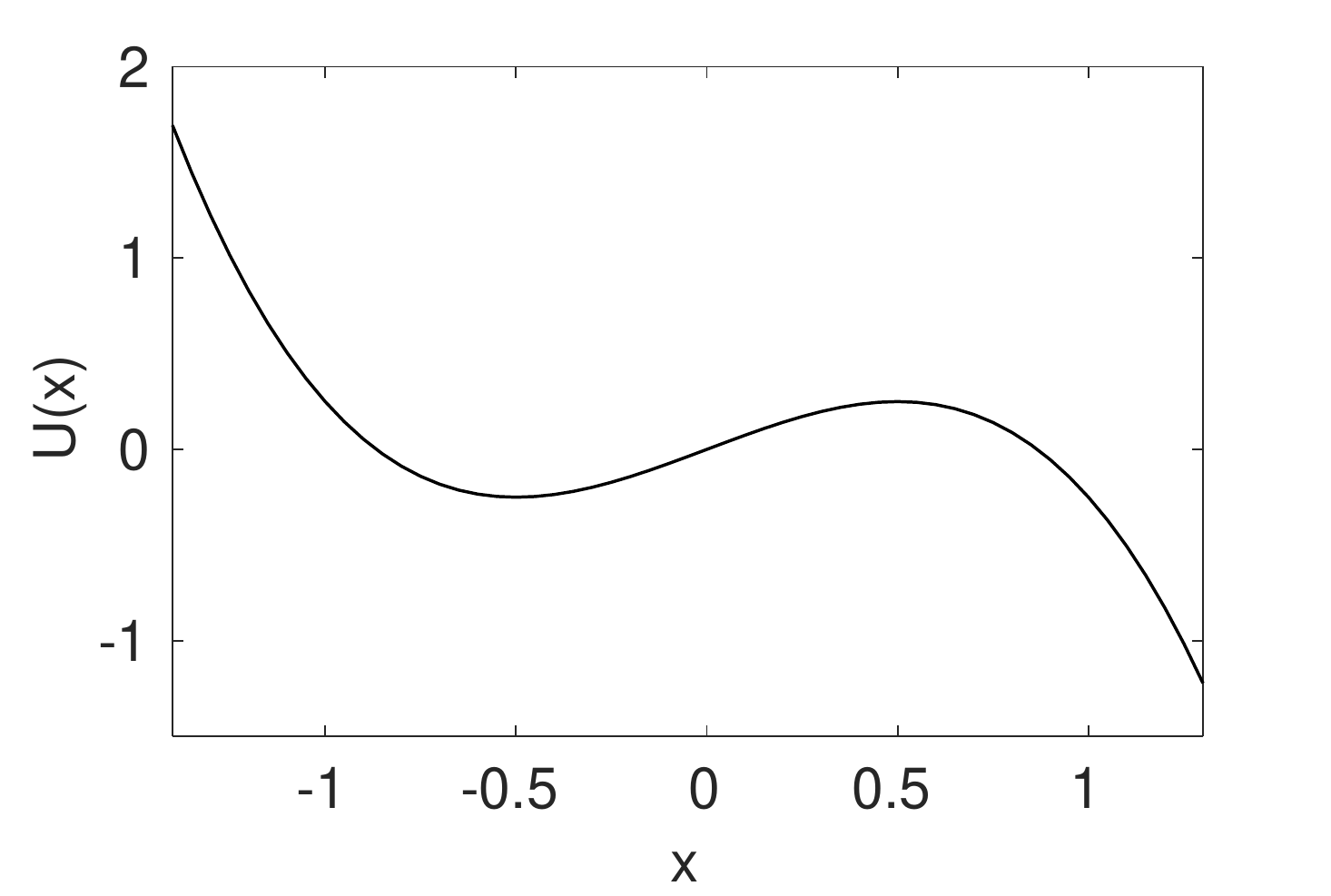}
\end{center}
\end{minipage}
\begin{minipage}{0.49\linewidth}
\begin{center}
\includegraphics[scale=0.42]{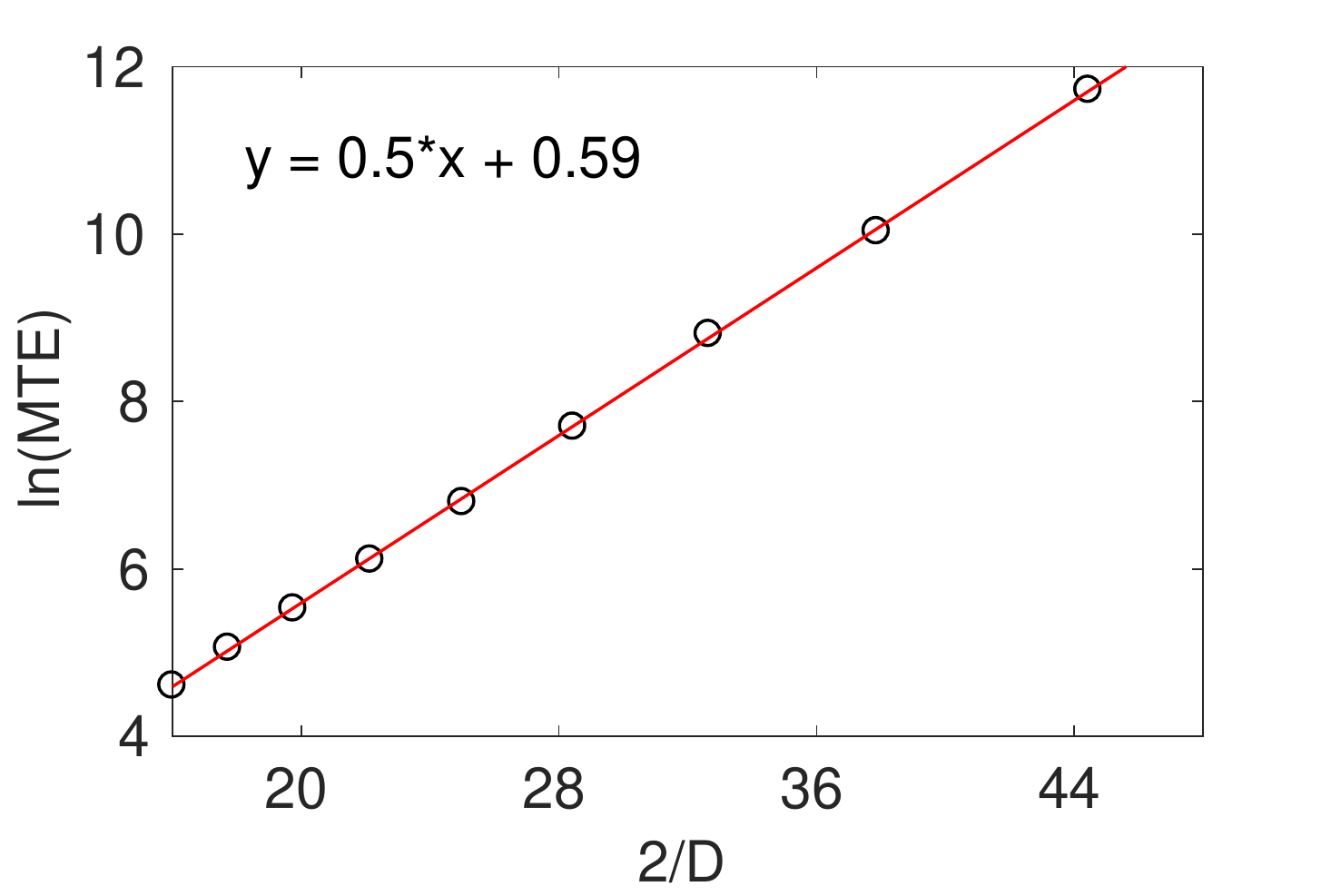}
\end{center}
\end{minipage}
\caption{\label{fig:ex}(a) Graph of the potential well given by
  $U(x)=-x^3+\frac{3}{4}x$ with associated well height $\Delta
  U=\frac{1}{2}$. (b) Logarithm of the numerically computed mean time to
  extinction (black circles) as a function of $2/D$ with a line of best fit (red line) passing through
the data. The equation of the best fit line is $y=0.5x+0.59$, whose slope
agrees perfectly with the analytical prediction of $\Delta U=\frac{1}{2}$.}
\end{figure}

A specific example is shown in Fig.~\ref{fig:ex}(a)-(b).
Figure~\ref{fig:ex}(a) shows the potential well given by
$U(x)=-x^3+\frac{3}{4}x$ with associated barrier height $\Delta
U=\frac{1}{2}$. We numerically integrate~(\ref{eq:Langevin2}) using the
Euler-Maruyama method (Appendix~\ref{s:appA}), and the time needed for a
particle to escape from the well is determined. By performing the computation
for 10,000 particles, the mean time to escape is computed for a variety of
noise intensities. By taking the logarithm of~(\ref{eq:MTE}), one has
$\ln{\tau}=\ln{K}+\frac{2\Delta U}{D}$. By plotting $\ln{\tau}$ versus $2/D$,
the analytical prediction is that the data should lie along a line with slope
$m=\Delta U$. Figure~\ref{fig:ex}(b) shows the numerical results for the
potential shown in Fig.~\ref{fig:ex}(a) along
with a line of best fit through the data. One can see excellent agreement
between the numerical and analytical results.

Returning to the SDE description, the optimal fluctuation equations given by~(\ref{e:Hamilt}) are easily shown in the gradient case under consideration to give
\be
\ddot{x} = \partialderiv{U}{x}\partialderiv[2]{U}{x},
\label{e:UpUpp}
\ee
which integrates once to give
\be
\dot{x}^2=\left(\partialderiv{U}{x}\right)^2+c.
\ee
In the case of infinite-time exit, the constant $c=0$, leading to an asymptotic dynamic behavior where the conjugate momentum $\lambda$ is either zero or acts against the gradient, i.e.,
\be
\dot{x} = -\partialderiv{U}{x} + \lambda = \mp\partialderiv{U}{x},
\label{e:withOrAgainst}
\ee
depending on whether or not the uncontrolled dynamics (i.e., the deterministic dynamics obtained by setting $\eta\equiv 0$ in~(\ref{eq:Langevin2})) leads $x$ in the direction of an exit.  If $x(\ti)<-1/2$, the uncontrolled dynamics takes infinite time to descend into the stable fixed point, after which the controlled dynamics follows the trajectory
\be
x(t) = \frac12\tanh(\frac32t),
\ee
taking infinitely long both to rise out of the stable fixed point and to approach the saddle, thereby effecting an exit.  The minimum action for any initial condition to the left of the stable fixed point is ${\mathcal S}=2\Delta U=1$.

To see how this changes when a finite time horizon is imposed, note that in this case~(\ref{e:UpUpp}) no longer implies~(\ref{e:withOrAgainst}).
\begin{figure}[t!]
\begin{minipage}{0.49\linewidth}
\begin{center}
\includegraphics[scale=0.3]{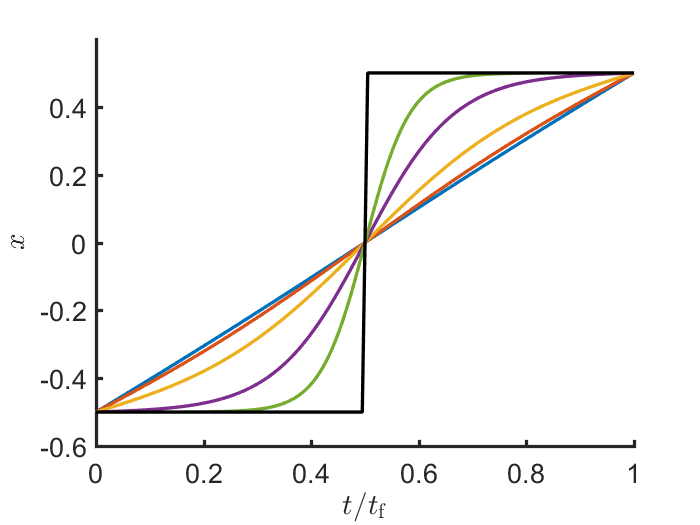}
\end{center}
\end{minipage}
\begin{minipage}{0.49\linewidth}
\begin{center}
\includegraphics[scale=0.3]{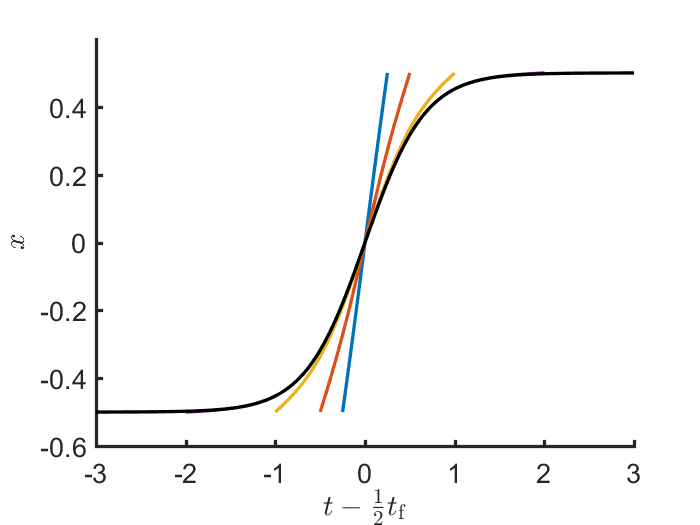}
\end{center}
\end{minipage}\\
\begin{minipage}{0.49\linewidth}
\begin{center}
\includegraphics[scale=0.3]{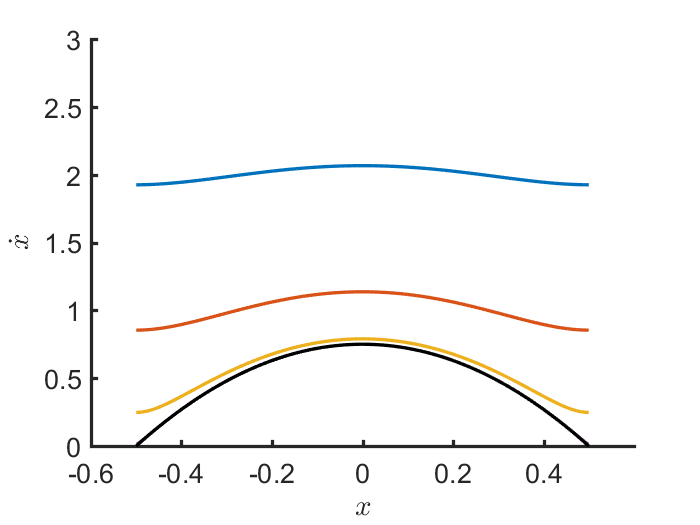}
\end{center}
\end{minipage}
\begin{minipage}{0.49\linewidth}
\begin{center}
\includegraphics[scale=0.3]{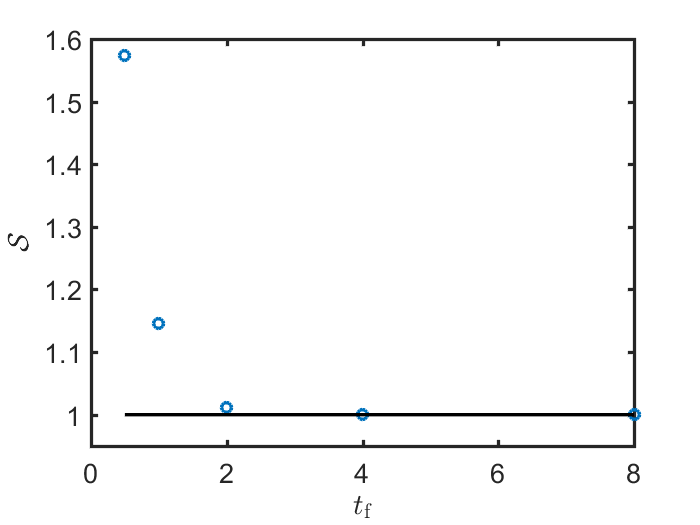}
\end{center}
\end{minipage}
\caption{\label{fig:finiteTime}Optimal paths computed from the stable fixed
  point to the saddle for final times $\tf =1/2$ (blue), $\tf =1$ (red), $\tf =2$ (orange), $\tf =4$ (purple), $\tf =8$ (green) and $\tf =\infty$ (black). (a) Top left: Optimal paths rescaled to unit exit time.  (b) Top right: Optimal paths with time shifted so that passages through zero coincide. (c) Bottom left: Phase portrait ($\dot{x}$ vs.\ $x$). (d) Bottom right: Minimum action ${\mathcal S}$ vs.\ $\tf$. }
\end{figure}
Figures~\ref{fig:finiteTime}(a)-(d) show the optimal paths and their phase portraits associated with $\tf\in\{1/2, 1, 2, 4, 8\}$ and compares them to the infinite-time optimal path.  Whereas for short time $\tf=1/2$ the optimal path is near a straight line path (more generally, a geodesic with respect to the metric imposed by $a(x)$), the optimal path quickly approaches the infinite-time path as $\tf$ increases.  Once these paths are indistinguishable, i.e., by $\tf=4$, increases in $\tf$ simply add to the time spent vanishingly close to either fixed point.  By this value of $\tf$, the action has also converged to its quasi-potential value of ${\mathcal S}=Q(-1/2,0)=1$.  The phase portrait depicted in Fig.~\ref{fig:finiteTime}(c) demonstrates again that for short times $\tf$ the optimal path tends to a constant velocity that generates the straight-line path depicted in Fig.~\ref{fig:finiteTime}(b); however, for large times the path through phase space collapses onto the anti-gradient flow described above.

Finally, we note that the reference to ``infinite-time optimal paths'' implies passage from, to, or through a hyperbolic fixed point. The optimization suggested by~(\ref{e:Sinf}) performed on paths connecting two points that are not hyperbolic fixed points of the dynamical system and that do not pass through one, e.g., $x(\ti)=-1/4$ and $x(\tf)=1/4$, yields an action determined by the quasi-potential $Q = 2\Delta U = 11/16$ and an optimal path time of
\be
\Delta t = \frac43\tanh^{-1}\frac12 =0.73.
\ee

\subsection{Master equation formalism for internal noise}
\label{s:mastereq}

Just as we considered the mean exit time due to external noise in Sec.~(\ref{s:background}), we would like to determine the mean time to extinction (MTE) and the optimal path for internal noise. Assuming again that we can model stochastic
population dynamics using a Markov process, we describe the evolution of the probability density using the
  master equation
\begin{equation}
\partialderiv{\rho (X,t)}{t}=\sum\limits_r \left  [ W_r(X-r)\rho (X-r,t) - W_r(X)\rho (X,t) \right ],\label{eq:ME}
\end{equation}
 where $\rho (X,t)$ is the probability of finding $X$ individuals at time $t$,
 and $W_r(X)$ is the transition rate from $X$ to $X+r$, where $r$ is a
  positive or negative integer increment.
The first term of~(\ref{eq:ME}) is the gain to state $X$ from state
$X-r$. The second term of~(\ref{eq:ME}) is the loss from state $X$ to
  other states.

The master equation is a large, or even infinite, set of differential
  equations. A diffusion approximation leading to a Fokker-Planck equation as in~(\ref{e:fFPE}) is often used but is very (exponentially)
  inaccurate~\cite{gamoto96,Doering2005,Kessler07,Assaf2007}. Assuming that
  the typical size of the population $X\sim K$ as determined by the maximum of
  the quasi-stationary probability distribution, we can use a WKBJ
  approximation of the master equation~\cite{Kubo1973,gang87,Dykman1994d,elgkam04,Kessler07}.
If the mean-field equation associated with the original stochastic problem is
$n$-dimensional, then the WKBJ method transforms the problem into one of classical mechanics in
$2n$-dimensional phase space where the doubling of dimension is due to the
conjugate momenta. It is important to note that these deterministic
dynamics are not the system size $K\to\infty$ mean-field dynamics of the original
stochastic model, but rather provide auxiliary dynamics that describe
processes not captured by the mean-field dynamics.

For large populations, the time to extinction is long. The mean time to
extinction is determined by the probability flux (extinction rate) into the
extinct state, and is determined by the tail of the quasi-stationary
probability distribution where $\frac{\partial \rho}{\partial t}\approx 0$. Therefore
\begin{equation}
0=\sum\limits_r \left [W_r(X-r)\rho (X-r) -  W_r(X)\rho (X) \right ].
\end{equation}
Let $X$ be scaled by $K$, the typical population size in the metastable state.  Using $x=X/K$, the transition rate $W_r(X)=W_r(Kx)$
can be represented as the following expansion in $K$,
\begin{equation}
W_r(Kx)=Kw_r(x)+u_r(x) +\mathcal{O}(1/K),
\label{e:transrates}
\end{equation}
where $x$ and the scaled transition rates $w_r$ and $u_r$ are
$\mathcal{O}(1)$. Additionally we can write,
\begin{equation}
\rho(X)\equiv \rho(Kx)=\pi(x).
\end{equation}

Then the scaled master equation is
\begin{equation}
\sum\limits_r \left [w_r(x-\frac{r}{K})\pi (x-\frac{r}{K}) -  w_r(x)\pi (x) \right ]=0.
\end{equation}

For $K\gg 1$ we approximate the scaled master equation using the
WKB approximation.
To apply the WKB approximation, we assume that
\begin{equation}
\pi(x)=A(x)\exp{[(-K\mathcal{S}(x))(1+\mathcal{O}(1/K))]}.
\end{equation}
We substitute the WKB ansatz along with
  $\mathcal{S}(x-\frac{r}{K})=\mathcal{S}(x)-\frac{r}{K}\mathcal{S}^{\prime}(x)+\frac{r^2}{2K^2}\mathcal{S}^{\prime\prime}(x)+\mathcal{O}(1/K^3)$ and
$A(x-\frac{r}{K})=A(x)-\frac{r}{K}A^{\prime}(x)+\mathcal{O}(1/K^2)$ into the
master equation, and at leading order one obtains a Hamilton-Jacobi equation with  Hamiltonian
\begin{equation}
\mathcal{H}(x,\lambda)=\sum\limits_r w_r(x)[\exp{(r\lambda)-1}]=0, \qquad
\lambda=\frac{\partial\mathcal{S}}{\partial x}.
\label{e:Hamil}
\end{equation}

Hamilton's equations are
\begin{subequations}
\begin{flalign}
\dot{x}=&\frac{\partial\mathcal{H}}{\partial \lambda}=\sum\limits_r rw_r(x)\exp{(r\lambda)},\\
\dot{\lambda}=&-\frac{\partial\mathcal{H}}{\partial x}=-\sum\limits_r
[\exp{(r\lambda)}-1]\frac{\partial w_r(x)}{\partial x}.
\end{flalign}
\label{e:Hameqnsi}
\end{subequations}
In summary, the WKB solution amounts to finding a zero-energy trajectory of an
effective mechanical system, and at least one of the solutions
to the zero-energy Hamiltonian is the optimal path.
 The $x$ dynamics along the $\lambda=0$ deterministic line are described by
\begin{equation}
\dot{x}=\frac{\partial \mathcal{H}(x,\lambda)}{\partial \lambda} \rvert_{\lambda=0}
\label{e:MFEq}
\end{equation}
which is the rescaled mean-field rate equation associated with
the original deterministic problem.  For a  single step process, this simplifies to
$ \dot{x}=  w_{1}(x) - w_{-1}(x)$.

In Sec.~\ref{s:extinction}, we will discuss the application of this master
equation formalism to two examples involving extinction in population biology.
For the simple models discussed in Sec.~\ref{s:extinction}, the deterministic
steady states are nodes. It is easy to show that the WKB method
transforms these steady state nodes in the original 1D deterministic system into
steady state saddle points in the 2D set of Hamilton's equations. This
allows for noise-induced escape from a metastable state
and provides a path to
extinction that did not exist in the original deterministic model.

\subsection{Mean Time to Extinction}
\label{s:MTEn}

In a single step process, such as the two examples considered in Sec.~\ref{s:extinction}, the optimal path $\lambda_{opt} (x)$ will always have the general form
\begin{equation} \label{e:gen_popt}
\lambda_{opt}(x)=-\ln { \left( w_{+1}(x) /w_{-1}(x)  \right)}.
\end{equation}
Using the definition of the conjugate momentum $\lambda=d{\cal S}/dx$, the action ${\cal S}_{opt}$
along the optimal path $\lambda_{opt} (x)$ is given by
\begin{equation}
\label{e:Sopt}
{\cal S}_{opt} = \int^{x} \lambda_{opt} (x) dx,
\end{equation}
where the limits of integration are determined by the appropriate equilibrium
points of the specific problem.
Therefore, the mean time to extinction (MTE) can be approximated by
\begin{equation}\label{eq:tau_ext_gen}
\tau = B \exp({K {\cal S}_{opt}}) ,
\end{equation}
where $B$ is a prefactor that depends on the system parameters and on the
population size. An accurate approximation of the MTE depends on obtaining
$B$.
The specific form of the prefactor differs for the two types of extinction
problems considered in Sec.~\ref{s:extinction}, and are provided in that section.

\subsection{Sampling for exits}
\label{s:exitSampling}

In the simple gradient flow discussed in Sec.~\ref{s:potentialWell} it is possible to obtain or approximate not only the exponential scaling law in $D$ of the exit probability over finite times, but also the normalization constant that provides the pre-exponential factor $K$ in~(\ref{eq:MTE}).  In flows that are more complex, including non-potential flows that do not satisfy detailed balance~\cite{vanKampen_book} and those of higher (or infinite) dimension, this is no longer possible analytically.  Nor is it practical computationally, due to the formation of boundary layers with the singular perturbation introduced as $D\rightarrow 0$ in Dynkin's equation~(\ref{e:Dynkin}) or the time-dependent Fokker-Planck equation~(\ref{e:fFPE}).

An alternate method to obtain the normalization constant uses
sampling through
Monte Carlo simulations.
This too is subject to the ``curse of dimensionality'' in cases of high dimension such as discretizations of stochastic PDEs; however, variance reduction techniques can be used to improve the numerical efficiency of obtaining an estimate.  We will focus here on importance sampling, a variance reduction method that uses information from the controlled dynamics of the deterministic system to inform a bias applied to the distribution from which the Monte Carlo simulations are drawn.

Suppose we wish to calculate the probability that, after time $\tf$, a stochastic process $x(t)$ evolving according to~(\ref{e:SDECov}) terminates outside of set $\Omega$, i.e., with $x(\tf)\notin \Omega$, where we assume for concreteness that the deterministic evolution $\dot{x} = f(x,t)$
yields $x(t)\in \Omega$ for $x(\ti) = \xinit$ and $t\leq \tf$.  In other words, we wish to calculate the probability that the stochastic forcing term $\sigma(x,t)\eta$ drives the state some distance away from its deterministic trajectory.  The standard Monte Carlo estimator for this probability $P$ is given by
\begin{equation}\label{equ:RegProbEst}
\Des{\hat{P}}{MC} = \frac{1}{N} \sum_{k=1}^{N} I(X_{\tf}^{(k)})
\end{equation}
where $\{X_{\tf}^{(k)}, k=1,\dots,N\}$ are $N$ independent simulations of~(\ref{e:SDE}) and $I(x)$ is an indicator function evaluating to 0 if $x\in \Omega$ and 1 otherwise.
Thus, $I(X_{\tf})$ forms a Bernoulli random variable, implying that the estimator $\Des{\hat{P}}{MC}$ has variance given by
\begin{equation}\label{equ:VarEst}
\Var{\Des{\hat{P}}{MC}} = \dfrac{P(1-P)}{N},
\end{equation}
with a relative error (coefficient of variation) of,
\begin{equation}\label{equ:CVar}
\CV{\Des{\hat{P}}{MC}} = \frac{\sqrt{\Var{\Des{\hat{P}}{MC}}}}{\Expect{\Des{\hat{P}}{MC}}} = \dfrac{\sqrt{1-P}}{\sqrt{NP}}.
\end{equation}
Ensuring that \eqref{equ:RegProbEst} is an accurate estimator for the true probability requires that
$\CV{\Des{\hat{P}}{MC}} \ll 1$.  However, if the set $\Omega$ and noise strength $D$ are such that $P \ll 1$
(i.e., $X_{\tf}\notin \Omega$ is a rare event), this requirement is approximately expressed by $N \gg 1/P$, which
implies that a very large (and often unattainable) number of Monte Carlo runs is required to obtain an accurate estimate.

The idea behind importance sampling for diffusions is to replace the original dynamics given by~(\ref{e:SDECov}) with
\be
\dot{x} = f(x,t) + \sigma(x,t)(b+\eta),
\label{e:biasedSDE}
\ee
so that the noise term $\eta$ is biased in such a way that the events of interest, $x({\tf})\notin \Omega$, are much more likely.  The particular form of bias suggested by~(\ref{e:biasedSDE}) is mean-biasing; other forms of biasing the noise are possible but less common.  Mean biasing effectively replaces increments $\Delta\eta$ drawn from the original distribution (e.g., Gaussian) $p(\Delta\eta)$ by shifted increments
\be
\Delta\tilde\eta = \Delta\eta + b\Delta t,
\ee
which is equivalent to drawing $\Delta\tilde\eta$ from biased distribution $\tilde p(\Delta\tilde\eta) = p(\Delta\tilde\eta-b\Delta t)$.
Despite the biased dynamics described by~(\ref{e:biasedSDE}), an unbiased estimator can be recovered from Monte Carlo simulations of this equation through an appropriate weighting referred to as the likelihood ratio
\be
l(\Delta\tilde\eta) = \frac{{\tilde p}(\Delta \tilde\eta)}{p(\Delta \tilde\eta)}
\ee
with the likelihood ratio for the $k^\mathrm{th}$ simulated path given by the product of these increments.
The unbiased estimator is then
\begin{equation}\label{equ:ISProbEst}
 \Des{\hat{P}}{IS}  = \frac{1}{N}\sum_{k=1}^{N} I({\tilde X}_{\tf}^{(k)})l^{(k)},\quad\mathrm{with}\quad l^{k} = \prod_j l(\Delta\tilde\eta^{(k)}_j)
\end{equation}
and with variance
\begin{equation}\label{equ:ISVarEst}
\mathbb{{\tilde V}}\hspace{-0.03in}\left[\Des{\hat{P}}{IS}\right] = \dfrac{\mathbb{{\tilde E}}\left[I({\tilde X}_{\tf})l^2\right]-P^2}{N} \approx \dfrac{\mathbb{E}\left[I({\tilde X}_{\tf})l\right]}{N},
\end{equation}
where ${\mathbb{\tilde E}}$ and ${\mathbb{\tilde V}}$ denote expectation and variance with respect to the biased noise process.

A ``good'' importance sampling density is one that produces a small variance,
i.e., one that satisfies $\tilde p \gg p$ for the events of interest,
$x({\tf})\notin \Omega$.  In the present example, the biased density is
determined by the control term in~(\ref{e:biasedSDE}).  In general,
determining the control is a difficult task and one that depends sensitively on the deterministic dynamics and noise covariance structure.
To motivate a particularly effective choice of $b$, we appeal again to the optimal control problem described above arising in the theory of large deviations.  In particular, note that~(\ref{e:biasedSDE}) essentially just adds noise to~(\ref{e:shootForward}) with control $b=\sigma^\dag \lambda$.  Thus, building a control that satisfies~(\ref{e:shootBackward}) provides paths that form the most likely routes to the rare events of interest in~(\ref{e:SDE}).

\section{Applications}
\label{s:applications}

\subsection{Phase slips in mode-locked lasers}
\label{s:phaseSlips}

\begin{figure}[t!]
\centerline{\includegraphics[width=.7\linewidth]{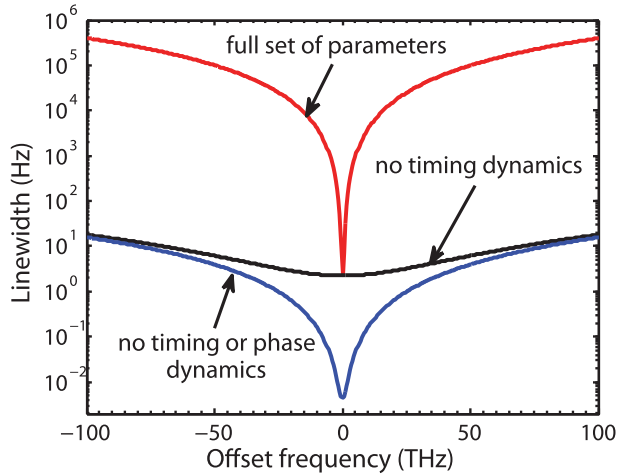}}
\caption{\label{f:bikini} Linewidth analysis of a mode-locked laser, demonstrating the dominant role played by timing and phase dynamics in determining linewidth.  Reproduced from Ref.~\cite{bao_pulse_2014}.}
\end{figure}

\begin{figure}[t!]
\centerline{\includegraphics[width=.7\linewidth]{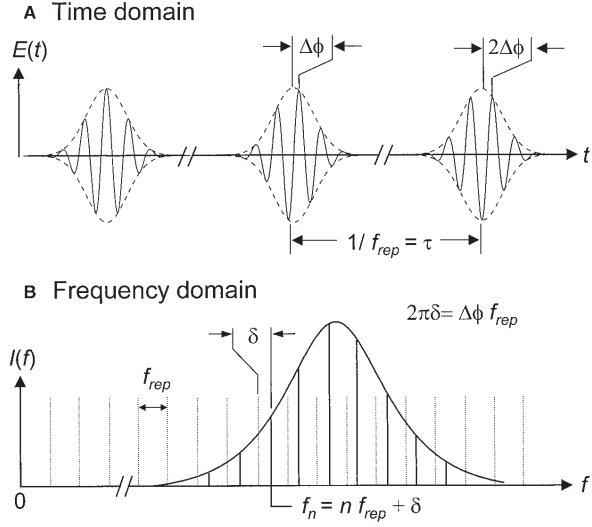}}
\caption{\label{f:mll} Sketch of electric field output from a mode-locked laser, in time (top) and in frequency space (bottom).  A phase slip occurs when $\Delta\phi$ undergoes a noise-driven rotation of $2\pi$. Reproduced from Ref.~\cite{jones_carrier-envelope_2000}.}\label{f:MLL}
\end{figure}

Mode-locked lasers provide an important example of systems subject to failure caused by large noise-induced phase excursions.  These lasers produce pulses of light whose electric field envelope's phase is locked to the phase of the underlying carrier wave.  In frequency space, this corresponds to a highly regular comb of laser lines with broad bandwidth and narrow linewidth~\cite{jones_carrier-envelope_2000}.  Important applications of the temporal and spectral properties of mode-locked lasers depicted in Fig.~(\ref{f:MLL}) include laser manipulation of ultrafast chemical reactions and extremely precise measurement of time intervals~\cite{jun_ye_optical_2003}.

A typical analysis of the effect of noise on laser operation is restricted to the induced linewidth obtained by linearizing the model equations for the laser about their desired operational state to produce Ornstein-Uhlenbeck dynamics; see for example~\cite{bao_pulse_2014},
\[
\dot{x} = -Ax + \eta,\quad x=(\Delta g,\Delta w,\Delta\omega,\Delta\tau,\Delta\theta)^\dag,
\]
where $\Delta g$, $\Delta w$, $\Delta\omega$, $\Delta\tau$, and $\Delta\theta$ represent fluctuations in saturated gain, pulse energy, central frequency, central pulse timing, and phase, respectively.  The conclusion suggested by Fig.~\ref{f:bikini} is that the timing and phase dynamics provide the dominant mechanisms by which noise affects linewidth.  As explained below, however, this analysis neglects another important source of uncertainty.

The generation of a frequency comb with narrow individual lines over a broad spectral region requires precise control over the carrier-envelope phase offset, depicted as $\Delta\phi$ in Fig.~\ref{f:mll}.  This control is implemented through interferometric feedback that typically only measures $\Delta \phi$ modulo full rotations of $2\pi$.  If a sequence of noise events successfully pushes $\Delta\phi$ to $\pm\pi$ (assuming without loss of generality that $\Delta\phi=0$ for the noise-free operational state) then the feedback drives $\Delta\phi$ to the neighboring equilibrium value of $\pm 2\pi$.  This situation is referred to as a phase slip and is an additional source of uncertainty in the laser's output.  A noise-driven damped pendulum provides a good conceptual framework for these two sources of uncertainty, where the linewidth is provided by the stationary measure obtained in the small-amplitude limit of a damped linear oscillator, while phase slips correspond to the occasional full rotations induced in the pendulum by very unlikely sequences of noise increments.

To account for these rare events, consider the heuristic mode-locked laser
model with  active modulation introduced in Ref.~\cite{donovan_rare_2007}:
\be
i\partialderiv{u}{t}+\frac12\partialderiv[2]{u}{\xi}+|u|^2u =
-c_0\cos(\omega \xi)u-ic_1u +ic_2\partialderiv[2]{u}{\xi} +id_1|u|^2u -id_2|u|^4u + i \eta(\xi,t),
\label{e:mll}
\ee
where $u$ is the normalized complex electric field envelope, $\xi$ is the averaged
transverse variable (e.g., time, for temporal pulses as seen on an
oscilloscope), and $t$ is the longitudinal variable (e.g., a continuous
approximation of the cavity round-trip index).  Constants $c_1$, $c_2$, $d_1$, and $d_2$
respectively represent coefficients of linear cavity loss, spectral filtering,
nonlinear gain, and gain saturation.  Constants $c_0$ and $\omega$ are the
amplitude and inverse width of active modulation intended to stabilize the
pulse with its center position at the origin, and 
$\eta(\xi,t)$ is space-time white noise, with
\be
{\mathbb E}[\eta^\dag(\xi,t)\eta(\xi',t')] = D\delta(\xi-\xi')\delta(t-t').
\label{e:ansatz}
\ee

The stabilization term in~(\ref{e:mll}) provides a restoring force on the pulse parameters that mitigates the impact of disturbances, including noise.  This effect is characterized below through a linearized analysis of the reduced system of pulse parameters.  The stabilization has a secondary effect, however, of producing additional equilibria that become dynamically accessible upon the introduction of noise.  Transitions to these equilibria represent phase slips in the context of mode-locked lasers, and their rate of occurrence is an important quantity.  We demonstrate below how this rate is computed approximately using the quasi-potential, and how the optimal paths provide insight into how phase slips arise.

When all coefficients are trivial, so that $c_0=c_1=c_2=d_1=d_2=D=0$,~(\ref{e:mll}) is the integrable nonlinear Schr\"odinger equation, which supports a four-parameter family of soliton solutions with functional form
\be
u_s(\xi,t) = A(t)\sech[A(t)(\xi-\Xi(t))]e^{i\Theta(t)+i\xi\Omega(t)},
\ee
where $\dot{A}=\dot{\Omega}=0$, $\dot{\Xi}=\Omega$, and $\dot{\Theta}=\frac12(A^2-\Omega^2)$.  In the perturbative limit where the coefficients are small, the dynamics of~(\ref{e:mll}) can be approximated well for short times by a diffusion process in the four parameters induced by the noise term $\eta$~\cite{iannone_nonlinear_1998}.  Related approximations made outside of the perturbative limit are not rigorously justifiable yet often agree well with numerics~\cite{kivshar_dynamics_1989}.  We employ an approach based on averaging the Lagrangian density
\be
\quad L[u_s,u_{s\xi},u_{st}] := \Im(u_{st}^*u_s) - \frac12|u_{s\xi}|^2+\frac12|u_s|^4+c_0\cos(\omega \xi)|u_s|^2.
\ee
of the variational (i.e., non-dissipative) terms in~(\ref{e:mll}) over the ansatz expressed by~(\ref{e:ansatz}).  The dissipative terms in~(\ref{e:mll}) are included through their projections against the tangent space of~(\ref{e:ansatz})~\cite{anderson_variational_2001,bale_variational_2008}, i.e.,
\begin{align}
\nabla_{y}{\cL}-&\partialderiv{}{t}\left(\nabla_{\dot y}{\cL}\right) =\nonumber\\
&2\Re\int i[-c_1u_s+c_2u_{s\xi\xi}+d_1|u_s|^2u_s-d_2|u_s|^4u_s+\eta(\xi,t)]\nabla_y{u_s^*}\,d\xi
\end{align}
where $\cL = \int L\,d\xi$ and $y=(A,\Omega,\Xi,\Theta)^T$.  Note that this Lagrangian $\cL$ and its density $L$ are associated with the variational model reduction technique used in this application, and are not related to the Lagrangian associated with the action minimization described in Sec.~\ref{s:LDT}.

\begin{figure}[t!]
\centerline{\includegraphics[width=\linewidth]{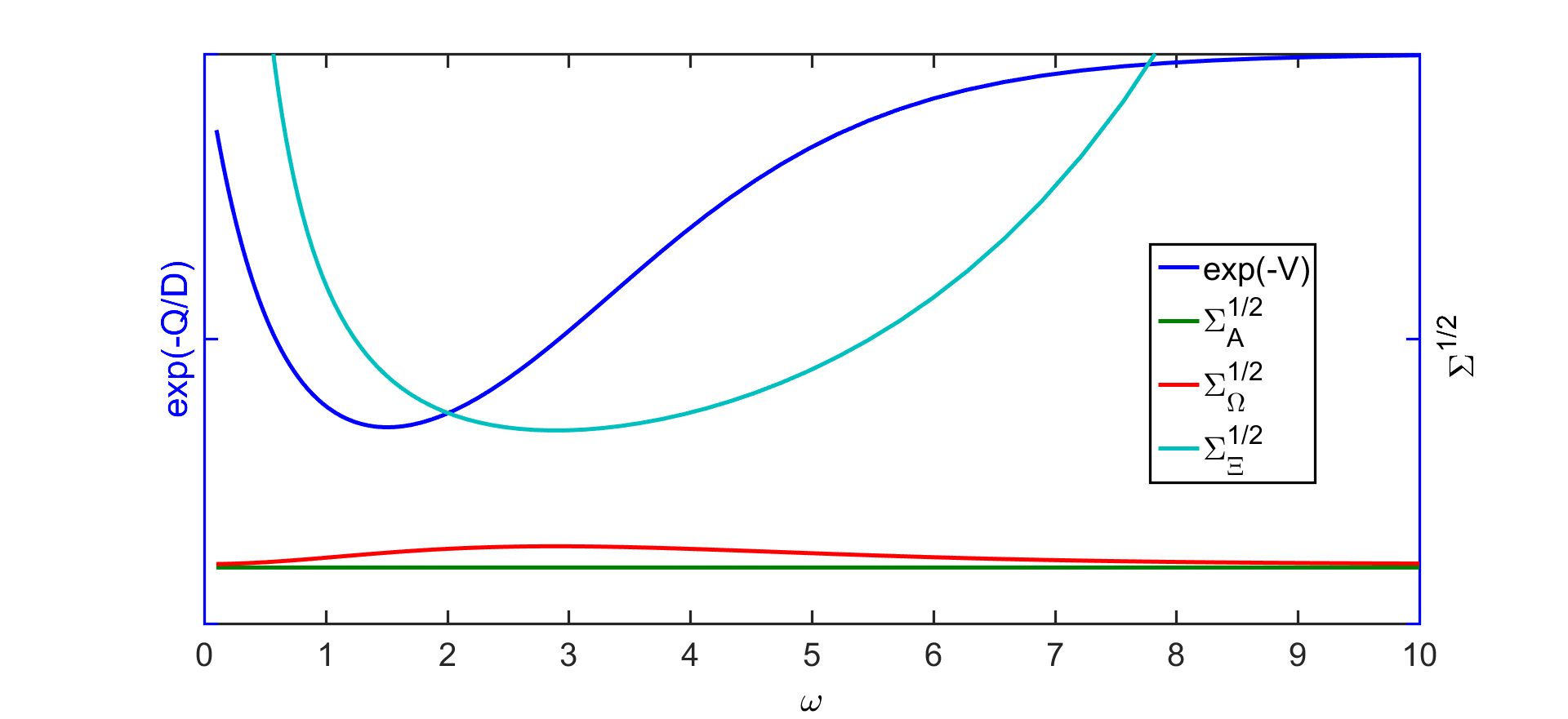}}
\caption{\label{f:uncert} Left axis: Expected phase slip rate according to large deviation theory.  Right axis: Square root of diagonal elements of covariance $\Sigma$, representing the root-mean-squared linewidth based on a conventional linearized analysis.  Minimum escape probability is achieved at a different value of $\omega$ than minimum position jitter $\Sigma_\Xi^{1/2}$.}
\end{figure}

Keeping only the first three components of $y$, i.e., $x_j=y_j$ for $j=1,\dots,3$, we have
\be
\dot{x} = f(x) + \sigma(x)\eta(t)
\label{e:mllSDE}
\ee
with drift
\be
f(x) = \bp -2c_1x_1 + (\frac43d_1-\frac23c_2)x_1^3-\frac{16}{15}d_2x_1^5 -2c_2x_1x_2^2\\
-\frac43c_2x_1^2x_2-\pi c_0\omega^2\csch(\pi\omega/2x_1)\sin(\omega x_3)/2x_1^3\\
x_2\ep
\ee
and diffusivity
\be
\sigma(x) = \bp \sqrt{x_1} & 0 & 0\\
0 & \sqrt{x_1/3} & 0\\
-x_3/\sqrt{x_1} & 0 & \sqrt{\pi^2/12x_1^3+x_3^2/x_1}\ep.
\ee
The phase dynamics for $y_4=\Theta$ are slaved to the other independent variables, with
\begin{align}
\dot\Theta &=-\frac{\pi\omega c_0}{A^3}\cos(\omega \Xi)\csch\left(\frac{\pi\omega}{2A}\right)\left(1+\frac{\pi\omega}{2A^2}\coth\left(\frac{\pi\omega}{2A}\right)\right)\nonumber\\
&-\Xi+\frac12(A^2-\Omega^2)+\left(\frac{\sqrt{12+\pi^2}}{6\sqrt{A}}\right)\eta_4.
\label{e:phaseSDE}
\end{align}
The noise term $\eta(t)\in{\mathbb R}^3$ in~(\ref{e:mllSDE}) satisfies
 \be
{\mathbb E}[\eta(t)\eta(t')^T] = DI\delta(t-t')
\ee
and is obtained from $\eta(\xi,t)$ in~(\ref{e:mll}) using the same averaged Lagrangian method, and essentially represents that portion of the stochastic driving term that directly influences the parameters in $x$.  Matrix $I$ is the $3\times 3$ identity.

\begin{figure}[t!]
\begin{minipage}{0.49\linewidth}
\begin{center}
\includegraphics[scale=0.32]{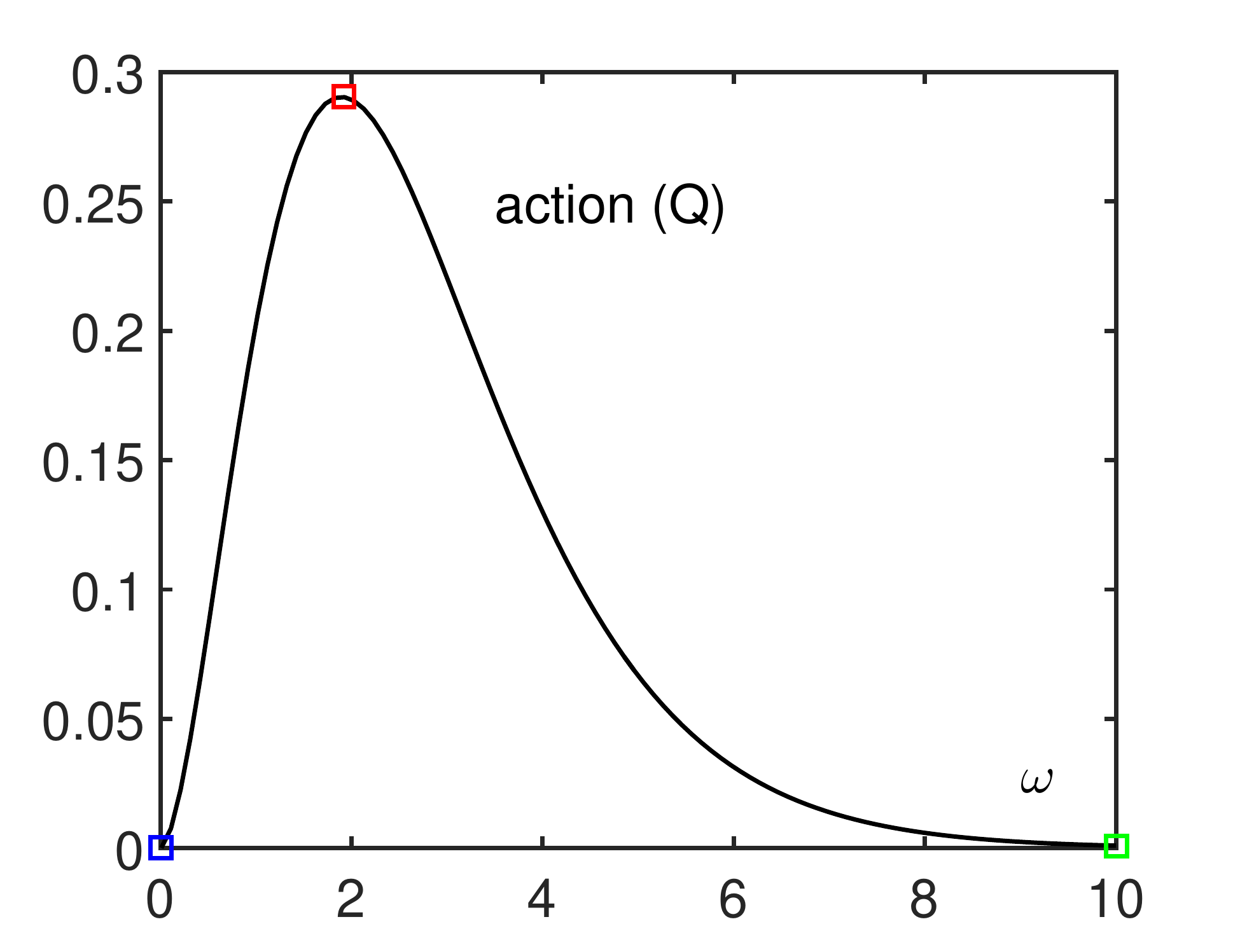}
\end{center}
\end{minipage}
\begin{minipage}{0.49\linewidth}
\begin{center}
\includegraphics[scale=0.32]{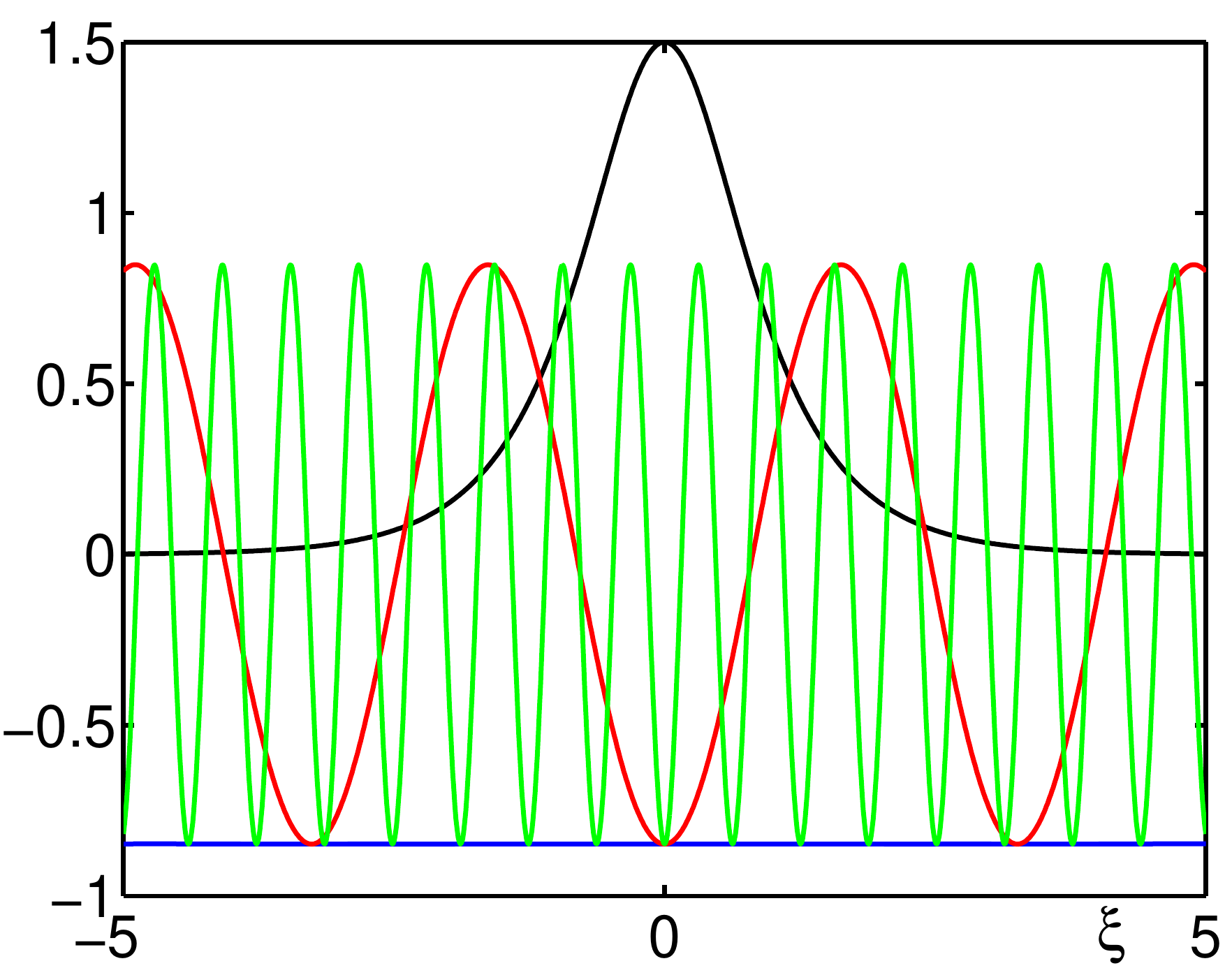}
\end{center}
\end{minipage}
\caption{\label{f:varyingOmega}(a) Minimum action connecting stable fixed
  point at origin ($n=0$) to saddle ($n=1$), plotted vs.\ parameter $\omega$
  from~(\ref{e:mll}).  (b) Laser pulse (black) and active feedback for
  $\omega=0$ (blue), $\omega=1.9$ (red), and $\omega=10$ (green),
  corresponding to the squares in (a).}
\end{figure}

In the deterministic ($D=0$) limit, the nontrivial equilibria are $x^{(n)} = (A_0,\; 0,\; n\pi/\omega)$ for $n\in{\mathbb Z}$, with
\[
A_{0\pm}^2 = \frac{5}{16d_2}\left(2d_1-c_2\pm\sqrt{(2d_1-c_2)^2-\frac{96}5c_1d_2}\right).
\]
Of these, the equilibria with $A_{0+}$ are stable nodes or spirals for $n$ even and saddles for $n$ odd.  The equilibria with
$A_{0-}$ are always saddles.
To repeat the analysis from Ref.~\cite{bao_pulse_2014}, linearization about $u_0$ gives steady state covariance
\be
\Sigma = \frac{D}{2\pi}\int ds \left[(is I-M)^{-1}\sigma(x^0)\right]\left[(is I-M)^{-1}\sigma(x^0)\right]^{\dag},
\ee
where
\be
M = \bp 8c_1-\frac43(2d_1-c_2)A_0^2 & 0 & 0\\
0 & -\frac43c_2A_0^2 & -\frac{(-1)^n\pi c_0\omega^3}{2A_0^3}\csch(\frac{\pi\omega}{2A_0})\\
0 & 1 & 0\ep.
\ee

As discussed above, however, the feedback mechanism requires that we consider two sources of uncertainty in the laser output, the first associated with laser linewidth and the second associated with the expected rate of phase slips.
A complete analysis of the phase slip rate requires a consideration of the phase dynamics expressed by ~(\ref{e:mllSDE}) and~(\ref{e:phaseSDE}), including contributions from the dynamics of all other parameters.  For simplicity of presentation, we consider here transitions between the equilibria described above, noting that the drift term in~(\ref{e:phaseSDE}) changes by $2\pi/\omega$ between neighboring stable equilibria.  In appropriate parameter regimes, this is the dominant contributing factor to phase slips.

As described in Sec.~\ref{s:escapeTime}, the expected time for transition between neighboring equilibria is
$\exp(Q(x^0,x^2)/D)$, where $Q$ is the quasi-potential.
Figure~\ref{f:uncert} plots these distinct types of uncertainty against one of the active feedback parameters, $\omega$, where the reciprocal of the transition time has been expressed as a rate.  It is clear from the figure that one cannot in general expect to be able to minimize both linewidth and phase slip rate simultaneously; practical considerations inform the balance that must be struck between these sources of uncertainty~\cite{moore_trade-off_2014}.

Although the action $Q$ provided by the GMAM computation is the quantity that determines the mean exit rate estimate, it is also instructive to examine the paths associated with minimum action for various choices of physical constants.
\begin{figure}[t!]
\begin{minipage}{0.32\linewidth}
\begin{center}
\includegraphics[scale=.23]{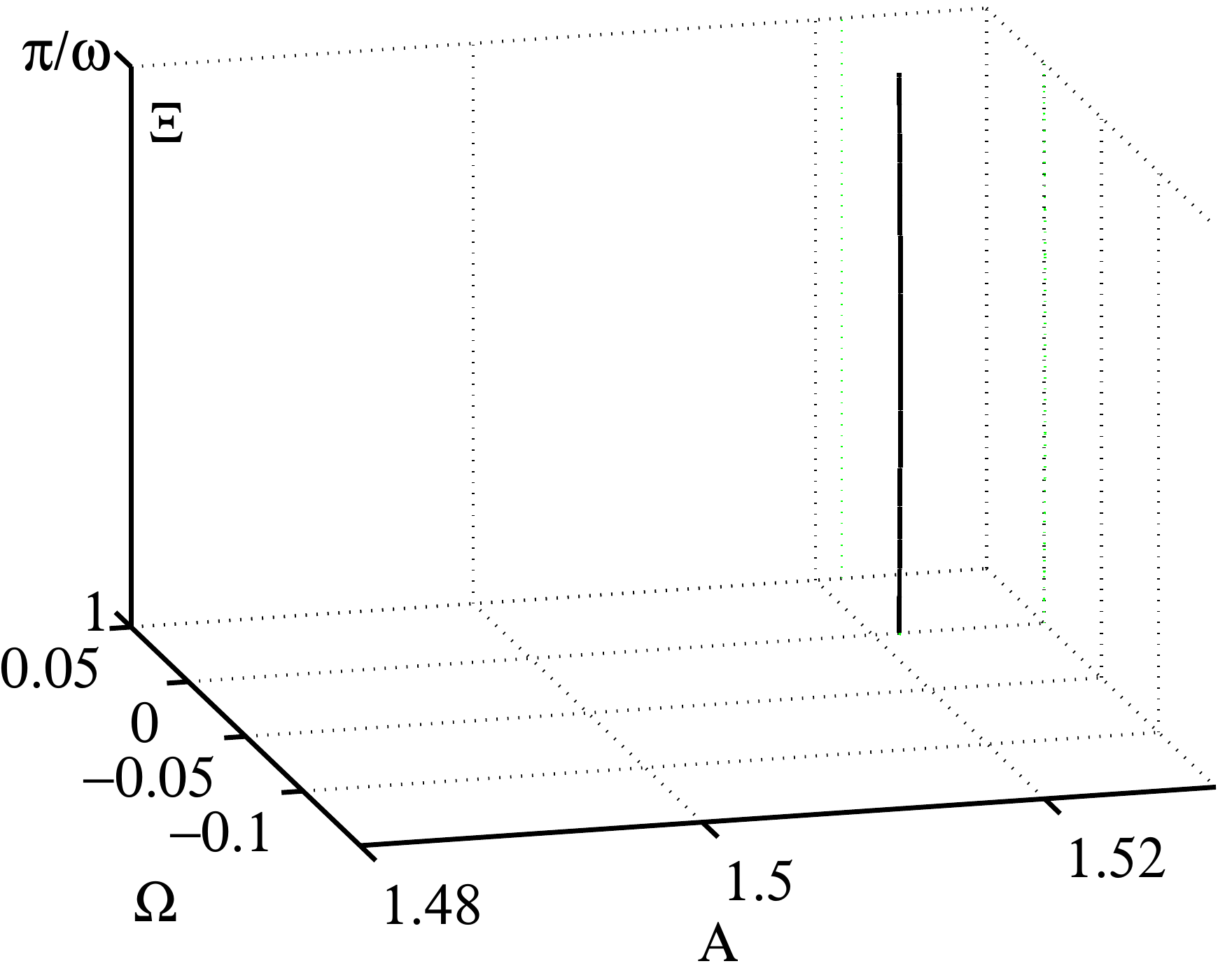}
\end{center}
\end{minipage}
\begin{minipage}{0.32\linewidth}
\begin{center}
\includegraphics[scale=.23]{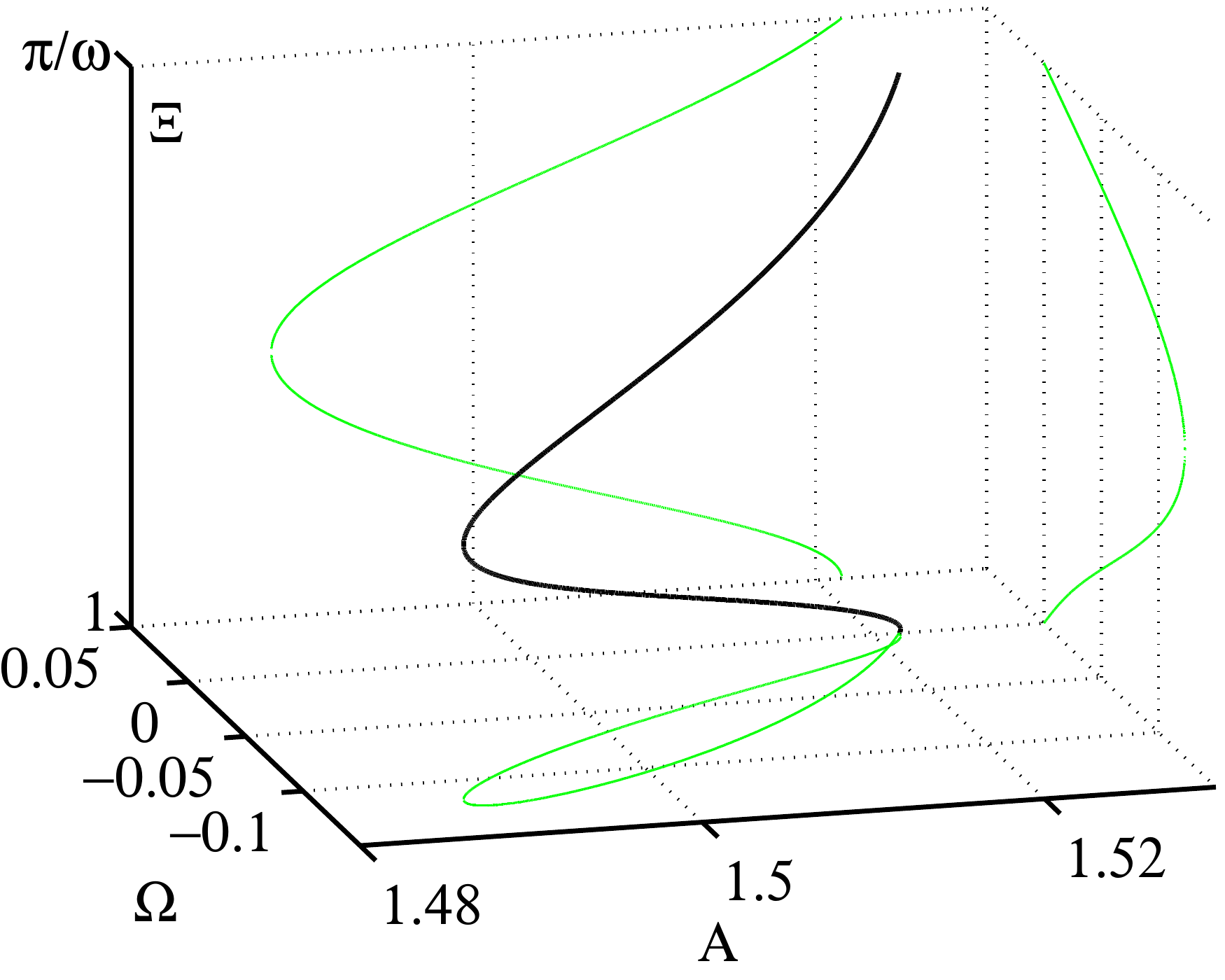}
\end{center}
\end{minipage}
\begin{minipage}{0.32\linewidth}
\begin{center}
\includegraphics[scale=.23]{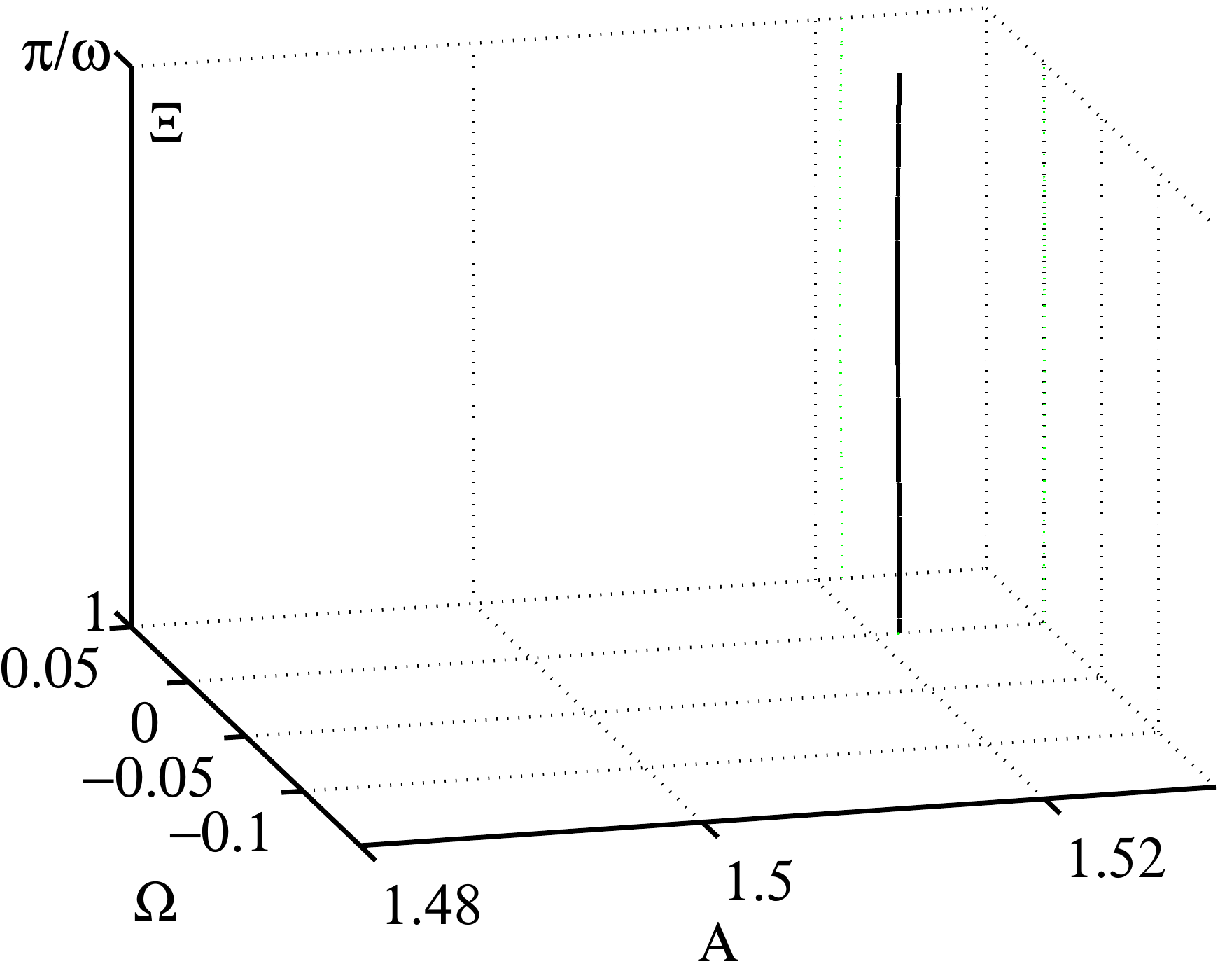}
\end{center}
\end{minipage}
\caption{\label{f:omegaPaths} Minimum action paths from stable fixed point ($n=0$) to saddle ($n=1$) for (a) $\omega=0$, (b) $\omega=1.9$, and (c) $\omega=10$, corresponding to squares in Fig.~\ref{f:varyingOmega}(a).}
\end{figure}

\begin{figure}[b!]
\begin{minipage}{0.49\linewidth}
\begin{center}
\includegraphics[scale=0.3]{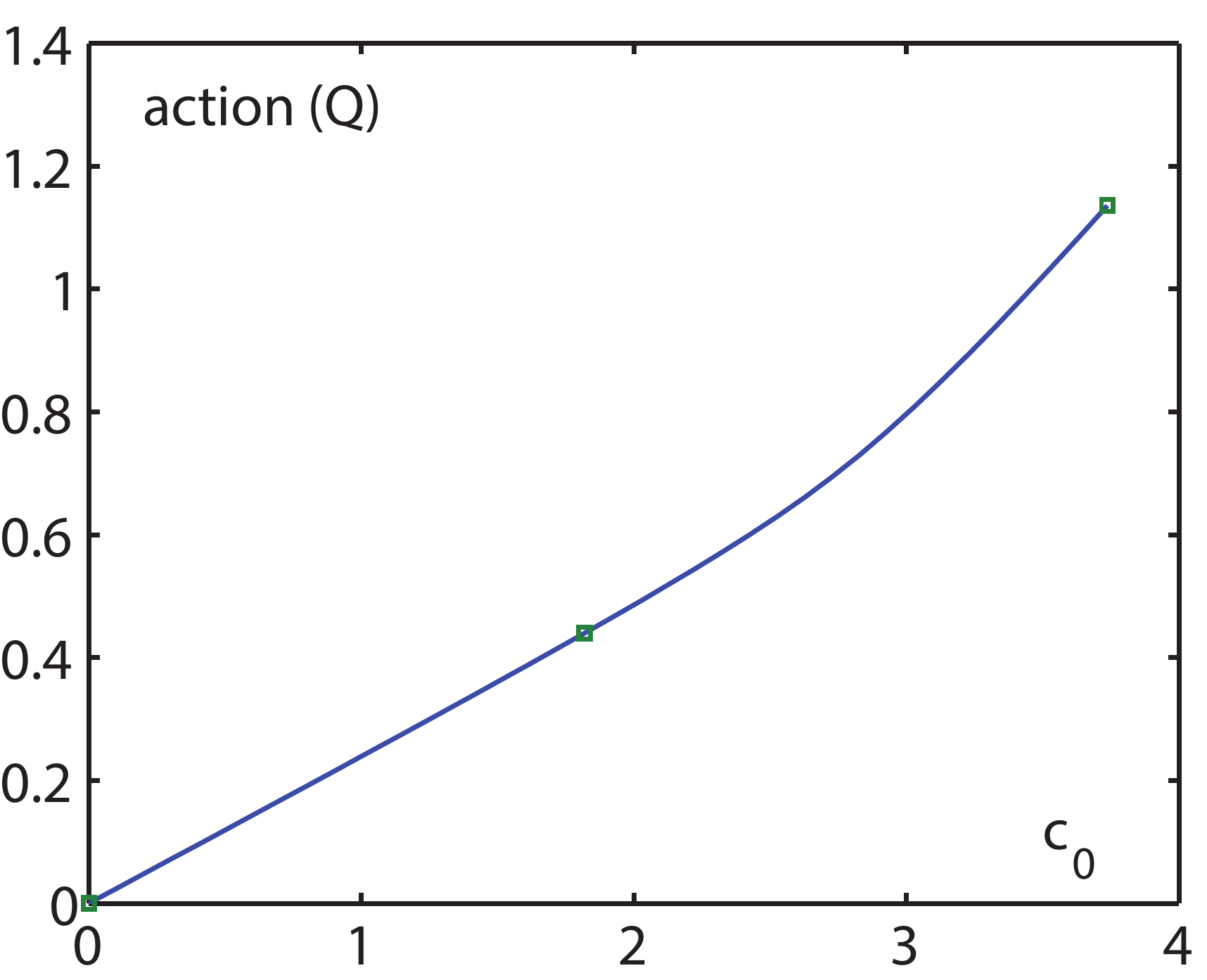}
\end{center}
\end{minipage}
\begin{minipage}{0.49\linewidth}
\begin{center}
\includegraphics[scale=0.3]{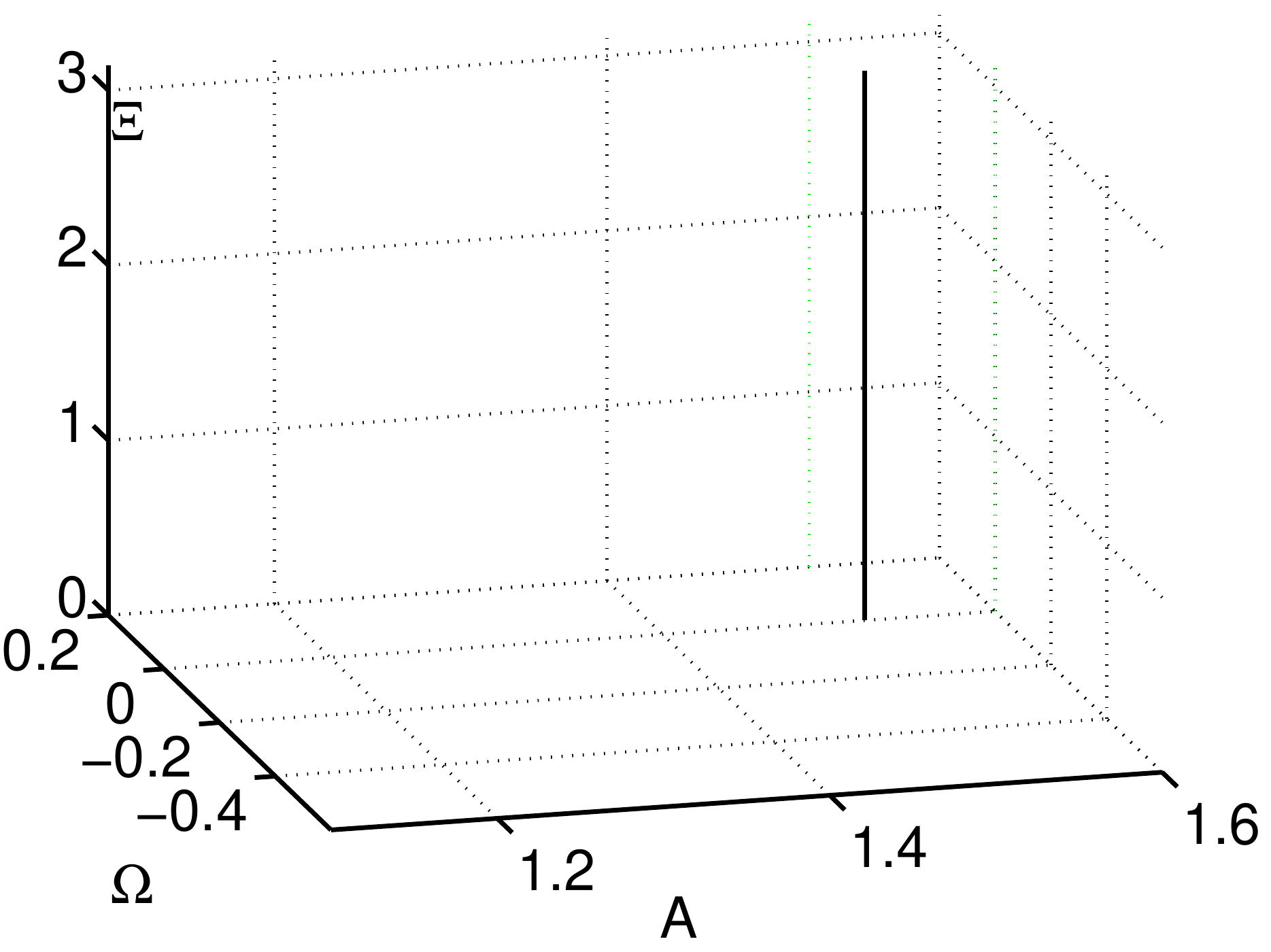}
\end{center}
\end{minipage}
\begin{minipage}{0.49\linewidth}
\begin{center}
\includegraphics[scale=.3]{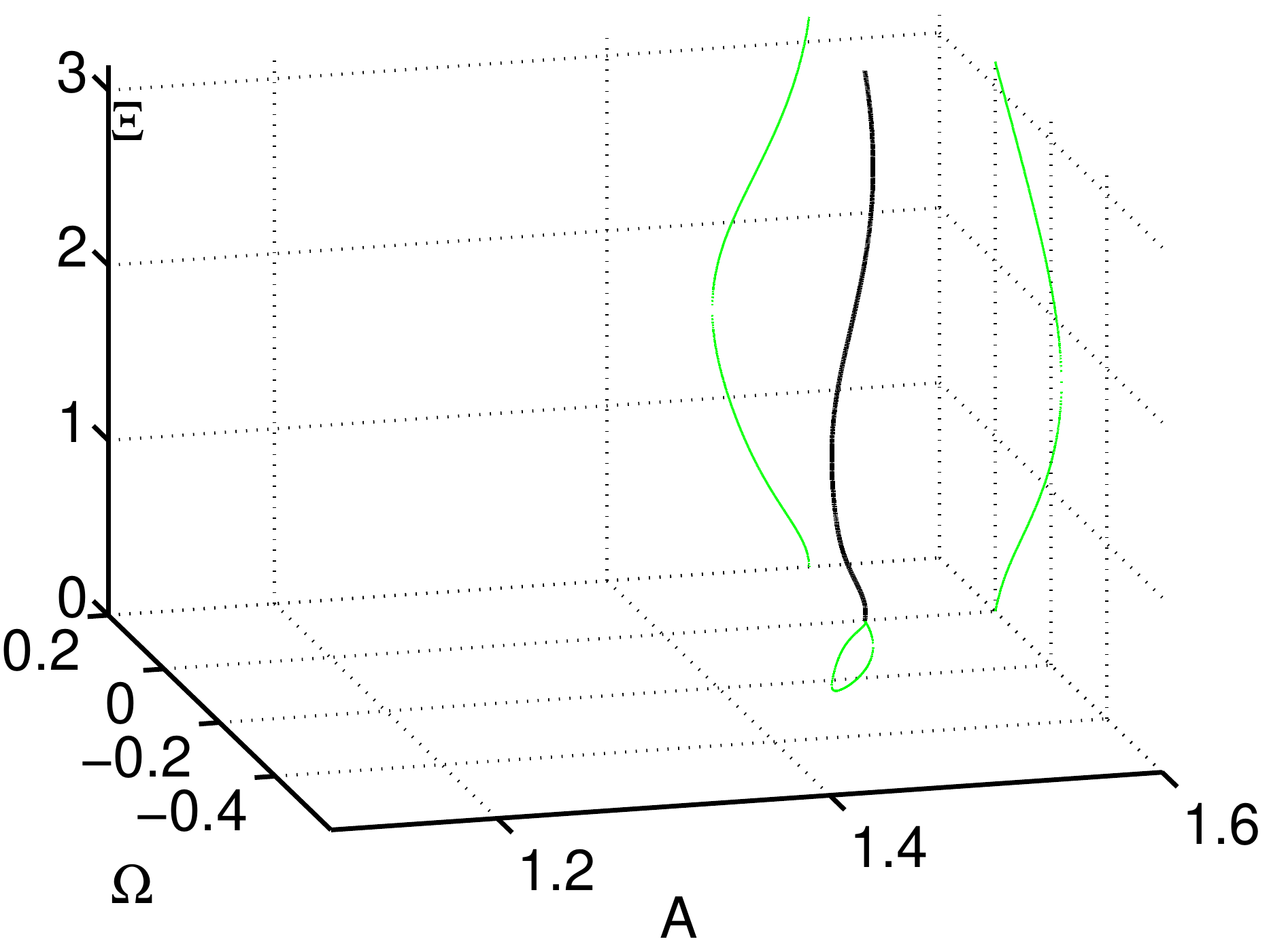}
\end{center}
\end{minipage}
\begin{minipage}{0.49\linewidth}
\begin{center}
\includegraphics[scale=.3]{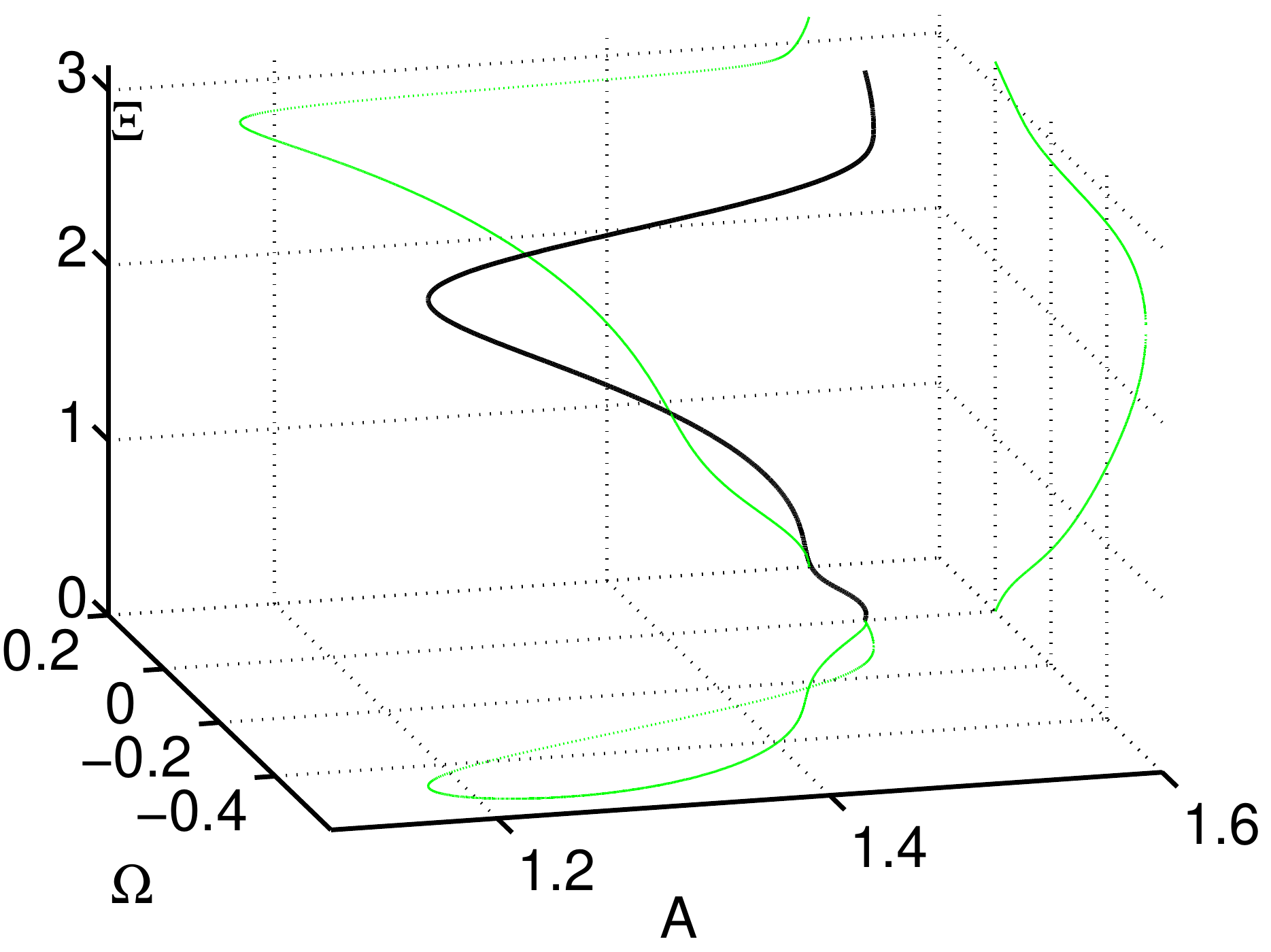}
\end{center}
\end{minipage}
\caption{\label{f:varyingc0}(a) Minimum action path connecting stable fixed point at origin ($n=0$) to saddle ($n=1$), plotted vs.\ parameter $c_0$ from~(\ref{e:mll}).  Minimum action paths for squares on (a) are plotted in (b), (c), and (d).}
\end{figure}

Figure~\ref{f:varyingOmega} graphs the minimum action as a function of the inverse width of the active feedback trap, represented by $\omega$ in~(\ref{e:mll}).  For small values of $\omega$, the trap is so wide as to be essentially constant over the support of the pulse.  This leaves~(\ref{e:mll}) effectively invariant with respect to translations of the pulse offering no restoring force to noise that moves the pulse's position from the origin.  The action for a change in position along this neutrally stable manifold is therefore zero.  This is reflected in Fig.~\ref{f:omegaPaths}(a), where the most likely transition path involves simple translation in $\xi$.  Conversely, for large values of $\omega$, as seen in the green curve, the trap's wavelength is much shorter than the pulse's width, so that its ability to stop the pulse from moving is limited.  Neighboring equilibria are also very close together in this limit.  The action therefore approaches zero and the most likely transition path between equilibria is again a linear translation, as shown in Fig.~\ref{f:omegaPaths}(c).  Between these values of $\omega$, however, is a range of values where variations of the trap are commensurate in scale with the support of the soliton, such that the trap offers maximal resistance to large excursions and generates a large value of action.  Figure~\ref{f:omegaPaths}(b) shows that this increased action is associated with an exit path that exploits the internal dynamics of the SDE much more significantly.  This picture is qualitatively similar to the impact of varying $\omega$ on the standard deviation, or jitter, in pulse position, but the minimum is achieved at a different value of $\omega$, as illustrated in Fig.~\ref{f:uncert}.

Figure~\ref{f:varyingc0} illustrates a different dependence of action on parameter value $c_0$, the amplitude of the active feedback trap for fixed $\omega=1$.  In this case, the action increases monotonically with $c_0$, and the minimum action paths are seen to become more and more tortuous as the trap becomes more and more effective at blocking the transition mechanism of simple translation.
\subsection{Importance sampling for large soliton walks}
\label{s:IS}

\begin{figure}[t!]
        \centerline{\includegraphics[width=\textwidth]{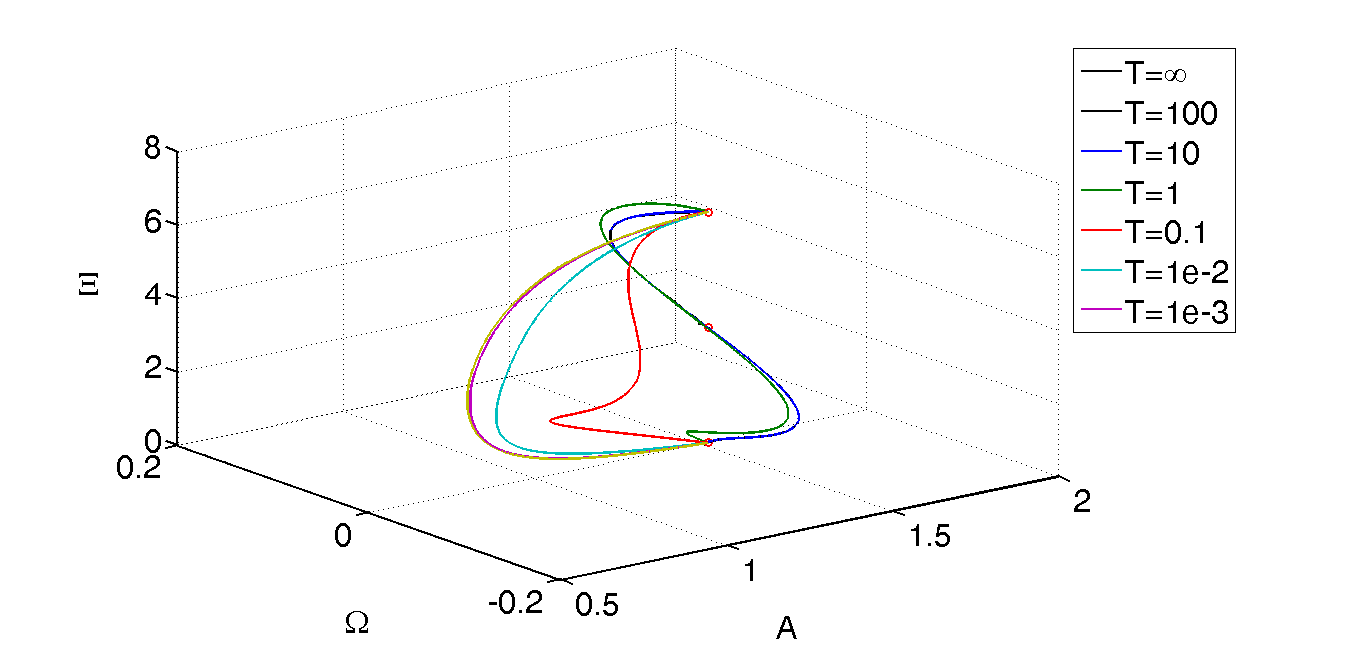}}
        \caption{\label{f:pathsVaryingT}Minimum action paths for the SDE given
          by~(\ref{e:mllSDE}) as the horizon time $\tf$ is varied.}
\end{figure}

In the laser dynamics example in Sec.~\ref{s:phaseSlips}, the minimum action calculated for passage from a stable fixed point to a neighboring saddle is used to approximate the average rate at which noise is expected to drive the laser from its desired operating point.  The paths associated with these minimum action computations are informative in their own right since they illustrate the dynamic paths taken on these most probable excursions.  In fact, these paths provide the information necessary to implement importance sampling in Monte Carlo simulations designed to compute correct expectations of exit probabilities in finite time as described in Sec.~\ref{s:exitSampling}.  Figure~\ref{f:pathsVaryingT} demonstrates that these paths change significantly as the time horizon for an exit is varied, from a geodesic path with respect to a norm dictated by the diffusion tensor to the limiting path corresponding to the quasi-potential calculations from Sec.~\ref{s:phaseSlips}.

For exit probabilities on finite times, traditional large deviation theory provides an asymptotic scaling law,
\be
\lim_{D\rightarrow 0}D\ln P\left(x(\tf)\notin \Omega\right) = -\inf_{x(\tf)\notin \Omega}{\mathcal{S}}_{\tf}[x],
\ee
where the probability $P$ contains an undetermined prefactor.  For the reduced
model provided by~(\ref{e:mllSDE}), computing this prefactor is at least
possible using more substantial analysis of the Fokker-Planck
equation~\cite{matkowsky_exit_1977}.  Given that the model reduction from the
stochastic PDE~(\ref{e:mll}) is itself imperfect, one would ideally seek to
perform Monte Carlo simulations on the original model to verify the accuracy
of estimates obtained through the reduced model.  To compute exit
probabilities in the context of small noise, a naive Monte Carlo approach is
computationally prohibitive, requiring the use of importance-sampled Monte Carlo simulations based on the paths computed above~\cite{moore_method_2008,moore_soliton_2005,spiller_computing_2005}.

\begin{figure}[t!]
    \centerline{\includegraphics[width=\textwidth]{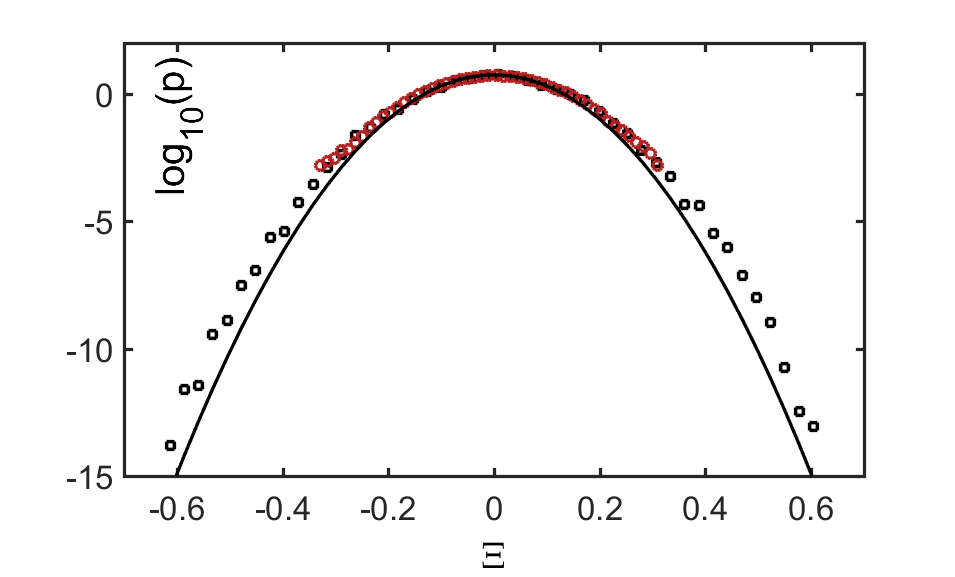}}
    \caption{\label{f:NLStim_eps} Probability density $p$ associated with the
      diffusion of position $\Xi$ obtained from importance-sampled Monte Carlo
      simulations of the stochastic PDE~(\ref{e:filter}) (black squares) and
      from unbiased simulations (brown circles).  The black curve is the Gaussian
      probability density function predicted by linearizing the
      reduced ODE given by~(\ref{e:filterSDE})-(\ref{e:diffusivity}).}
\end{figure}

For example, a simplified version of~(\ref{e:mll}) relevant to optical communications includes spectral filtering, compensatory gain ($c_1<0$), and white noise:
\begin{align}
i\partialderiv{u}{t}+&\frac12\partialderiv[2]{u}{\xi}+|u|^2u = -ic_1u +ic_2\partialderiv[2]{u}{\xi}  + i\eta(\xi,t).
\label{e:filter}
\end{align}
Filters have been used as a mechanism for reducing jitter in the position of solitons by adding a restoring force to the spatial frequency, which is tied to the soliton's speed through the dispersion relation~\cite{iannone_nonlinear_1998}.  Damping fluctuations in spatial frequency is therefore an effective way to suppress position jitter, as illustrated in Fig.~\ref{f:NLStim_eps}.
The same reduction method used above produces the SDE given by
\be
\dot{x} = f(x) + \sigma(x)\eta
\label{e:filterSDE}
\ee
with $x = (A,\,\Omega,\,\Xi)^T$, where the drift is
\be
f(x) = \bp -2c_1x_1 -\frac23c_2x_1^3-2c_2x_1x_2^2\\
-\frac43c_2x_1^2x_2\\
x_2
\ep,
\ee
and diffusivity is
\be
\sigma(x) = \bp \sqrt{x_1} & 0 & 0 \\
0 & \sqrt{x_1/3} & 0 \\
-x_3/\sqrt{x_1} & 0 & \sqrt{\pi^2/12x_1^3+x_3^2/x_1}\\
\ep,
\label{e:diffusivity}
\ee
and where we have again omitted the dynamics in $\Theta$.
The GMAM method presented in Sec.~\ref{s:escapeTime} can be applied to compute
the optimal path from a pulse centered at the origin with amplitude,
frequency, and position coordinates $(A,\Omega,\Xi)=(1,0,0)$ to a final
position $\Xi(\tf) = \Xi_f$ with the other two final coordinates left
unspecified.  By picking different values of $\Xi_f$, one can reconstruct a
histogram of probabilities from which the entire probability density function
(pdf) for the final position $\Xi$ can be computed well down into the low-probability tails.
Figure~\ref{f:NLStim_eps} illustrates such a reconstructed pdf, where 50,000
importance-sampled Monte Carlo simulations of the stochastic
PDE~(\ref{e:filter}) have been used to compute pdf tail probabilities on the
order of $10^{-12}$ and below.  Also shown is the Gaussian pdf predicted by a linearization of the ODE reduction given by~(\ref{e:filterSDE})-(\ref{e:diffusivity}).  This curve and simulations of the nonlinear ODEs (not shown) demonstrate that the low-order reduction fails to accurately compute probabilities in the tails of the pdf, but Fig.~\ref{f:NLStim_eps} clearly shows that it is effective at generating appropriately biased simulations of the original PDE model.

\subsection{Extinction in biological populations}
\label{s:extinction}

We now consider internal noise and how it may induce the rare event of extinction of a species or of an infectious disease. Extinction can be good (disease) or bad (species); either way, it is important to understand why and how often it can be expected to happen for different model parameters.  We consider two examples: (1) the Susceptible-Infectious-Susceptible (SIS) epidemic model, and (2) the Allee effect model.

\begin{figure}[t!]
\begin{center}
\begin{tikzpicture}[ultra thick]
\filldraw[black](0,0) circle (0.075) node[below=0.5cm,align=center] {$x_0$}
node[below=1.0cm,align=center] {Extinct\\ state};
\filldraw[black](5,0) circle (0.075) node[below=0.5cm,align=center] {$x_1$}
node[below=1.0cm,align=center] {Endemic\\ state};
\draw[->] (0,0) to (-1,0)  node[left=0.5cm] {(a)};
\draw[middlearrow={<}] (5,0) to (7,0);
\draw[middlearrow={>}] (0,0) to (2,0);
\draw[middlearrow={>}] (3,0) to (5,0);
\draw (2,0) to (3,0);
\end{tikzpicture}\\
\vspace{2cm}
\begin{tikzpicture}[ultra thick]
\filldraw[black](0.5,0) circle (0.075) node[below=0.5cm,align=center] {$x_0$}
node[below=1.0cm,align=center] {Extinct\\ state};
\filldraw[black](5.5,0) circle (0.075) node[below=0.5cm,align=center] {$x_2$}
node[below=1.0cm,align=center] {Carrying\\ capacity};
\filldraw[black](3,0) circle (0.075) node[below=0.5cm,align=center] {$x_1$}
node[below=1.0cm,align=center] {Allee\\ threshold};
\draw[middlearrow={<}] (-0.5,0) to (-1,0)  node[left=0.5cm] {(b)} ;
\draw (-0.5,0) to (0.5,0);
\draw[middlearrow={<}] (0.5,0) to (3,0);
\draw[middlearrow={>}] (3,0) to (5.5,0);
\draw[middlearrow={<}] (6.5,0) to (7,0);
\draw (5.5,0) to (6.5,0);
\end{tikzpicture}
\end{center}
\caption{\label{fig:ext-scenario}(a) Topology of the deterministic SIS model.
There are two fixed points: a stable endemic
   state and an unstable extinct state. (b) Topology of the deterministic
   Allee effect model. There are three fixed points: a stable
   carrying capacity, a stable extinct state, and an unstable Allee threshold.}
\end{figure}
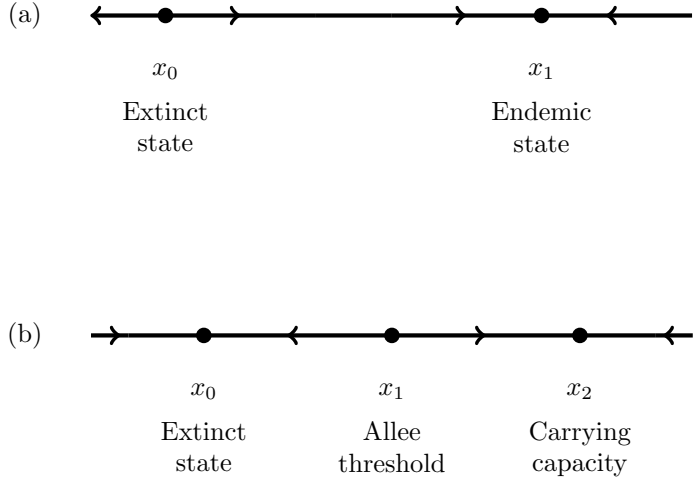

The deterministic SIS model has two fixed points. One is an extinct state where no
infectious individuals are present, while the other is an endemic state where the infection
is maintained. The stability of these two fixed points is determined by the
value of the reproductive number $R_0$. The reproductive number can be thought
of as the average number of new infectious individuals that one infectious
individual generates over the course of the infectious period, in an entirely
susceptible population. For $R_0 > 1$, the extinct state is unstable while the
endemic state is stable, as shown in Fig.~\ref{fig:ext-scenario}(a). Note that
since the model is deterministic, a population at the attracting endemic state
can never go extinct.

The deterministic Allee effect model has three fixed points. One can see in
Fig.~\ref{fig:ext-scenario}(b) that the extinct state $x_0$ is stable, as is
the carrying capacity $x_2$. The Allee threshold $x_1$ is unstable, so when
initial values lie between $x_1$ and $x_2$, the deterministic solution will
increase to $x_2$, while for initial values less than $x_1$ the deterministic
solution decreases to the extinct state. Similar to the deterministic SIS, a
population at the attracting carrying capacity can never go extinct.

\begin{figure}[t!]
\begin{center}
\includegraphics[scale=0.5]{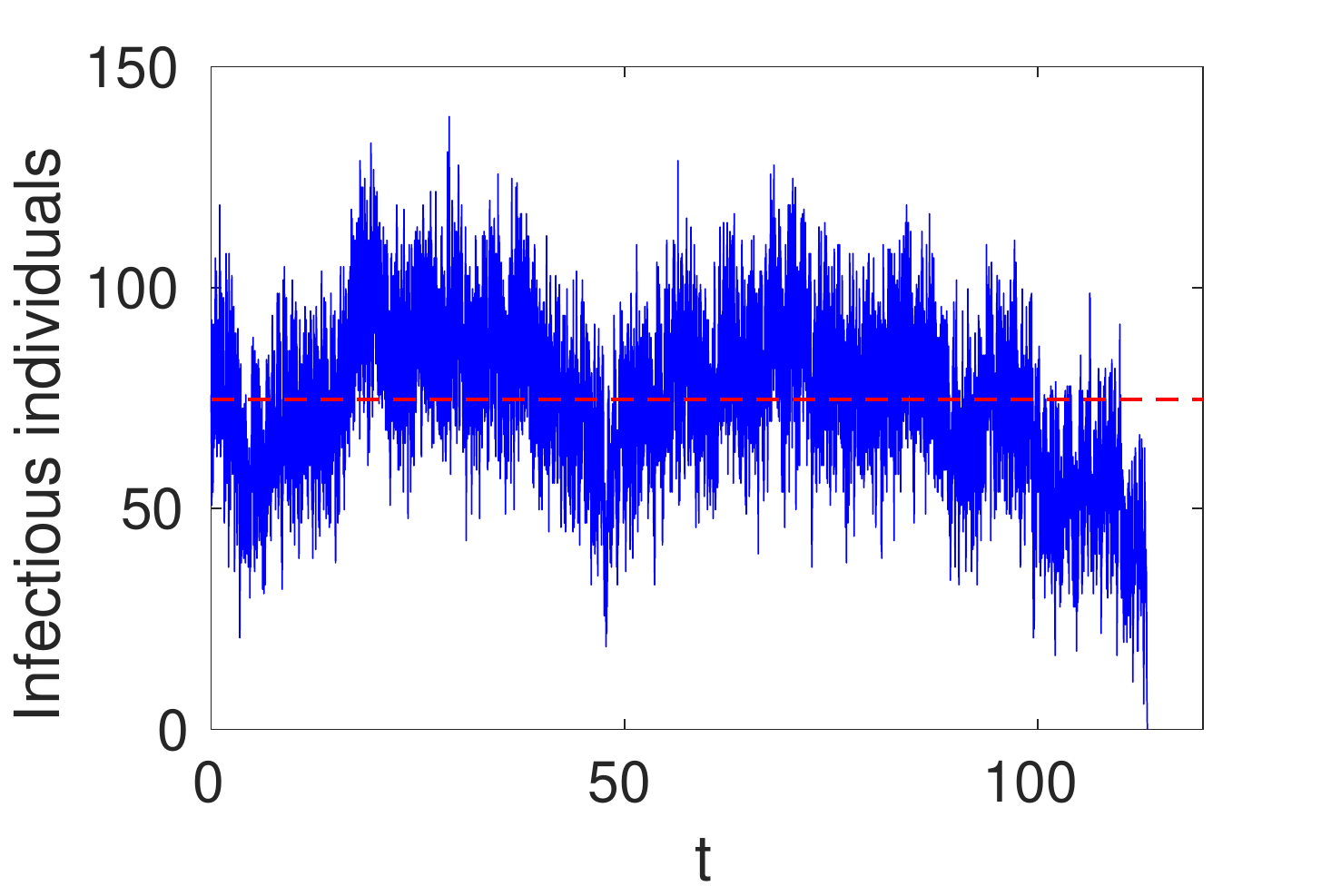}
\end{center}
\caption{\label{fig:SISa} Stochastic realization of an SIS system showing
  extinction after fluctuating about the endemic state for a long time.}
\end{figure}

To capture extinction events, we must consider a stochastic model with
internal noise that represents the random interactions of individuals in the
population. Figure~\ref{fig:SISa} shows a time series of infectious individuals
for a stochastic SIS model. One can see fluctuations about the endemic state
for a long period of time until eventually the noise-induced rare event
(extinction) occurs. Similarly, Fig.~\ref{fig:int-switch} shows a time series
of individuals for a stochastic Allee effect model. One sees the population fluctuating about the carrying
capacity for a long period of time until the rare extinction event
occurs. Details of the Monte Carlo method used in both simulations can be
found in Appendix~\ref{s:appA}.

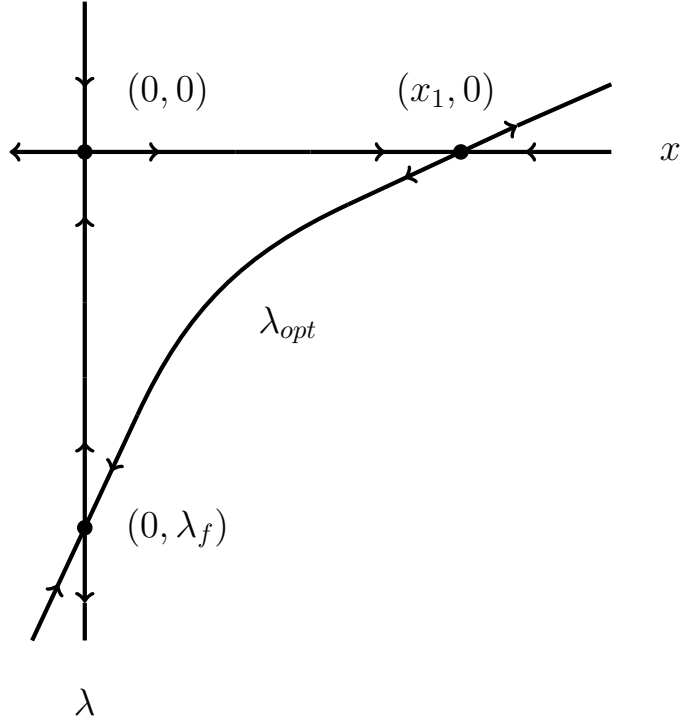
\begin{figure}[t!]
\begin{center}
\begin{tikzpicture}[ultra thick]
\filldraw[black](0,0) circle (0.075);
\filldraw[black](5,0) circle (0.075) node[above=0.4cm,xshift=-0.2cm] {\Large $(x_1,0)$};
\draw[->] (0,0) to (-1,0) ;
\draw[middlearrow={<}] (5,0) to (7,0) node[right=0.5cm] {\Large $x$};
\draw[middlearrow={>}] (0,0) to (2,0) node[above=0.4cm,xshift=-0.9cm] {\Large $(0,0)$};
\draw[middlearrow={>}] (3,0) to (5,0);
\draw (2,0) to (3,0);
\draw[middlearrow={<}] (0,0) to (0,2);
\draw[middlearrow={<}] (0,0) to (0,-2);
\filldraw[black](0,-5) circle (0.075) node[right=0.4cm] {\Large $(0,\lambda_f)$};
\draw[->] (0,-5) to (0,-6);
\draw[middlearrow={<}] (0,-3) to (0,-5);
\draw (0,-2) to (0,-3);
\draw[middlearrow={>}] (0.7,-3.5) to (0,-5);
\draw[middlearrow={>}] (-0.7,-6.5) to (0,-5);
\draw[middlearrow={>}] (5,0) to (3.5,-0.7);
\draw (3.5,-0.7) to[out=205,in=65] (0.7,-3.5) node[xshift=2cm,yshift=1.2cm] {\Large $\lambda_{opt}$};
\draw[->] (5,0) to (5.75,0.35);
\draw (0,-6) to (0,-6.5) node[below=0.5cm] {\Large $\lambda$};
\draw (5.75,0.35) to (7,0.9);
\end{tikzpicture}
\end{center}
\caption{\label{fig:OPtopoA} Steady states of Hamilton's equations~(\ref{e:Hameqnsi}) and
  zero-energy trajectories of the Hamiltonian~(\ref{e:Hamil}) for the SIS model as well as
  other models with deterministic mean-field topology of the type shown
  in~\ref{fig:ext-scenario}(a). The
extinct state $x_0=0$ is a repelling point of the deterministic mean-field
equation, while the endemic state $x_1$ is an attracting point. In this
extinction scenario, the optimal path to extinction $\lambda_{opt}$ is a heteroclinic
trajectory with non-zero momentum that connects the equilibrium point
$(x,\lambda)=(x_1,0)$ with a fluctuational extinct state
$(x,\lambda)=(0,\lambda_f)$.}
\end{figure}

The two examples considered in this section have different mean-field
topologies as seen in Fig.~\ref{fig:ext-scenario}. For the SIS model, the
extinct state $x_0$ is a repelling point of the deterministic mean-field
equation, while the endemic state $x_1$ is an attracting point. In this
extinction scenario, the optimal path to extinction is a heteroclinic
trajectory with non-zero momentum that connects the equilibrium point
$(x,\lambda)=(x_1,0)$ with a fluctuational extinct state
$(x,\lambda)=(0,\lambda_f)$. This new fluctuational extinct state is an
equilibrium point of Hamilton's equations given
by~(\ref{e:Hameqnsi}). Figure~\ref{fig:OPtopoA} shows the optimal path topology in
the expanded 2D space.

\begin{figure}[t!]
\begin{center}
\begin{tikzpicture}[ultra thick]
\filldraw[black](0.5,0) circle (0.075);
\filldraw[black](7,0) circle (0.075);
\filldraw[black](2,0) circle (0.075);
\draw[middlearrow={>}] (-1,0) to (0.5,0) node[above=0.3cm, xshift=-0.7cm] {\Large $(0,0)$};
\draw[middlearrow={<}] (0.5,0) to (2,0) node[above=0.3cm, xshift=0.7cm] {\Large $(x_1,0)$};
\draw (3.5,0) to (5.5,0);
\draw[middlearrow={>}] (2,0) to (3.5,0);
\draw[middlearrow={<}] (7,0) to (8.5,0) node[right=0.5cm] {\Large $x$};
\draw[middlearrow={>}] (5.5,0) to (7,0) node[above=0.3cm, xshift=-0.7cm] {\Large $(x_2,0)$};
\draw[->] (7,0) to (6.8,-.6);
\draw[->] (7,0) to (7.18,0.6);
\draw[middlearrow={>}] (2.55,-0.86) to (2,0);
\draw (6.8,-0.6) to[out=-110,in=-60,distance=4cm] (2.55,-0.86) node[xshift=4cm,yshift=-2.5cm] {\Large $\lambda_{opt}$};
\draw[middlearrow={>}] (1.5,0.8) to (2,0);
\draw (0.5,1) to (0.5,-4) node[below=0.5cm] {\Large $\lambda$};
\end{tikzpicture}
\end{center}
\caption{\label{fig:OPtopoB} Steady states of Hamilton's equations~(\ref{e:Hameqnsi}) and
  zero-energy trajectories of the Hamiltonian~(\ref{e:Hamil}) for the Allee effect model as well as
  other models with deterministic mean-field topology of the type shown
  in~\ref{fig:ext-scenario}(b). The extinct
state $x_0=0$ is an attracting point of the deterministic mean-field equation.
Additionally, there is an intermediate repelling point $x_1$ that lies between the extinct state and another attracting state $x_2$.
In this extinction scenario, the optimal
path to extinction, is composed of two segments.  The first segment is a
heteroclinic trajectory with non-zero momentum that connects the equilibrium
point $(x,\lambda)=(x_2,0)$, with the intermediate equilibrium point
$(x,\lambda)=(x_1,0)$.  The second segment consists of the
segment along $\lambda=0$ from $x_1$ to the extinct state $x_0$.}
\end{figure}
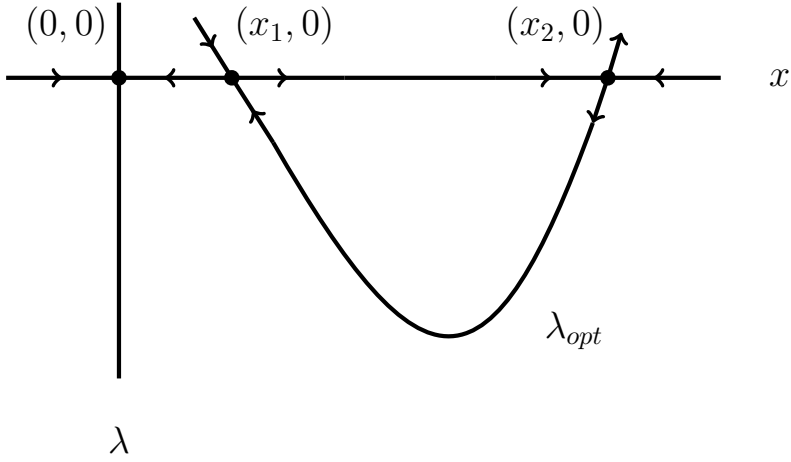

For the Allee effect model, the extinct
state $x_0$ is an attracting point of the deterministic mean-field equation.
Additionally, there is an intermediate repelling point $x_1$ that lies between the extinct state and another attracting state $x_2$.
In this extinction scenario, the optimal
path to extinction, is composed of two segments.  The first segment is a
heteroclinic trajectory with non-zero momentum that connects the equilibrium
point $(x,\lambda)=(x_2,0)$, with the intermediate equilibrium point
$(x,\lambda)=(x_1,0)$.  The second segment consists of the
segment along $\lambda=0$ from $x_1$ to the extinct state $x_0$.  Figure~\ref{fig:OPtopoB} shows the optimal path topology in
the expanded 2D space.

\subsubsection{SIS epidemic model}
\label{s:SISmodel}

The SIS problem has mean-field equations given as
\begin{subequations}
\begin{flalign}
\dot{S}&=\mu K-\frac{\beta}{K}SI+\gamma I -\mu S,\\
\dot{I}&=\frac{\beta}{K}SI-\gamma I -\mu I,
\end{flalign}
\end{subequations}
where $\beta$ is the contact rate, $\gamma$ is the recovery rate, and $\mu$ is
the birth/death rate. Using a constant population assumption $S+I=K$ leads to~\cite{Krylef89,Nasell1996,anddje98,Nasell1999}
\begin{equation}
\dot{I}=\frac{\beta}{K}(K-I)I-\gamma I -\mu I.\label{eq:sis}
\end{equation}

Steady states of~(\ref{eq:sis}) are $I=0$, and $I=K\left
  (1-\frac{1}{R_0}\right )$, where $R_0 =\frac{\beta}{\mu +\gamma}$ is the
reproductive number. The endemic state is stable for $R_0 >1$, as
seen in Fig.~\ref{fig:ext-scenario}(a). Therefore, deterministically, there is
no way for the disease to go extinct as mentioned previously. However, as
shown in Fig.~\ref{fig:SISa}, the internal noise can in fact induce a large fluctuation
which brings the population into the extinct state. Employing the theory
described in Secs.~\ref{s:mastereq} and~\ref{s:MTEn}, one can find the optimal path to extinction which connects the deterministic
endemic state to a new fluctuational extinct state.

We constrain the population size so that $S=K-I$. This allows one to consider the dynamics of the
constrained SIS model in terms of only the infectious individuals $I$. Rescaling
time by $\mu +\gamma$, the
mean-field equation of the constrained SIS model becomes
\be
\dot{I}=\frac{R_0}{K}(K-I)I- I.
\ee
In rescaling time by $\mu +\gamma$,
the corresponding
stochastic population model is represented by the transition processes of
infection, and removal due to recovery~\cite{Doering2005,dyscla08,sbdl09}. The associated rates $W(X;r)$ are given
as\\

\begin{center}
\begin{tabular}{ll}
Infection: & \hspace{1cm} $W(I;1)=(R_0 /K) I(K-I)$\\
Recovery: & \hspace{1cm} $W(I;-1)= I$\\
&
\end{tabular}
\end{center}
and the master equation is then
\begin{flalign}
\frac{\partial \rho (I,t)}{\partial t}&= [(I+1)\rho (I+1,t)-I\rho
(I,t)]\\
&+\frac{R_0}{K}[(I-1)(K-(I-1))\rho(I-1,t)-I(K-I)\rho (I,t)].
\end{flalign}

The scaled
transition rates in~(\ref{e:transrates}) are given as
\begin{equation}
\begin{array}{lll}
w_{+1}(x) = R_0(1-i)i, &~~~~~&
w_{-1}(x) =   i, \\
&&\\
u_{+1}(x) = 0, &~~~~~&
u_{-1}(x) =   0,
\end{array}
\label{e:SIS-rates}
\end{equation}
where $i=I/K$ is the fraction of infectious individuals in the population.

Substitution of~(\ref{e:SIS-rates}) into~(\ref{e:Hamil}) leads to the Hamiltonian
given as
\begin{equation}
\mathcal{H}(I,\lambda)=R_0 (1-i)i(e^\lambda -1) + i(e^{-\lambda} -1).
\label{e:HSIS}
\end{equation}
Solutions to $\mathcal{H}(I,\lambda)=0$ are
\begin{equation}
i=0, \quad \lambda=0, \quad {\rm and} \quad \lambda(i)=-\ln{(R_0 (1-i))}.
\end{equation}
The third solution is $\lambda_{opt}$ and can also be found
using~(\ref{e:gen_popt}). Taking derivatives of~(\ref{e:HSIS}) with respect to
$x$ and $\lambda$ (see~(\ref{e:Hameqnsi})) leads to the following
system of Hamilton's equations:
\begin{subequations}
\label{e:HEsis}
\begin{flalign}
\dot{i}&=\frac{\partial \mathcal{H}}{\partial \lambda}=R_0(1-i)ie^{\lambda} -ie^{-\lambda},\\
\dot{\lambda}&=-\frac{\partial \mathcal{H}}{\partial i}=-R_0
(1-2i)(e^{\lambda} -1) -(e^{-\lambda} -1).
\end{flalign}
\end{subequations}

\begin{figure}[t!]
\begin{center}
\includegraphics[scale=0.5]{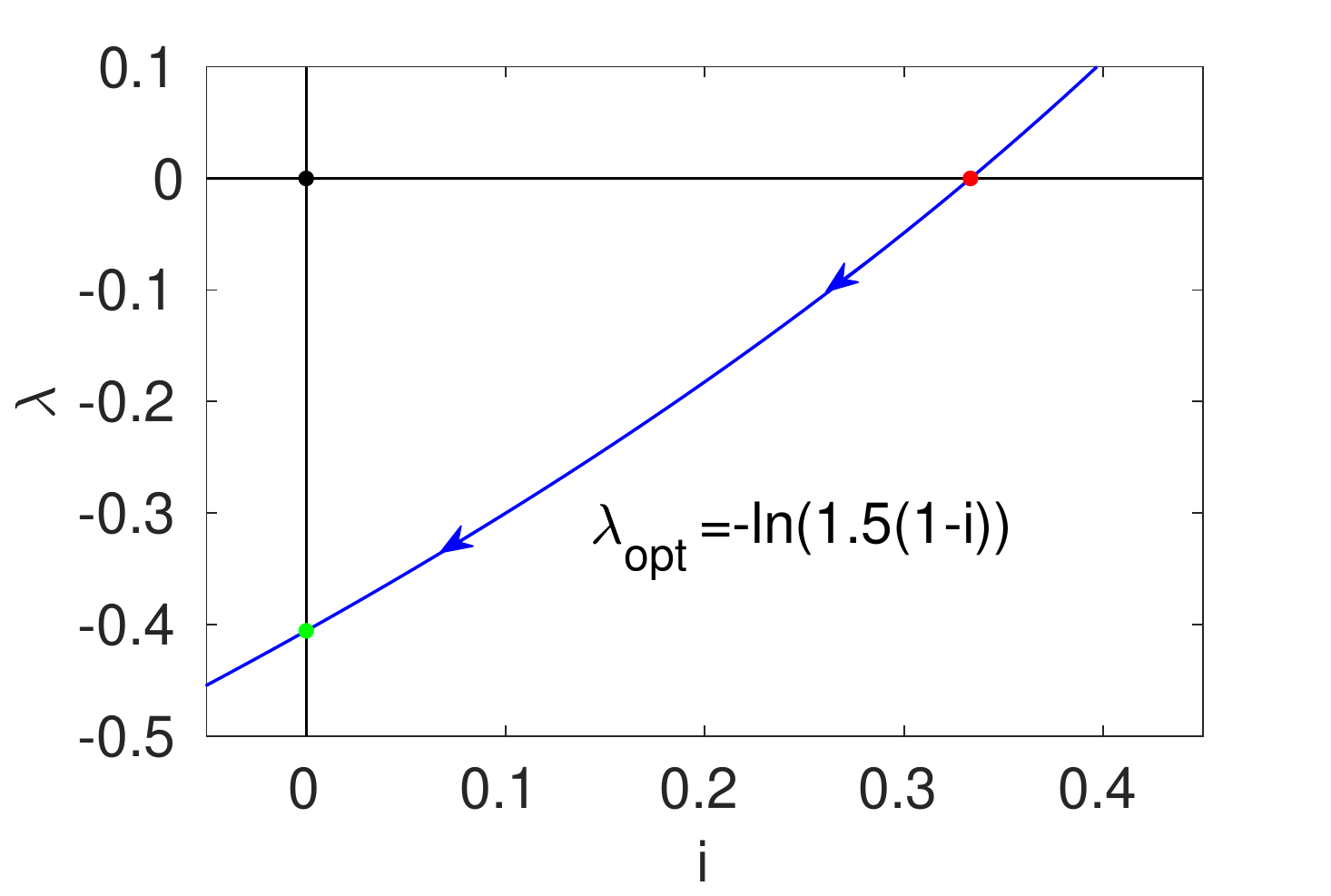}
\end{center}
\caption{\label{fig:1dsis-op} Steady states of Hamilton's equations~(\ref{e:HEsis}) and
  zero-energy curves of the SIS model Hamiltonian~(\ref{e:HSIS}) for $R_0=1.5$. The
  extinct state is located at $(0,0)$ (black dot), the endemic state is
  located at $(1/3,0)$ (red dot), and the fluctuational extinct state is
  located at $(0,-\ln{(1.5)})$ (green dot). The optimal path
  $\lambda_{opt}=-\ln{(1.5(1-i))}$ (blue curve) runs from the endemic state to the
  fluctuational extinct state.}
\end{figure}

This system of Hamilton's equations has three steady states given by
an extinct state $(i,\lambda)=(0,0)$, an endemic state $(i,\lambda)=(1-1/R_0,0)$, and a
fluctuational extinct state $(i,\lambda)=(0,-\ln{(R_0)})$. These steady states along
with the zero energy curves of the Hamiltonian are shown in
Fig.~\ref{fig:1dsis-op}.

The action along the optimal path given by~(\ref{e:Sopt}) is
\begin{equation}
\mathcal{S}_{\rm opt}=\int\limits_{1-\frac{1}{R_0}}^{0} \lambda(i)\, di = \ln{(R_0)}-1+\frac{1}{R_0}.
\end{equation}
The general form of the MTE for single-step problems with a similar topology to that
shown in Fig.~\ref{fig:ext-scenario}(a) is
\begin{equation}\label{e:tau_sis}
\tau = \frac{\sqrt{2 \pi R} \exp \left(\int_{x_0}^{x_1}\left(\frac{u_{+1}(x)}{w_{+1}(x)}-\frac{u_{-1}(x)}{w_{-1}(x)} \right)dx\right)}{(R-1)w_{+1}(x_1)\sqrt{K\lambda'_{opt}(x_1)}} \exp  \left(K  \int_{x_0}^{x_1} \ln \left( \frac{w_{+1}(x)}{w_{-1}(x)} \right) dx \right),
\end{equation}
where $R=w_{+1}^{\prime}(0)/w_{-1}^{\prime}(0)$.
The SIS problem therefore has a MTE given as 
\begin{equation}
\tau=B\exp{(K\mathcal{S}_{\rm opt})}=\frac{R_0}{(R_0
  -1)^2}\sqrt{\frac{2\pi}{N}}\exp{\left [K\left
      (\ln{(R_0)}-1+\frac{1}{R_0}\right )\right ]}.
\end{equation}

\begin{figure}[t!]
\begin{center}
\includegraphics[scale=0.5]{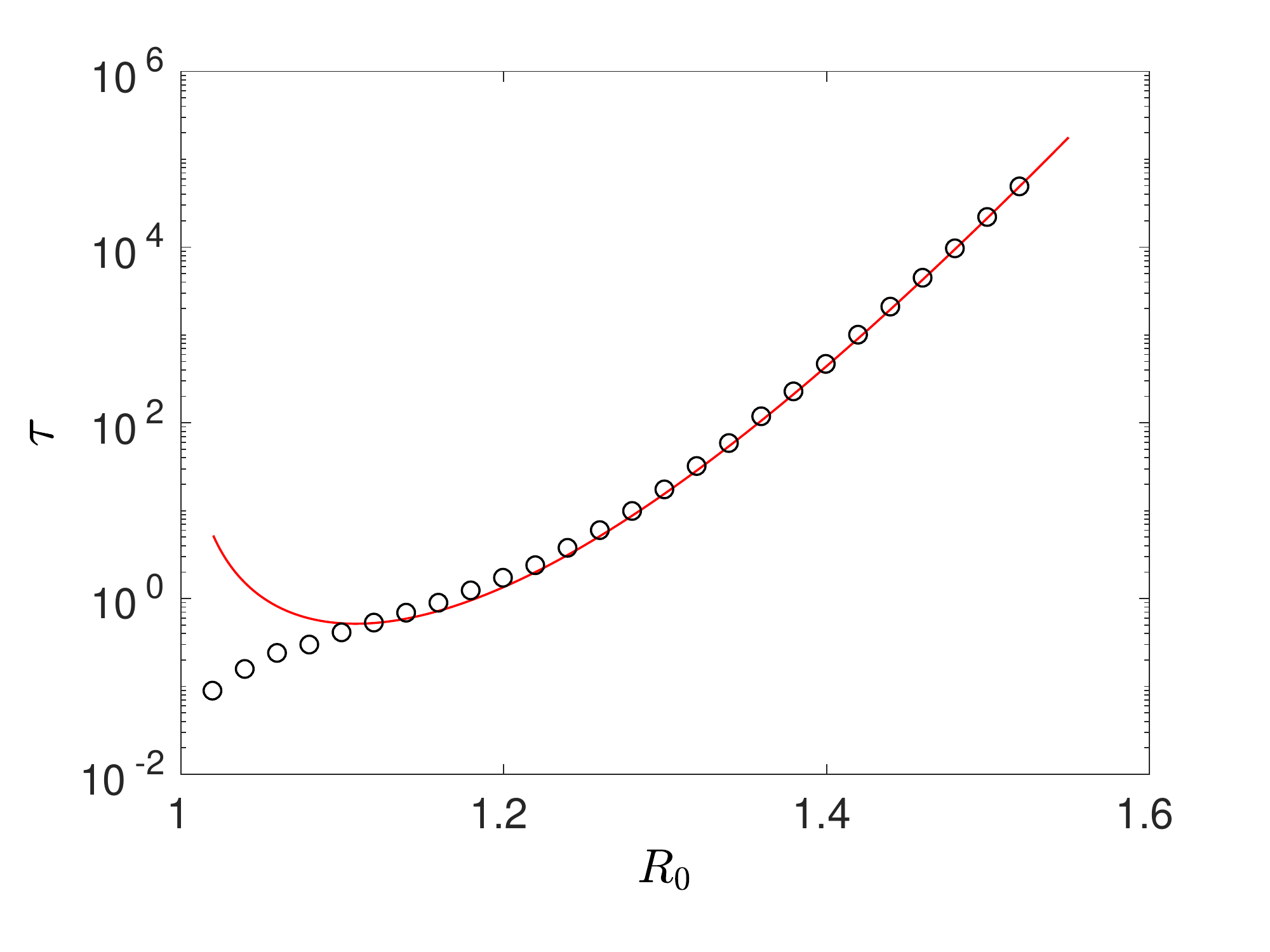}
\end{center}
\caption{\label{fig:1dsis-mte} Comparison of analytical (red curve) and
  numerical (black circles) mean time to extinction versus
  $R_0$ for SIS model.}
\end{figure}

The analytical MTE is confirmed using
numerical simulations. By numerically computing thousands of stochastic realizations and the associated
extinction times, one can calculate the MTE.
Figure~\ref{fig:1dsis-mte} shows a comparison of the analytical and numerical
mean time to extinction as a function of reproductive number $R_0$. There is
excellent agreement except at low $R_0$ where the quasi-stationary assumption
breaks down.

\subsubsection{Allee effect model}
\label{s:Allee}

The Allee effect problem has a mean-field equation given as
\begin{equation}
\dot{x}=\frac{\alpha x^2}{2}-\mu x -\frac{\sigma x^3}{6},\label{e:allee-mf}
\end{equation}
where $\mu$ is the death rate of low-density populations, $\alpha$ is the
growth rate when the population is large enough, and $\sigma$ is a negative
growth rate for an overcrowded population.

The steady states of~(\ref{e:allee-mf}) are $x_0=0$, and
$x_{1,2}=\frac{3\alpha \mp \sqrt{9\alpha ^2 -24\sigma\mu}}{2\sigma}$, where
$x_0$ is stable, $x_1$ is unstable, and $x_2$ is stable as shown in
Fig.~\ref{fig:ext-scenario}(b)~\cite{Assaf2010,Nieddu2014}. Therefore, deterministically, there is
no way for the disease to go extinct, as mentioned previously. However, as
shown in Fig.~\ref{fig:int-switch} the internal noise can in fact induce a large fluctuation
which brings the population to a vicinity of the repelling fixed point of the deterministic rate equation. From there, the population travels essentially deterministically to the extinct state. Employing the theory
described in Secs.~\ref{s:mastereq} and~\ref{s:MTEn}, one can find the optimal path to extinction.

\begin{figure}[t!]
\begin{center}
\includegraphics[scale=0.5]{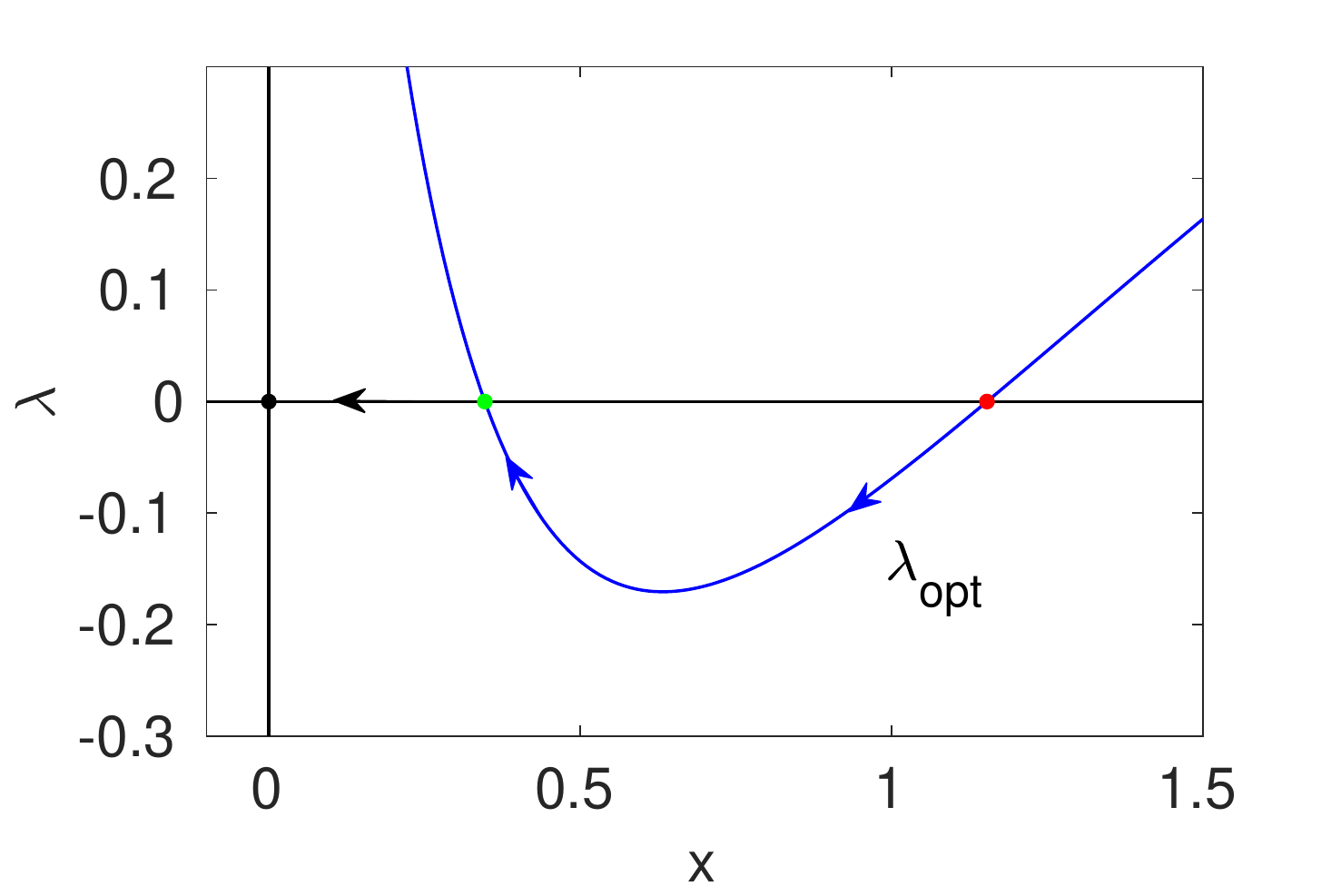}
\end{center}
\caption{\label{fig:allee-opt} Steady states of Hamilton's equations~(\ref{e:HamAll}) and
  zero-energy curves of the Hamiltonian~(\ref{e:Hallee}) for $\mu=0.2$,
  $\alpha =1.5$, and $\sigma=3$. The
  extinct state is located at $(0,0)$ (black dot), the Allee threshold state is
  located at $(0.35,0)$ (green dot), and the carrying capacity state is
  located at $(1.15,0)$ (red dot). The optimal path
  $\lambda_{opt}=-\ln{[(1.2+3x^2)/4.5x]}$ (blue curve) runs from the carrying
  capacity state to the Allee threshold state, and from there continues
  deterministically to the extinct state.}
\end{figure}

The corresponding stochastic population model is represented by the following
transition processes and associated rates $W(X;r)$.

\begin{center}
\begin{tabular}{ll}
&\\
${\displaystyle  X \overset{\mu}{\longrightarrow}\varnothing}$  &
\hspace{1cm} $\mu X$,\\
${\displaystyle 2X \overset{\alpha/K}{\longrightarrow}3X}$  & \hspace{1cm} $\alpha \frac{X(X-1)}{2K}$,\\
${\displaystyle 3X \overset{\sigma/K^2}{\longrightarrow}2X}$  &  \hspace{1cm}
$\sigma \frac{X(X-1)(X-2)}{6K^2}$.\\
&\\
\end{tabular}
\end{center}
The first two transitions involving the death rate $\mu$ and the growth rate
$\alpha$ are required to capture the Allee effect. The
negative growth rate $\sigma$ allows for an overcrowded population to decline
to the carrying capacity $K$~\cite{Assaf2010,Nieddu2014}.

The scaled
transition rates in~(\ref{e:transrates}) are given as
\begin{equation}
\begin{array}{lll}
w_{+1}(x) = \frac{\alpha x^2}{2}, &~~~~~&
w_{-1}(x) =   \mu x + \frac{\sigma x^3}{6}, \\
&&\\
u_{+1}(x) = -\frac{\alpha x}{2}, &~~~~~&
u_{-1}(x) =   -\frac{\sigma x^2}{2}.
\end{array}
\label{e:Allee_wr}
\end{equation}

Substitution of~(\ref{e:Allee_wr}) into~(\ref{e:Hamil}) leads to the Hamiltonian
given as
\begin{equation}
\mathcal{H}(x,\lambda)=\frac{\alpha x^2}{2}(e^{\lambda} -1) +\left (\mu x+\frac{\sigma
    x^3}{6} \right )(e^{-\lambda} -1).
\label{e:Hallee}
\end{equation}
Solutions to $\mathcal{H}(x,\lambda)=0$ are
\begin{equation}
x=0, \quad \lambda=0, \quad {\rm and} \quad \lambda(x)=\ln{\left ( \frac{6\mu +\sigma x^2}{3\alpha x}\right )}.
\end{equation}
The third solution is $\lambda_{opt}$ and can also be found
using~(\ref{e:gen_popt}). Taking derivatives of~(\ref{e:Hallee}) with respect to
$x$ and $\lambda$ (see~(\ref{e:Hameqnsi})) lead to the following
system of Hamilton’s equations:
\begin{flalign}
\label{e:HamAll}
\dot{x}&=\frac{\partial \mathcal{H}}{\partial \lambda}=\frac{\alpha x^2}{2}e^{\lambda}
-\left (\mu
x +\frac{\sigma x^3}{6}\right )e^{-\lambda},\\
&\nonumber\\
\dot{\lambda}&=-\frac{\partial \mathcal{H}}{\partial x}=-\alpha x(e^{\lambda} -1) -\left (\mu
+\frac{\sigma x^2}{2}\right )(e^{-\lambda} -1).
\end{flalign}

This system of Hamilton's equations has three steady states given as
$(x,\lambda)=(x_0,0)=(0,0)$, $(x,\lambda)=(x_1,0)$, and $(x,\lambda)=(x_2,0)$ where $x_0$,
$x_1$, and $x_2$ are the steady states of the mean-field equation~(\ref{e:allee-mf}) provided above.  These steady states along
with the zero energy curves of the Hamiltonian are shown in
Fig.~\ref{fig:allee-opt}.

The general form of the MTE for  single-step problems with a similar topology to that
shown in Fig.~\ref{fig:ext-scenario}(b) is
 \begin{equation}\label{e:tau_ext20}
\tau = \frac{2 \pi \exp \left(\int_{x_1}^{x_2}\left(\frac{u_{+1}(x)}{w_{+1}(x)}-\frac{u_{-1}(x)}{w_{-1}(x)} \right)dx\right)}{w_1(x_2)\sqrt{|\lambda'_{opt}(x_1)|\lambda'_{opt}(x_2)}} \exp  \left(K  \int_{x_2}^{x_1} \ln \left( \frac{w_{-1}(x)}{w_{+1}(x)} \right) dx \right) .
\end{equation}
It is worth noting that the derivation of~\eqref{e:tau_ext20} involves
  matching the solution from $x_2$ to $x_1$ asymptotically with the deterministic solution
  from $x_1$ to $x_0$. Because this latter solution is associated with $\lambda=0$,
  its final form does not involve an integral from $x_1$ to
  $x_0$. Nevertheless, the deterministic contribution is in fact included in~\eqref{e:tau_ext20}. Details regarding the derivation of the prefactor
  $B$ in~(\ref{e:tau_sis}) and~(\ref{e:tau_ext20}) as well as details
  regarding the MTE for multi-step problems can be found in~\cite{Assaf2010}.

\begin{figure}[t!]
\begin{center}
\includegraphics[scale=0.5]{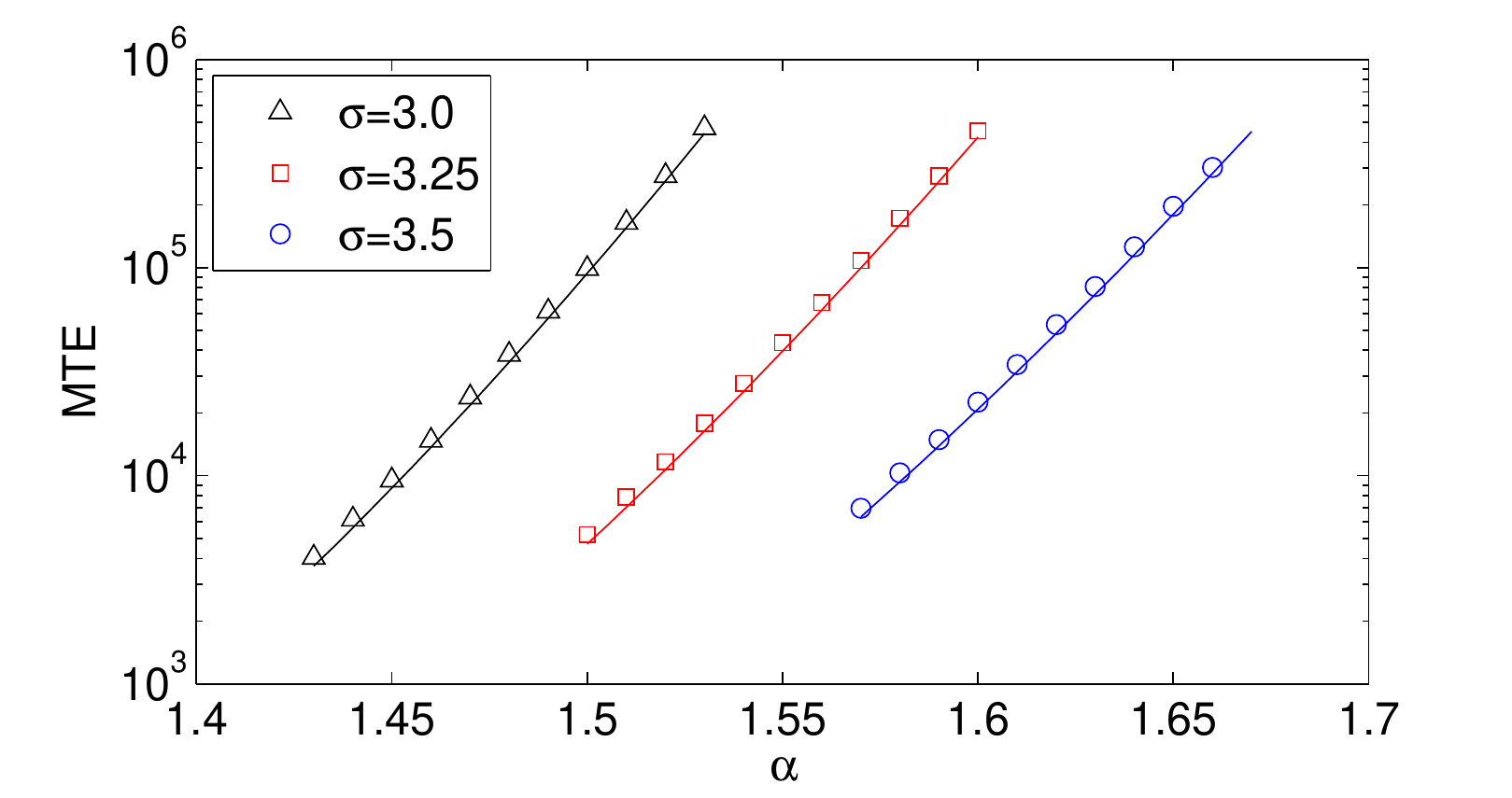}
\end{center}
\caption{\label{fig:allee-mte} Comparison of analytical and numerical mean time to extinction versus
  $\alpha$ for various $\sigma$ for the stochastic Allee effect problem.}
\end{figure}

The analytical MTE for the Allee effect problem is found using~(\ref{e:tau_ext20}) and is confirmed using
numerical simulations. By numerically computing thousands of stochastic realizations and the associated
extinction times, one can calculate the MTE. Figure~\ref{fig:allee-mte} shows the
comparison between the analytical and the numerical mean time to extinction as
a function of $\alpha$ for various choices of $\sigma$~\cite{Nieddu2014}.

\section{Conclusions}

The study of rare events in noise-driven physical and biological models presents a unique opportunity to combine tools from different arenas of applied mathematics, including variational calculus, optimal control, singular perturbation theory, WKBJ theory, asymptotics of integrals, and sampling techniques, to provide considerable insight into ostensibly intractable problems.  The examples included here are just a sample of the important questions that can be answered by an intuitive combination of theory and numerics.  The applicability of these tools can only grow with the explosion of interest in stochastic partial differential equations and geometric growth in computing power.  It is the authors' hope that this primer facilitates the entry of new researchers into this exciting and rapidly growing field.

\appendix
\section{Numerical Methods for Stochastic Differential Equations}
\label{s:appA}

The simplest numerical method to generate a random walk approximating a diffusion process is the stochastic
Euler  method, sometimes referred to as the Euler-Maruyama method. Given the It\^o
stochastic differential equation expressed in~(\ref{e:SDECov}),
the stochastic Euler scheme is given by
\be
X_{n+1}=X_n+f(X_n)\Delta t + \sigma(X_n)\Delta \eta_n,\label{eq:sE1}
\ee
where the time increment is $\Delta t=t_{n+1}-t_n$ and the noise increment is
$\Delta \eta_n = \eta_{t_{n+1}}-\eta_{t_{n}}$.
To implement the stochastic Euler scheme, we note that the noise increments
are independent Gaussian random variables with the following first and second moments:
\be
\mathbb{E}(\Delta \eta_n)=0,\qquad \mathbb{E}((\Delta \eta_n)^2)=D\Delta t,
\ee
so that
\be
\Delta \eta_n \sim {\mathcal N}(0,D\Delta t).\label{eq:Nd1}
\ee
In practice, one can solve~(\ref{eq:sE1}) numerically, using a random number
generator to draw noise values from the distribution given by~(\ref{eq:Nd1}).
Alternatively, one can numerically solve
\be
X_{n+1} = X_n + f(X_n)\Delta t + \sqrt{D\Delta t}\;\sigma(X_n)\;\hat{\eta}_{n},
\ee
where
\be
\hat{\eta}_n \sim {\mathcal N}(0,1).
\ee
For the Ornstein-Uhlenbeck process expressed by $f(x)=-x$ with $D\sigma\equiv 1$,
this becomes
\be
X_{n+1} = (1-\Delta t)X_n + \sqrt{\Delta t}\;\eta_{n},\qquad \eta_n \sim {\mathcal N}(0,1).
\ee

To generate a solution of a stochastic equation, where the noise is internal,
we use the Doob-Gillespie algorithm (also known as the Gillespie algorithm or
Gillespie's stochastic simulation algorithm (SSA))~\cite{doob1942,gillespie_exact_1977}. The algorithm is a type of Monte Carlo
method that was originally proposed by Kendall~\cite{kendall1950} for simulating
birth-death processes and was popularized by Gillespie~\cite{gillespie_exact_1977} as a useful method for simulating
chemical reactions based on molecular collisions. The results of a Gillespie
simulation is a stochastic trajectory that represents an exact sample from the
probability function that solves the master equation. Therefore the method can
be used to simulate population dynamics where molecular collisions are
replaced by individual events and interactions including birth, death, and
infection.

Let ${\bf x}=(x_1,\ldots ,x_n)^T$ denote the state variables of a system,
where $x_i$ provides the number of individuals in
state $x_i$ at time $t$.
The first step of the algorithm is to initialize the number of
individuals in the population compartments ${\bf x}_0$. For a given state
${\bf x}$ of the system, one calculates the transition rates
(birth rate, death rate, contact rate, etc.) denoted as $a_i({\bf x})$ for
$i=1\ldots l$, where $l$ is the number of transitions. Thus the sum of all
transition rates is given by $a_0=\sum\limits_{i=1}^{l}a_i({\bf x})$.

Random numbers are generated to determine both the next event to
occur as well as the time at which
the next event will occur. One simulates the time $\tau$ until the next
transition by drawing from an exponential distribution with mean $1/a_0$. This
is equivalent to drawing a random number $r_1$ uniformly on $(0,1)$ and computing $\tau=(1/a_0)\ln{(1/r_1)}$. During each random time step exactly one event
    occurs. The probability of any particular event taking place is equal to
    its own transition rate divided by the sum of all transition rates
    $a_i({\bf x})/a_0$. A second random number $r_2$ is drawn uniformly on
    $(0,1)$, and it is used to determine the transition event that occurs. If
    $0<r_2 < a_1({\bf x})/a_0$, then the first transition occurs; if
    $a_1({\bf x})/a_0 <r_2 < (a_1({\bf x})+a_2({\bf x}))/a_0$, then the second
    transition occurs, and so on.
Lastly, both the time step and the
number of individuals in each compartment are updated, and the process
is iterated until the disease goes extinct or until the simulation time has
been exceeded.

\section{Numerical Methods for Optimal Path Computations}
\label{s:appB}

Although the two-point boundary value problem given by~(\ref{e:shootForward}) and~(\ref{e:shootBackward}) in Sec.~(\ref{s:escapeTime}) naturally suggests a shooting method for numerical solution, an alternative approach used in the present work is to combine the state and costate equations into a single second-order differential equation,
\begin{align}
\ddot{x}_i &= \partialderiv{f_i}{x_j}\dot{x}_j + (\dot{x}_l-f_l)(a^{-1})_{lj}\partialderiv{a_{ij}}{x_k}\dot{x}_k - a_{ij}\partialderiv{f_k}{x_j}(a^{-1})_{kl}(\dot{x}_l-f_l)\nonumber\\
 &+ \frac12a_{ij}(\dot{x}_k-f_k)\partialderiv{(a^{-1})_{kl}}{x_j}(\dot{x}_l-f_l),\quad i=1,\dots,d,
\label{e:tpbvp}
\end{align}
where $a:=\sigma\sigma^T$ and we have used the summation convention over doubled indices.  The supplementary boundary conditions for $\lambda$ are re-expressed as Robin conditions for $x$ at $t=\tf$.
Introduction of an artificial time $T$ allows a simple initial condition for $x(t,T)$, generally a linear interpolant, to relax to the
optimal path through
\begin{align}
\partialderiv{x_i}{T} = &\partialderiv[2]{{x}_i}{t} - \partialderiv{f_i}{x_j}\partialderiv{{x}_j}t - \left(\partialderiv{{x}_l}t-f_l\right)(a^{-1})_{lj}\partialderiv{a_{ij}}{x_k}\partialderiv{{x}_k}t + a_{ij}\partialderiv{f_k}{x_j}(a^{-1})_{kl}\left(\partialderiv{{x}_l}t-f_l\right)\nonumber\\
 &- \frac12a_{ij}\left(\partialderiv{{x}_k}t-f_k\right)\partialderiv{(a^{-1})_{kl}}{x_j}\left(\partialderiv{{x}_l}t-f_l\right).
\label{e:relax}
\end{align}
In the present work, this PDE was solved using central differences in space and a semi-implicit method in time, with the second derivative handled implicitly and the remaining terms on the right-hand side of~(\ref{e:relax}) handled explicitly.  The PDE was evolved in artificial time $T$ until the maximum of the residual dropped below threshold.

The relaxation described in~\ref{e:relax} occurs over curves parameterized by true time $t$.  This parameterization is not convenient in cases where the true time of travel over the trajectory is infinite, such as for quasi-potential or mean exit time computations, or where the speed is highly nonuniform and the action is therefore concentrated in a particular segment of the path.  It is convenient in such cases to re-parameterize the curve in terms of arclength.  Following Ref.~\cite{heymann_geometric_2008}, we see that the action
\begin{align}
{\mathcal S}_{\tf} &= \frac12\int_{t_0}^{\tf} \|\dot{x}-f(x)\|^2\,dt\\
&=\frac12\int_{t_0}^{\tf}\left((\|\dot{x}\|-\|f\|)^2+2\|\dot{x}\|\|f\|-2(\dot{x},f)\right)\,dt\\
&\geq \int_{t_0}^{\tf}\left(\|\dot{x}\|\|f\| - (\dot{x},f)\right)\,dt\\
&= 2\int_{t_0}^{\tf} \|\dot{x}\|\|f\|\sin^2\frac12\theta(t)\,dt
\end{align}
where $\theta$ is the angle between $\dot{x}$ and $f(x)$.
In the case of infinite-time paths, time can be rescaled to force
\be
\|\dot{x}\|=\|f(x)\|,
\ee
so that the above inequality becomes an equality, and
\be
V(0,y) = 2\inf_{\gamma}\int_{\gamma}\|f(x(\alpha))\|\sin^2\frac12\theta(\alpha)\,d\alpha,
\label{e:arcPar}
\ee
i.e., the time-parameterized minimization can be replaced by a minimization
using arclength as a parameter, where in (\ref{e:arcPar}) the admissible paths
are absolutely continuous.  This approach to
the numerical computation of Wentzell-Freidlin action minimizers is referred
to as the geometric minimum-action method (GMAM), and requires a modification of~\ref{e:relax} to simultaneously compute the time parameterization.  The reader is referred to Ref.~\cite{heymann_geometric_2008} for more details.

A similar method that employs a direct, fully explicit
iterative scheme is the iterative action minimization
method (IAMM)~\cite{LindSch13}.
The IAMM is useful in the general situation where a path connecting steady states $C_a$ and $C_b$
starts at $C_a$ at $t=-\infty$ and ends at $C_b$ at $t=+\infty$. A time parameter $t$ exists such that $-\infty<t<\infty$. For this method, we require a numerical approximation of the time needed to leave the region of $C_a$ and arrive in the region of $C_b$.  Therefore, we define a time $T_{\epsilon}$ such that $-\infty<-T_{\epsilon}<t<T_{\epsilon}<\infty$.  Additionally, $C(-T_{\epsilon}) \approx C_a$ and $C(T_{\epsilon}) \approx C_b$. In other words, the solution stays very near the equilibrium $C_a$ for $-\infty <t\leq -T_{\epsilon}$, has a transition region from $-T_{\epsilon}<t<T_{\epsilon}$, and then stays near $C_b$ for $T_{\epsilon}<t<+\infty$. The interval $[-T_{\epsilon},T_{\epsilon}]$ is discretized into $n$ segments using a uniform step size $h=(2T_{\epsilon})/n$ or a suitable non-uniform step size $h_k$. The corresponding time series is $t_{k+1}=t_k+h_k$.

The derivative of the function value $\mathbf{q}_k$ is approximated using central finite differences by the operator $\delta_h$ given as
\begin{equation}
\frac{d}{dt}\mathbf{q}_k \approx \delta_h\mathbf{q}_k \equiv
\frac{h^2_{k-1}\mathbf{q}_{k+1} + (h^2_k - h^2_{k-1})\mathbf{q}_k -
  h^2_k\mathbf{q}_{k-1}  }{h_{k-1}h^2_k + h_kh^2_{k-1}}, \quad k=0,\ldots ,n.
\label{e:centraldiff}
\end{equation}
Clearly, if a uniform step size is chosen then Eq.~(\ref{e:centraldiff}) simplifies to the familiar form given as
\begin{equation}
\frac{d}{dt}\mathbf{q}_k \approx \delta_h\mathbf{q}_k \equiv
\frac{\mathbf{q}_{k+1} - \mathbf{q}_{k-1}  }{2h}, \quad k=0,\ldots ,n.
\label{e:centraldiff2}
\end{equation}

Thus, one can develop the system of nonlinear algebraic equations
\begin{equation}
\delta_h\mathbf{x}_k - \frac{\partial H(\mathbf{x}_k,\mathbf{p}_k)}{\partial \mathbf{p}}=0, \quad
\delta_h\mathbf{p}_k + \frac{\partial H(\mathbf{x}_k,\mathbf{p}_k)}{\partial
  \mathbf{x}}=0, \quad k=0,\ldots ,n,
\label{e:iammsystem}
\end{equation}
which is solved using a general Newton's method. We let
\begin{equation}
\mathbf{q}_j(\mathbf{x,p})=\{\mathbf{x}_{1,j}...\mathbf{x}_{n,j},\mathbf{p}_{1,j}...\mathbf{p}_{n,j}\}^T
\end{equation}
be an extended vector of $2nN$ components that contains the $j^{\rm th}$
Newton iterate, where $N$ is the number of populations. When $j=0$, $\mathbf{q}_0(\mathbf{x,p})$ provides the initial
``guess'' as to the location of the path that connects $C_a$ and $C_b$. Given the $j^{\rm th}$ Newton iterate $\mathbf{q}_j$, the new $\mathbf{q}_{j+1}$ iterate is found by solving the linear system
\begin{equation}
\mathbf{q}_{j+1}=\mathbf{q}_{j}-\frac{\mathbf{F}\left(\mathbf{q}_{j}\right)}{\mathbf{J}\left(\mathbf{q}_{j}\right)},\label{e:Nit}
\end{equation}
where $\mathbf{F}$ is the function defined by Eq.~(\ref{e:iammsystem}) acting on
$\mathbf{q}_{j}$, and $\mathbf{J}$ is the Jacobian. Equation~(\ref{e:Nit}) is
solved using LU decomposition with partial pivoting.

\section*{Acknowledgments}

We warmly thank Lora Billings for her insight
and Michael Khasin for reading and commenting on an early version of the
manuscript. We also gratefully acknowledge support from the National Science
Foundation awards CMMI-1233397 (EF), DMS-1418956 (EF), and DMS-1109278 (ROM).


\end{document}